\documentclass{article}

\usepackage{amsmath}
\usepackage{amssymb}
\usepackage{amsthm}

\usepackage{shuffle}

\usepackage{bbold}
\usepackage{caption}
\usepackage{hyperref}
\usepackage{mathtools}
\usepackage{float}
\usepackage{comment}
\usepackage{leftidx}
\usepackage{graphicx}
\usepackage{placeins}
\hypersetup{colorlinks=true}
\hypersetup{linkcolor=black,citecolor=black,urlcolor=black,filecolor=black,menucolor=black}

\newcommand{\Int}{\ensuremath{{\textstyle\int}}}
\newcommand{\Der}{\ensuremath{\partial}}
\newcommand{\id}{\ensuremath{{\mathrm{id}}}}
\newcommand{\E}{\ensuremath{{\mathrm{E}}}}

\DeclareMathOperator{\linspan}{span}
\DeclareMathOperator{\coeff}{coeff}

\DeclareMathOperator{\ann}{Ann}
\DeclareMathOperator{\diag}{diag}

\newcommand{\Spol}{\ensuremath{{\mathrm{SP}}}}

\newtheorem{theorem}{\bf Theorem}[section]
\newtheorem{definition}[theorem]{\bf Definition}
\newtheorem{lemma}[theorem]{\bf Lemma}
\newtheorem{corollary}[theorem]{\bf Corollary}

\theoremstyle{remark}
\newtheorem{remark}[theorem]{\bf Remark}
\newtheorem{example}[theorem]{\bf Example}

\usepackage[T1]{fontenc}

\title{The fundamental theorem of calculus in differential rings}
\author{Clemens G.~Raab\textsuperscript{a,}\footnote{Corresponding author}\ \footnote{Present address: Johann Radon Institute of Computational and Applied Mathematics, Austrian Academy of Sciences, Linz, Austria}\ \ and Georg Regensburger\textsuperscript{a,b}}
\hypersetup{pdfauthor={Clemens G. Raab and Georg Regensburger}}
\date{}

\begin{document}

\maketitle

\begin{center}
\textsuperscript{a}Institute for Algebra, Johannes Kepler University Linz, Austria\\
\textsuperscript{b}Institute of Mathematics, University of Kassel, Germany\\
\href{mailto:clemensr@algebra.uni-linz.ac.at}{\tt clemensr@algebra.uni-linz.ac.at}\\
\href{mailto:regensburger@mathematik.uni-kassel.de}{\tt regensburger@mathematik.uni-kassel.de}
\end{center}

\begin{abstract}
In this paper, we study the consequences of the fundamental theorem of calculus from an algebraic point of view.
For functions with singularities, this leads to a generalized notion of evaluation.
We investigate properties of such integro-differential rings and discuss many examples.
We also construct corresponding integro-differential operators and provide normal forms via rewrite rules.
They are then used to derive several identities and properties in a purely algebraic way, generalizing well-known results from analysis.
In identities like shuffle relations for nested integrals and the Taylor formula, additional terms are obtained that take singularities into account.
Another focus lies on treating basics of linear ODEs in this framework of integro-differential operators.
These operators can have matrix coefficients, which allow to treat systems of arbitrary size in a unified way.
In the appendix, using tensor reduction systems, we give the technical details of normal forms and prove them for operators including other functionals besides evaluation.
\end{abstract}


\section{Introduction}

Differential rings are a well-established algebraic structure for modelling differentiation by derivations, i.e.\ linear operations satisfying the Leibniz rule.
More recently, integro-differential algebras have been introduced to additionally model integration and point evaluation of continuous univariate functions by linear operations satisfying corresponding algebraic identities.
In contrast, the integro-differential rings introduced in this paper only use the two identities
\[
 \frac{d}{dx}\int_a^xf(t)\,dt = f(x) \quad\text{and}\quad \int_a^xf^\prime(t)\,dt = f(x)-f(a)
\]
of the fundamental theorem of calculus and the Leibniz rule as axioms.
In particular, this results in a generalized notion of evaluation that is only required to map to constants.
This allows to deal with evaluations of functions even if singularities or discontinuities are present.
For example, it is natural to consider integro-differential rings containing the rational functions leading to so-called hyperlogarithms.

It turns out that many analytic identities and general statements about ODEs, like variation of constants, have a purely algebraic proof in integro-differential rings that is independent of the analytic properties of concrete functions.
Likewise, computations with and transformations of linear integro-differential equations and initial conditions can be done in this algebraic setting as well.
Moreover, exploring algebraic consequences of the Leibniz rule and the fundamental theorem of calculus, we also find new results that introduce additional evaluation terms into identities like shuffle relations for nested integrals and the Taylor formula.
In short, we investigate the analytic operations of differentiation, integration, and evaluation from a purely algebraic point of view.
In this context, we implement in some sense the somewhat provocative statement of Rota \cite[p.~57]{Rota2001} that {\lq\lq}the algebraic structure sooner or later comes to dominate [\ldots]. Algebra dictates the analysis.{\rq\rq}

In order to study linear integro-differential equations and identities, we use the operator point of view, which treats identities of functions as identities of corresponding linear operators acting on them.
To this end, we algebraically construct the ring of integro-differential operators.
For simplifying expressions for operators, we use identities as rewrite rules.
As a key result, we work out a particular rewrite system that can be applied straightforwardly to obtain normal forms and to prove any algebraic identity of integro-differential operators.
Our construction of the ring of operators can be used for operators with scalar coefficients and for operators with matrix coefficients.
In particular, it allows to uniformly deal with scalar equations as well as with systems, even of undetermined size.
Preliminary versions of some results presented in this paper have already been presented by the authors at the conference {\lq\lq}Differential Algebra and Related Topics{\rq\rq} (DART~VII) in 2016.

A related approach was taken in the work of Danuta Przeworska-Rolewicz, see for example \cite{Przeworska-Rolewicz1988}.
In short, she considers a right invertible linear operator on a vector subspace as a generalization of derivation.
Right inverses and projections onto the kernel take the role of integration and evaluation, respectively.
In this linear setting, she develops algebraic generalizations of results for calculus and linear differential equations.
She refers to this as \emph{algebraic analysis}, see also \cite{Przeworska-Rolewicz2000} for historic context and references.
Several results, e.g.\ a Taylor formula, can be formulated already in this purely linear setting.
For other results, multiplication in a commutative algebra is considered.
In some statements, the Leibniz rule or weakened versions of it are needed.
However, she does not consider shuffle relations for nested integrals or additional evaluation terms in the Taylor formula.
Her treatment of operators is limited to properties and identities of given operators and she does not consider rings generated by them or normal forms.

Integro-differential algebras and operators over a field of constants were already introduced in \cite{Rosenkranz2005,RosenkranzRegensburger2008a}, see \cite{RosenkranzRegensburgerTecBuchberger2012} for a detailed overview and further references.
More general differential algebras with integration over rings were introduced in \cite{GuoKeigher2008}, see \cite{GuoRegensburgerRosenkranz2014} for a unified presentation and comparison.
In contrast, the integro-differential rings defined in \cite{HosseinPoorRaabRegensburger2018,HosseinPoor2018} require the integration to be linear over all constants, but the construction of corresponding operators introduced there allows noncommutative coefficients and constants.
In all these cases, integrations satisfy the Rota-Baxter identity \cite{Guo2012}, while integrations considered in this paper do not necessarily satisfy it.
Further references to the literature can be found in the respective sections of the present paper.

So far, all algebraic treatments of integro-differential operators in the literature restrict to multiplicative evaluations, i.e.\ evaluation of a product is the product of the individual evaluations.
Assuming only the Leibniz rule \eqref{eq:LeibnizRule} and the identities \eqref{eq:SectionAxiom} and \eqref{eq:EvaluationDef} of the fundamental theorem of calculus, we deal with non-multiplicative evaluations quite naturally in this paper.
In Section~\ref{sec:IDR}, we introduce (generalized) integro-differential rings following this principle.
Analyzing the relations imposed, we show for example that any linear projection onto constants may be used as the evaluation of such an integro-differential ring.
We also present many examples, both with multiplicative and with non-multiplicative evaluation.
Remark~\ref{rem:IDalgebras} discusses the differences among our definition and the definitions in the literature.

In Section~\ref{sec:products}, we investigate identities satisfied by integrals and their products.
This reveals a generalization of the Rota-Baxter identity for integration that contains an additional evaluation term.
Focusing on nested integrals, we discover generalized shuffle relations with elaborate additional terms involving nested integrals of lower depth.
We also characterize properties of repeated integrals of $1$, which form the smallest integro-differential ring with given ring of constants.

Linear integro-differential operators with coefficients in arbitrary (generalized) integro-differential rings are constructed algebraically in Section~\ref{sec:IDO} by generators and relations.
Working at the operator level enables statements of broader applicability, since operators not only act on integro-differential rings but also on more general modules.
We compute a complete set of rewrite rules, including \eqref{eq:opFTC1} and \eqref{eq:opFTC2} based on the fundamental theorem of calculus, to simplify such operators to normal form.
A precise analysis of the uniqueness of these normal forms is presented in the appendix only, since it requires a refined construction relying on tensor rings (like in \cite{HosseinPoorRaabRegensburger2018,HosseinPoor2018} for the multiplicative case).
After a largely self-contained introduction to tensor reduction systems in the appendix, this is carried out in Section~\ref{sec:TensorIDOFunctionals} allowing also other functionals besides evaluation.
In Section~\ref{sec:IDOaction}, we collect properties of integro-differential operators with coefficients from an integral domain (e.g.\ analytic functions).
In particular, we also characterize the action of such operators when evaluation is multiplicative.

In the remaining sections, we illustrate how computations in the ring of integro-differential operators can be used to prove and generalize well-known results from analysis.
For example, variation of constants remains valid for arbitrary integro-differential rings, which is the focus of Section~\ref{sec:equationalprover}.
In particular, we also detail how integro-differential operators with noncommutative coefficients can be used for proving statements about systems of arbitrary or even undetermined size.
In Section~\ref{sec:generalizations}, we discuss how results from analysis need to be modified for allowing the induced evaluation to be non-multiplicative.
First, we look at the formula for variation of constants, which for multiplicative evaluations automatically satisfies homogeneous initial conditions, and include an extra term to retain this property in general.
Finally, we present a generalization of the Taylor formula with integral remainder term that is valid also for generalized evaluations.

To conclude this introduction, we briefly mention some additional applications of the (generalized) integro-differential rings and operators introduced in this paper.
Based on the generalized shuffle relations presented in Section~\ref{sec:products}, we work out normal forms of nested integrals for symbolic integration and construct the integro-differential closure of commutative differential rings in an upcoming paper \cite{RaabRegensburger}.
Furthermore, the normal forms for integro-differential operators including other functionals, as presented in the appendix, provide a foundation to study systems of linear boundary problems with generalized evaluations, similar to the analysis in \cite{RosenkranzRegensburger2008a} for scalar problems with multiplicative evaluations.

\paragraph{Conventions}
Throughout the paper, rings are implicitly assumed to have a unit element and to be different from the zero ring.
Unless stated otherwise, rings are not assumed to be commutative and can be of arbitrary characteristic.
Nevertheless, for easier reading, we use the notions of \emph{modules} and \emph{linear maps} from the commutative setting to refer to \emph{bimodules} and \emph{bimodule-homomorphisms} over noncommutative rings.
In addition, we use operator notation for linear maps, e.g.\ the Leibniz rule for the derivative of products then reads $\Der{fg}=(\Der{f})g+f\Der{g}$.

\section{Integro-differential rings}
\label{sec:IDR}

To uniformly deal with differentiation of various kinds of functions, we use a few basic abstract notions.
Recall from differential algebra that a \emph{derivation} on a ring $\mathcal{R}$ is an additive map $\Der\colon\mathcal{R}\to\mathcal{R}$ that satisfies the Leibniz rule
\begin{equation}\label{eq:LeibnizRule}
 \Der{fg}=(\Der{f})g+f\Der{g}
\end{equation}
for all $f,g\in\mathcal{R}$.
Then, $(\mathcal{R},\Der)$ is called a \emph{differential ring} and $f\in\mathcal{R}$ is called a \emph{constant} in this differential ring if and only if $\Der{f}=0$.
It is easy to see that the set of constants forms a subring of $\mathcal{R}$ and $\Der$ is linear w.r.t.\ the ring of its constants.
For further theory of differential rings see e.g.\ \cite{Kaplansky1957}.
In addition, we introduce the following notions of integration and evaluation in differential rings.

\begin{definition}\label{def:evaluations}
 Let $(\mathcal{R},\Der)$ be a differential ring and let $\mathcal{C}$ be its ring of constants. We call a $\mathcal{C}$-linear map $\Int\colon \mathcal{R}\to\mathcal{R}$ an \emph{integration} on $\mathcal{R}$, if
 \begin{equation}\label{eq:SectionAxiom}
  \Der\Int{f}=f
 \end{equation}
 holds for all $f \in \mathcal{R}$. A $\mathcal{C}$-linear functional $e\colon \mathcal{R}\to\mathcal{C}$ which acts on $\mathcal{C}$ as the identity is called an \emph{evaluation} on $\mathcal{R}$.
\end{definition}

In other words, integrations on differential rings are right inverses of the derivation that are linear over the constants and evaluations on differential rings are $\mathcal{C}$-linear projectors onto the ring of constants $\mathcal{C}$.

\begin{definition}\label{def:IDRing}
 Let $(\mathcal{R},\Der)$ be a differential ring and let $\Int\colon\mathcal{R}\to\mathcal{R}$ be an integration on $\mathcal{R}$.
 We call $(\mathcal{R},\Der,\Int)$ a \emph{(generalized) integro-differential ring} and we define the \emph{(induced) evaluation} $\E$ on $\mathcal{R}$ by
 \begin{equation}
 \label{eq:EvaluationDef}
  \E{f} := f-\Int\Der{f}.
 \end{equation}
 If in addition $\mathcal{R}$ is a field or skew field, then we also call $(\mathcal{R},\Der,\Int)$ a \emph{(generalized) integro-differential (skew) field}, respectively.
\end{definition}

This extends the definition of integro-differential rings in \cite{HosseinPoorRaabRegensburger2018} by dropping the additional requirement that the induced evaluation should be multiplicative.
In the present paper, the notion of integro-differential rings always refers to Definition~\ref{def:IDRing}.
The following lemma shows that in any integro-differential ring, the (induced) evaluation $\E$ is indeed an evaluation as defined in Definition~\ref{def:evaluations}. Moreover, the ring $\mathcal{R}$ can be decomposed as direct sum of constant and non-constant ``functions''.

\begin{lemma}\label{lem:decompositionofring}
Let $(\mathcal{R},\Der,\Int)$ be an integro-differential ring with constants $\mathcal{C}$.
Then, for all $f \in \mathcal{R}$ and $c \in \mathcal{C}$, we have $\E{f} \in \mathcal{C}$, $\E\Int{f}=0$, and $\E{c}=c$.
Moreover,
\begin{equation*}
\label{eq:decompositionofR}
\mathcal{R} = \mathcal{C} \oplus \Int\mathcal{R}
\end{equation*}
as direct sum of $\mathcal{C}$-modules.
\end{lemma}
\begin{proof}
 First, we compute $\Der\E{f} = \Der(f-\Int\Der{f}) = \Der{f}-\Der{f} = 0$ and $\E\Int{f} = \Int{f}-\Int\Der\Int{f} = 0$ for $f \in \mathcal{R}$ as well as $\E{c} = c-\Int\Der{c} = c$ for $c \in \mathcal{C}$.
 For any $f \in \mathcal{R}$, we have $f = \E{f}+f-\E{f} = \E{f}+\Int\Der{f}$ and hence $\mathcal{R} = \mathcal{C} + \Int\mathcal{R}$. Let $f \in \mathcal{C} \cap \Int\mathcal{R}$ and $g \in \mathcal{R}$ such that $f=\Int{g}$. Then, $0 = \Der{f} = \Der\Int{g} = g$, which implies $f=0$. Hence, the sum $\mathcal{R} = \mathcal{C}+\Int\mathcal{R}$ is direct.
\end{proof}

By the previous lemma, any integration induces an evaluation by $\id-\Int\Der$. Conversely, any evaluation $e$ can be used to define an integration that has $e$ as its induced evaluation, as the following theorem shows.
It is easy to see that two different integrations cannot have the same induced evaluation.
Altogether, on differential rings with $\Der\mathcal{R}=\mathcal{R}$, there is a one-to-one correspondence of integrations and evaluations.

\begin{theorem}
\label{thm:InducedIntegral}
 Let $(\mathcal{R},\Der)$ be a differential ring such that $\Der\mathcal{R}=\mathcal{R}$ and let $e$ be an evaluation on $\mathcal{R}$. Define $\Int\!_e\colon \mathcal{R}\to\mathcal{R}$ by
 \[
  \Int\!_e{f}:=g-eg
 \]
 for all $f \in \mathcal{R}$, where $g \in \mathcal{R}$ is such that $\Der{g}=f$. Then $(\mathcal{R},\Der,\Int\!_e)$ is an integro-differential ring and the induced evaluation is $\E=e$.\par
 Moreover, any integration $\Int$ on $\mathcal{R}$ can be obtained from its induced evaluation $\E$ via this construction: $\Int = \Int\!_{\E}$.
\end{theorem}
\begin{proof}
 Let $\mathcal{C}$ be the ring of constants of $(\mathcal{R},\Der)$.
 First, we show that $\Int\!_e$ is well-defined. If $g,\tilde{g} \in \mathcal{R}$ are such that $\Der{g}=\Der\tilde{g}$, then with $c:=\tilde{g}-g \in \mathcal{C}$ we have $\tilde{g}-e\tilde{g} = g+c-eg-ec = g-eg$, since $ec=c$ by definition of $e$.
 For showing $\mathcal{C}$-linearity of $\Int\!_e$, we let $c \in \mathcal{C}$ and $f_1,f_2,g_1,g_2 \in \mathcal{R}$ with $\Der{g_i}=f_i$. Then, $\Der(cg_1+g_2)=cf_1+f_2$ together with $\mathcal{C}$-linearity of $\id-e$ implies $\mathcal{C}$-linearity of $\Int\!_e$.
 Consequently, $(\mathcal{R},\Der,\Int\!_e)$ is an integro-differential ring, since we also have $\Der\Int\!_e{f}=f$ by construction. The induced evaluation is given by $\E f = f-\Int\!_e\Der f = f-(f-ef) = ef$ for $f \in \mathcal{R}$.\par
 Let $\Int:\mathcal{R}\to\mathcal{R}$ be any $\mathcal{C}$-linear right inverse of $\Der$, $\E$ its induced evaluation, and $f \in \mathcal{R}$. Then, we have $\Int\!_{\E}f = \Int{f}-\E\Int{f} = \Int{f}$ by definition of $\Int\!_{\E}$ and by Lemma~\ref{lem:decompositionofring}.
\end{proof}

In particular, if $(\mathcal{R},\Der,\Int)$ is an integro-differential ring and $e$ is any evaluation on $\mathcal{R}$, then the integration that induces $e$ can be given in terms of $\Int$ by
\begin{equation}\label{eq:IntByEval}
 \Int\!_e:=\Int-e\Int.
\end{equation}
This implies that the difference of two integrations $\Int\!_1,\Int\!_2$ on the same differential ring can be given as $\Int\!_1-\Int\!_2 = \E_2\Int\!_1 = -\E_1\Int\!_2$ in terms of the induced evaluations $\E_1,\E_2$.
More generally, if $(\mathcal{R},\Der,\Int)$ is an integro-differential ring with constants $\mathcal{C}$ and $e:\mathcal{R}\to\mathcal{C}$ is only $\mathcal{C}$-linear, then $\Int\!_e$ defined by \eqref{eq:IntByEval} can be easily seen to be an integration on $\mathcal{R}$ and its induced evaluation is given by $e+(\id-e)\E$, which agrees with $e$ if and only if the latter is an evaluation (i.e.\ $e1=1$).
\par
By Lemma~\ref{lem:decompositionofring}, $\Int\mathcal{R}$ is a direct complement of $\mathcal{C}$ in an integro-differential ring $\mathcal{R}$. Conversely, any direct complement of $\mathcal{C}$ gives rise to an evaluation on $\mathcal{R}$, which in turn induces an integration by Theorem~\ref{thm:InducedIntegral}. More precisely, we have the following characterization of integro-differential rings.

\begin{corollary}\label{cor:characterizeIDRings}
 Let $(\mathcal{R},\Der)$ be a differential ring with ring of constants $\mathcal{C}$.
 Then, $(\mathcal{R},\Der)$ can be enriched into an integro-differential ring if and only if $\Der\mathcal{R}=\mathcal{R}$ and $\mathcal{C}$ is a complemented $\mathcal{C}$-module in $\mathcal{R}$.
 Moreover, if $\Der\mathcal{R}=\mathcal{R}$, there exists a one-to-one correspondence between direct complements of $\mathcal{C}$ in $\mathcal{R}$ and integrations on $(\mathcal{R},\Der)$.
\end{corollary}

This characterization shows that on an integro-differential ring $(\mathcal{R},\Der,\Int)$, in general, there are many other integrations that make $\mathcal{R}$ into an integro-differential ring with the same derivation. In contrast, the following characterization shows that, in general, $\Der$ is the only derivation that turns $\mathcal{R}$ into an integro-differential ring with the same integration.

\begin{lemma}
 Let $\mathcal{R}$ be a ring, let $\mathcal{C}$ be a subring of $\mathcal{R}$, and let $\Int:\mathcal{R}\to\mathcal{R}$ be a $\mathcal{C}$-linear map. Then, there exists a derivation $\Der$ on $\mathcal{R}$ such that $(\mathcal{R},\Der,\Int)$ is an integro-differential ring with constants $\mathcal{C}$ if and only if the following conditions hold.
 \begin{enumerate}
  \item $\Int$ is injective.
  \item $\mathcal{R} = \mathcal{C} \oplus \Int\mathcal{R}$
  \item $(\Int{f})\Int{g} - \Int{(\Int{f})g} - \Int{f\Int{g}} \in \mathcal{C}$ for all $f,g\in\mathcal{R}$.
 \end{enumerate}
 Moreover, this derivation is unique if it exists.
\end{lemma}
\begin{proof}
 First, it is easy to see, from injectivity of $\Int$ and by $\mathcal{R} = \mathcal{C} \oplus \Int\mathcal{R}$, that there exists a $\mathcal{C}$-linear map $\Der:\mathcal{R}\to\mathcal{R}$ such that $\ker\Der=\mathcal{C}$ and $\Der\Int=\id$ and that this map is unique.
 To show that $\Der$ is indeed a derivation, we verify the Leibniz rule on two arbitrary elements of $\mathcal{R}$. By $\mathcal{R} = \mathcal{C} \oplus \Int\mathcal{R}$, we write these two elements as $c+\Int{f}$ and $d+\Int{g}$ with $c,d\in\mathcal{C}$ and $f,g\in\mathcal{R}$. Now, we compute
 \begin{multline*}
  \Der(c+\Int{f})(d+\Int{g})-(\Der(c+\Int{f}))(d+\Int{g})-(c+\Int{f})\Der(d+\Int{g})\\
  = \Der(\Int{f})\Int{g}-f\Int{g}-(\Int{f})g\\
  = \Der\big((\Int{f})\Int{g} - \Int{(\Int{f})g} - \Int{f\Int{g}}\big),
 \end{multline*}
 which is zero by $\ker\Der=\mathcal{C}$ and the last assumption on $\Int$.\par
 Conversely, if $(\mathcal{R},\Der,\Int)$ is an integro-differential ring with constants $\mathcal{C}$, then injectivity of $\Int$ follows from the definition \eqref{eq:SectionAxiom} and the other two conditions on $\Int$ follow from Lemma~\ref{lem:decompositionofring} and Theorem~\ref{thm:mRB}.
\end{proof}

It is straightforward to equip the matrix ring over an integro-differential ring with an integro-differential ring structure.
Such noncommutative integro-differential rings are relevant when working with linear systems, see Section~\ref{sec:linearsystems}.

\begin{lemma}\label{lem:matrixring}
 Let $(\mathcal{R},\Der,\Int)$ be an integro-differential ring with constants $\mathcal{C}$ and let $n\ge1$. Then, $(\mathcal{R}^{n\times n},\Der,\Int)$ is an integro-differential ring with constants $\mathcal{C}^{n\times n}$, where operations $\Der,\Int,\E$ act on matrices by applying the corresponding operation entrywise in $\mathcal{R}$.
\end{lemma}
\begin{proof}
 Clearly, $\Der,\Int,\E$ are $\mathcal{C}^{n\times n}$-linear on $\mathcal{R}^{n\times n}$ satisfying $\Der\Int{A}=A$ and $\Int\Der{A}=A-\E{A}$ for all $A \in \mathcal{R}^{n\times n}$. Moreover, it is straightforward to verify that $\Der$ is a derivation on $\mathcal{R}$ with constants $\mathcal{C}^{n\times n}$.
\end{proof}

\begin{example}\label{ex:multiplicative}
Basic examples for commutative integro-differential rings are univariate polynomials $\mathcal{C}[x]$ and formal power series $\mathcal{C}[[x]]$ over a commutative ring $\mathcal{C}$ with $\mathbb{Q}\subseteq\mathcal{C}$.
The derivation is given by $\Der=\frac{d}{dx}$ with ring of constants $\mathcal{C}$ and integration is defined $\mathcal{C}$-linearly by
\[
 \Int x^n=\frac{x^{n+1}}{n+1}
\]
for all $n \in \mathbb{N}$.
The induced evaluation extracts the constant coefficient and corresponds to evaluation at $0$.
If $\mathcal{C} \subseteq \mathbb{C}$, the integration $\Int$ corresponds to integration $\int_0^x$ from $0$.
Also the rings of complex-valued smooth or analytic functions on a (possibly unbounded) interval $I \subseteq \mathbb{R}$ together with derivation $\Der=\frac{d}{dx}$ and integration $\Int=\int_a^x$, for fixed $a \in I$, are integro-differential rings.
Then, the induced evaluation is the evaluation of functions at the point $a$.
In particular, the ring of exponential polynomials on the real line generated by polynomials and exponential functions is closed under differentiation and integration and hence is an integro-differential ring as well.
Algebraically, for any field $\mathcal{C}$ of characteristic zero, we can consider the ring of exponential polynomials $\mathcal{C}[x,e^{\mathcal{C}x}]$, where $e^{cx}e^{dx}=e^{(c+d)x}$ for all $c,d \in \mathcal{C}$ and $e^{0x}=1$, together with derivation $\Der=\frac{d}{dx}$ and with the integration that is induced by evaluation at $0$ based on Theorem~\ref{thm:InducedIntegral}.
\qed
\end{example}

\begin{example}\label{ex:HurwitzSeries}
A basic example of integro-differential rings of arbitrary characteristic are Hurwitz series, which are closely related to formal power series and have been defined in~\cite{Keigher1975,Keigher1997}, with derivation $\Der(a_0,a_1,\dots)=(a_1,a_2,\dots)$ and integration given by $\Int(a_0,a_1,\dots)=(0,a_0,a_1,\dots)$, see also~\cite{KeigherPritchard2000}.
\qed
\end{example}

The examples mentioned so far have the special property that the evaluation of the integro-differential ring is multiplicative, as in the usual definition of integro-differential algebras.
However, for certain differential rings (in particular for differential fields, cf.\ Corollary~5 in \cite{RosenkranzRegensburgerTecBuchberger2012}), it is not possible to define a multiplicative evaluation for the following reason.
\begin{remark}\label{rem:nonmult}
 If the induced evaluation is multiplicative, one can see easily that no element of $\Int\mathcal{R}$ can have a multiplicative inverse.
 Since otherwise we would have $\E{f}\frac{1}{f}=\E1=1$ and $(\E{f})\E\frac{1}{f}=0\E\frac{1}{f}=0$ for such $f\in\Int\mathcal{R}$.
\qed
\end{remark}
Analytically, if evaluation should correspond to evaluation of functions at a fixed point $a$ for functions that are continuous at $a$, then requiring multiplicativity of evaluation means that functions with poles at $a$ cannot be considered.
In particular, if a function $f$ has a pole of order $m$ at $a$, then evaluation of the product $(x-a)^mf$ gives a nonzero value, but the factor $(x-a)^m$ evaluates to zero at $a$.

In the following theorem, based on results from the literature, we briefly characterize when the induced evaluation is multiplicative.
From \eqref{eq:EvaluationDef} it immediately follows that the identity \eqref{eq:diffRB} is equivalent to multiplicativity of $\E$.
Since in integro-differential rings $\Int$ is $\mathcal{C}$-linear by definition, we have the following characterization of integro-differential rings with multiplicative evaluation.
\begin{theorem}\label{thm:multiplicative}
 Let $(\mathcal{R},\Der,\Int)$ be an integro-differential ring. Then the following properties are equivalent.
 \begin{enumerate}
  \item $\E$ is multiplicative, i.e.\ for all $f,g \in \mathcal{R}$ we have
   \begin{equation}\label{eq:multiplicative}
    \E fg = (\E f)\E g.
   \end{equation}
  \item $\Int$ satisfies the Rota-Baxter identity, i.e.\ for all $f,g \in \mathcal{R}$ we have
   \begin{equation}\label{eq:RB}
    (\Int{f})\Int{g} = \Int{(\Int{f})g} + \Int{f\Int{g}}.
   \end{equation}
  \item The hybrid Rota-Baxter identity holds, i.e.\ for all $f,g \in \mathcal{R}$ we have
   \begin{equation}\label{eq:diffRB}
    (\Int\Der f)\Int\Der g = (\Int\Der f)g + f\Int\Der g - \Int\Der fg.
   \end{equation}
 \end{enumerate}
\end{theorem}
\begin{proof}
 Since $(\mathcal{R},\Der)$ is a differential $\mathbb{Z}$-algebra of weight $0$ and $\Int$ is $\mathbb{Z}$-linear, this immediately follows from items (b), (g), and (a) of Theorem~2.5 in \cite{GuoRegensburgerRosenkranz2014}.
\end{proof}

\begin{remark}\label{rem:IDalgebras}
 In the literature, integro-differential $\mathcal{K}$-algebras (of weight $0$) over a commutative ring $\mathcal{K}$ with unit element are defined as differential $\mathcal{K}$-algebras where the additional map $\Int$ is only required to be $\mathcal{K}$-linear, but has to satisfy the hybrid Rota-Baxter axiom \eqref{eq:diffRB} in addition to \eqref{eq:SectionAxiom}, see \cite{GuoRegensburgerRosenkranz2014}.
 Analogously, the more general notion of differential Rota-Baxter $\mathcal{K}$-algebras (of weight $0$) imposes the Rota-Baxter identity \eqref{eq:RB} instead of \eqref{eq:diffRB} in addition to \eqref{eq:SectionAxiom}.
 On a differential Rota-Baxter $\mathcal{K}$-algebra with constants $\mathcal{C}$, a $\mathcal{K}$-linear map $\E$ is defined by \eqref{eq:EvaluationDef} as well and, by Theorem~2.5 in \cite{GuoRegensburgerRosenkranz2014}, properties \eqref{eq:multiplicative}, \eqref{eq:diffRB}, and $\mathcal{C}$-linearity of $\Int$ are equivalent, see also Proposition~10 in \cite{RosenkranzRegensburgerTecBuchberger2012}.
\par
 By the previous theorem, any (generalized) integro-differential ring with constants $\mathcal{C}$ is an integro-differential $\mathcal{K}$-algebra for $\mathcal{K}=\mathbb{Z}$ and for $\mathcal{K}=\mathcal{C}\cap{Z(\mathcal{R})}$ (i.e.\ $\mathcal{K}=\mathcal{C}$, if $\mathcal{R}$ is commutative) if any of the equivalent conditions holds.
 Conversely, any differential Rota-Baxter $\mathcal{K}$-algebra (of weight $0$) is an integro-differential ring if and only if $\Int$ is linear over the constants $\mathcal{C}$.
 In particular, this is automatically the case for integro-differential $\mathcal{K}$-algebras (of weight $0$) by Proposition~10 in \cite{RosenkranzRegensburgerTecBuchberger2012}.
 For concrete differential Rota-Baxter algebras (of weight $0$) where $\Int$ is not $\mathcal{C}$-linear, see Example~3 in \cite{RosenkranzRegensburger2008a} and the algebraic analog of piecewise functions constructed in \cite{RosenkranzSerwa2019}.
 \qed
\end{remark}

As a basic example for integro-differential rings with non-multiplicative evaluation, we extend the polynomial ring $\mathcal{C}[x]$ over a commutative ring $\mathcal{C}$ with $\mathbb{Q}\subseteq\mathcal{C}$ by adjoining the multiplicative inverse $x^{-1}$.
In order to have a surjective derivation $\Der=\frac{d}{dx}$, we also need to adjoin the logarithm $\ln(x)$ as in the following example.
\begin{example}\label{ex:Laurent}
On $\mathcal{C}[x,x^{-1},\ln(x)]$, with $\mathbb{Q}\subseteq\mathcal{C}$ and $\Der=\frac{d}{dx}$, we can define the $\mathcal{C}$-linear integration recursively as follows.
\[
 \Int{x^k\ln(x)^n} :=
 \begin{cases}\frac{x^{k+1}}{k+1}&k\neq-1\wedge{n=0}\\[\smallskipamount]
 \frac{x^{k+1}}{k+1}\ln(x)^n-\frac{n}{k+1}\Int{x^k\ln(x)^{n-1}}&k\neq-1\wedge{n>0}\\[\smallskipamount]
 \frac{\ln(x)^{n+1}}{n+1}&k=-1\end{cases}
\]
The same recursive definition also works on the larger ring $\mathcal{C}((x))[\ln(x)]$ of formal Laurent series with logarithms, where every element can be written in the form $\sum_{k=-m}^\infty\sum_{n=0}^mc_{k,n}x^k\ln(x)^n$ for some $m\in\mathbb{N}$ and $c_{k,n}\in\mathcal{C}$.
In both cases, the induced evaluation acts by
\[
 \E\sum_{k=-m}^\infty\sum_{n=0}^mc_{k,n}x^k\ln(x)^n = c_{0,0}
\]
and is not multiplicative as expected by Remark~\ref{rem:nonmult}.
Moreover, for the integro-differential subrings of polynomials or formal power series, this evaluation corresponds to the usual multiplicative evaluation at $0$.
\qed
\end{example}

\begin{example}\label{ex:hyperlogs}
Rational functions together with nested integrals of rational functions also form an integro-differential ring with non-multiplicative evaluation.
Algebraically, if $\mathcal{C}$ is a field of characteristic zero, this can be understood as an integro-differential subring of $\mathcal{C}((x))[\ln(x)]$ with integration $\Int$ as in the previous example.
In fact, this ring is the smallest integro-differential ring containing $\mathcal{C}(x)$ and is generated as a $\mathcal{C}(x)$-vector space by $1$ and all nested integrals $\Int{f_1}\Int{f_2}\Int\dots\Int{f_n}$ of arbitrary depth $n\ge1$, where $f_i\in\mathcal{C}(x)$ are proper and have irreducible denominators.
In particular, if $\mathcal{C}=\mathbb{C}$, the integrands can be chosen as $f_i=\frac{1}{x-a_i}$ with $a_i\in\mathcal{C}$.
These kind of nested integrals are called \emph{hyperlogarithms} \cite{LappoDanilevski1928} and have been investigated already in \cite{Kummer1840}.
In \cite{GuoRegensburgerRosenkranz2014}, the free integro-differential algebra (having multiplicative evaluation) generated by $\mathbb{C}(x)$ has been constructed at the somewhat unnatural expense that $\Int1 \not\in \mathbb{C}(x)$ and the constants of the resulting differential ring contain much more than just $\mathbb{C}$.
\qed
\end{example}

\begin{example}\label{ex:Dfinite}
Another example of an integro-differential ring that contains the rational functions are the $D$-finite functions \cite{Stanley1980}. They are characterized as solutions of linear differential equations with rational function coefficients and indeed form a differential ring with surjective derivation $\frac{d}{dx}$.
Since the constants of this differential ring are given by $\mathbb{C}$, there exists an integration by Corollary~\ref{cor:characterizeIDRings}.
\qed
\end{example}

\begin{example}\label{ex:transseries}
All the examples considered above are just rings, not fields.
In contrast, transseries $\mathbb{R}[[[x]]]$ are an explicit construction of a differential field (with field of constants $\mathbb{R}$) that is closed under taking antiderivatives, see~\cite{Hoeven2006,Edgar2010,AschenbrennerDriesHoeven2017} and references therein.
Thus, $\mathbb{R}[[[x]]]$ can be turned into an integro-differential field by Corollary~\ref{cor:characterizeIDRings} whose evaluation necessarily is non-multiplicative by Remark~\ref{rem:nonmult}.
\qed
\end{example}

As shown by Corollary~\ref{cor:characterizeIDRings}, in general, there are many different choices for an integration $\Int$ in order to turn a differential ring with $\Der\mathcal{R}=\mathcal{R}$ into an integro-differential ring. On the same differential ring, for some integrations the induced evaluation is multiplicative and for others it is not.
It may even be the case that a canonical choice of $\Int$ yields a non-multiplicative evaluation while there are other choices that would give a multiplicative evaluation.
In particular, this is the case for exponential polynomials, for example.
They were mentioned above with a multiplicative evaluation, while the following canonical definition of $\Int$ gives rise to a non-multiplicative one.
\begin{example}
The ring of exponential polynomials $\mathcal{C}[x,e^{\mathcal{C}x}]$ over a field $\mathcal{C}$ of characteristic zero is $\mathcal{C}$-linearly generated by terms of the form $x^ke^{cx}$ with $k\in\mathbb{N}$ and $c\in\mathcal{C}$. Apart from the evaluation-based integration on $\mathcal{C}[x,e^{\mathcal{C}x}]$ mentioned above, it is quite natural to define a $\mathcal{C}$-linear integration $\Int$ recursively as follows.
\[
 \Int{x^ke^{cx}} :=
 \begin{cases}\frac{x^{k+1}}{k+1}&c=0\\[\smallskipamount]
 \frac{1}{c}e^{cx}&k=0\wedge{c\neq0}\\[\smallskipamount]
 \frac{1}{c}x^ke^{cx}-\frac{k}{c}\Int{x^{k-1}e^{cx}}&k>0\wedge{c\neq0}\end{cases}
\]
In terms of the Pochhammer symbol $(a)_k:=a{\cdot}(a+1){\cdot}\dots{\cdot}(a+k-1)$, $\Int$ can be given explicitly as
\[
 \Int{x^ke^{cx}}=\sum_{i=0}^k(-k)_{k-i}c^{i-k-1}x^ie^{cx}.
\]
Then, the induced evaluation $\E{f}:=f-\Int\Der{f}$ acts by
\[
 \E x^ke^{cx} =
 \begin{cases}1&k=c=0\\
 0&\text{otherwise}\end{cases}
\]
and is not multiplicative since we have $\E{e^{cx}e^{-cx}}=\E1=1$ but $\E{e^{\pm cx}}=0$ for any $c\neq0$, for example.
On the other hand, $\mathcal{C}[x,e^x]$ is an integro-differential subring with multiplicative evaluation, for example.
With this integration, for instance, the subset $\mathcal{C}[x,e^x]e^x$ is closed under addition, multiplication, derivation, and integration and, hence, could be viewed as an integro-differential subring without unit element and having multiplicative evaluation and the zero ring as its constants.
\qed
\end{example}

All the examples with explicit integration discussed so far contain an integro-differential subring on which the induced evaluation is multiplicative.
In general, however, this need not be the case as the following example shows.
\begin{example}
On $\mathcal{C}[x]$ with the usual derivation and $\mathbb{Q}\subseteq\mathcal{C}$, for example, we can define a $\mathcal{C}$-linear integration by $\Int x^n = \frac{x^{n+1}}{n+1}+c$ for all $n\in\mathbb{N}$ for any fixed $c \in Z(\mathcal{C})$.
Such an integration induces the evaluation $\E{f}=f(0)-cf^\prime(1)$ on $\mathcal{C}[x]$, which is not multiplicative if $c\neq0$ (e.g.\ $f=x$ and $g=x^2-2x$ yield $\E{g}=0$ and $\E{fg}=c$).
\qed
\end{example}

\section{Products of nested integrals}
\label{sec:products}

In integro-differential rings with multiplicative evaluation the standard Rota-Baxter identity \eqref{eq:RB} allows to write the product of integrals as a sum of two nested integrals.
For nested integrals, this leads to shuffle identities \cite{Ree1958} where a product of two nested integrals  is expressed as a sum of nested integrals.
More generally, for Rota-Baxter operators with weight and corresponding shuffle products involving additional terms, see \cite{Guo2012} and references therein.
In general integro-differential rings, i.e.\ if $\E$ is not multiplicative, additional terms involving the evaluation arise in the identities \eqref{eq:RB} and \eqref{eq:diffRB} and also in the shuffle identities.
Note that all evaluation terms are evaluations of products of integrals.
Therefore, they vanish if $\E$ is multiplicative, since $\E\Int{f}=0$ for all $f\in\mathcal{R}$.
\begin{theorem}\label{thm:mRB}
 Let $(\mathcal{R},\Der,\Int)$ be an integro-differential ring. Then the \emph{Rota-Baxter identity with evaluation}
 \begin{equation}\label{eq:mRB}
  (\Int{f})\Int{g} = \Int{f}\Int{g} + \Int(\Int{f})g + \E(\Int{f})\Int{g}
 \end{equation}
 holds for all $f,g \in \mathcal{R}$ as well as
 \begin{equation}\label{eq:mdiffRB}
  (\Int\Der f)\Int\Der g = (\Int\Der f)g + f\Int\Der g - \Int\Der fg - \E(\Int\Der f)\Int\Der g.
 \end{equation}
\end{theorem}
\begin{proof}
 Using \eqref{eq:EvaluationDef}, we can effect the decomposition of $\mathcal{R}$ shown in Lemma~\ref{lem:decompositionofring}.
 For $(\Int{f})\Int{g}$, we thereby obtain the decomposition $\Int\Der(\Int{f})\Int{g} + \E(\Int{f})\Int{g}$, which implies \eqref{eq:mRB} by the Leibniz rule for the derivation $\Der$.
 By $\Int\Der f = f - \E f$, $\Int\Der g = g - \E g$, and $\Int\Der fg = fg - \E fg$, one can write \eqref{eq:mdiffRB} in a form that can be easily verified using the fact that $\E$ is an evaluation.
\end{proof}

As a first application, we show that in every integro-differential ring the repeated integrals of $1$ give rise to an integro-differential subring.
It is the smallest integro-differential ring with the same constants.
In characteristic zero, this ring consists of the univariate polynomials with coefficients in the constants and, for nonzero characteristic, it consists of a finite version of Hurwitz series~\cite{Keigher1997}.
\begin{theorem}\label{thm:polynomials}
 Let $(\mathcal{R},\Der,\Int)$ be an integro-differential ring with constants $\mathcal{C}$.
 For all $n \ge 1$, let $x_n:=\Int^n1 \in \mathcal{R}$ and let $x_0:=1$.
 Then, $1,x_1,x_2,\ldots$ commute with all elements of $\mathcal{C}$ and are $\mathcal{C}$-linearly independent.
 The $\mathcal{C}$-module
 \[
  \mathcal{P} := \linspan_\mathcal{C}\{1,x_1,x_2,\ldots\}
 \]
 is an integro-differential subring of $\mathcal{R}$.
 If $\mathbb{Q} \subseteq \mathcal{R}$, then $\mathcal{P}=\mathcal{C}[x_1]$.
\par
 Moreover, $\E$ is multiplicative on $\mathcal{P}$ if and only if $\E{x_mx_n}=0$ for $m,n\ge1$.
 Assuming $\E$ is multiplicative on $\mathcal{P}$, then we have $x_mx_n=\binom{m+n}{m}x_{m+n}$ and, if in addition $\mathbb{Q} \subseteq \mathcal{R}$, $x_n=\frac{1}{n!}x_1^n$ for $m,n \in \mathbb{N}$.
\end{theorem}
\begin{proof}
 Since $\Int$ is $\mathcal{C}$-linear, every element of $\mathcal{C}$ commutes with $x_i$ for every $i \in \mathbb{N}$, even if $\mathcal{R}$ or $\mathcal{C}$ is noncommutative.
 So, every element of $\mathcal{P}$ is of the form $\sum_{i=0}^nc_ix_i$, for some $c_i \in \mathcal{C}$.
 To show $\mathcal{C}$-linear independence of $x_0,x_1,\dots$, let $n\in\mathbb{N}$ be minimal such that there are $c_0,\dots,c_n \in \mathcal{C}$ with $c_n\neq0$ and $\sum_{i=0}^nc_ix_i=0$.
 Then, $\sum_{i=0}^{n-1}c_{i+1}x_i=\Der\sum_{i=0}^nc_ix_i=0$ would imply $n=0$ by minimality of $n$.
 Because this would yield $c_0=0$, we conclude that $x_0,x_1,\dots$ are $\mathcal{C}$-linearly independent.
\par
 Obviously, $\mathcal{P}$ is closed under $\Der$ and $\Int$ since both operations are $\mathcal{C}$-linear with $\Der{x_n} \in \mathcal{P}$ and $\Int{x_n}=x_{n+1}$ for all $n \in \mathbb{N}$. For showing that $\mathcal{P}$ is closed under multiplication, it suffices to show that $x_mx_n \in \mathcal{P}$ for all $m,n\ge1$. We proceed by induction on the sum $n+m$. For $m=n=1$ we have $x_1^2 = 2x_2+\E{x_1^2}$ by \eqref{eq:mRB}. For $m+n>2$ we have
 \[
  x_mx_n = \Int{x_{m-1}x_n}+\Int{x_mx_{n-1}}+\E{x_mx_n}
 \]
 by \eqref{eq:mRB}. By the induction hypothesis, $x_{m-1}x_n$ and $x_mx_{n-1}$ are in $\mathcal{P}$. Hence, the same is also true after applying $\Int$, which completes the induction. Altogether, $\mathcal{P}$ is an integro-differential subring of $\mathcal{R}$.
\par
 For showing $\mathcal{P}=\mathcal{C}[x_1]$, it is sufficient to prove that every $x_n$ is contained in $\mathcal{C}[x_1]$. By \eqref{eq:mRB}, we obtain $x_{n+1} = x_1x_n - \Int{x_{n-1}x_1} - \E{x_1x_n}$ for all $n\ge1$. Therefore, $x_n \in \mathcal{C}[x_1]$ follows by induction, if $\mathcal{C}[x_1]$ is closed under $\Int$. Assuming $\mathbb{Q} \subseteq \mathcal{R}$, we verify that $\Int{x_1^n} - \frac{1}{n+1}x_1^{n+1} \in \mathcal{C}$ for all $n\ge1$ by applying $\Der$ to it, which shows $\Int{x_1^n} \in \mathcal{C}[x_1]$ for all $n\ge1$.
\par
 Moreover, any $x_n$ with $n\ge1$ satisfies $\E{x_n}=0$. So, if $\E$ is multiplicative on $\mathcal{P}$, then trivially $\E{x_mx_n} = (\E{x_m})\E{x_n} = 0$ for $m,n\ge1$. Conversely, since $\mathcal{P}$ is generated by $1,x_1,x_2,\dots$ as a $\mathcal{C}$-module, we know that $\Int\mathcal{P}$ is generated by $x_1,x_2,\dots$ as a $\mathcal{C}$-module. Hence, applying $\E$ to the product of two elements of $\Int\mathcal{P}$ gives $0$, if $\E{x_mx_n}=0$ for all $m,n\ge1$. Altogether, using the decomposition $\mathcal{P}=\mathcal{C}\oplus\Int\mathcal{P}$ given by Lemma~\ref{lem:decompositionofring}, we conclude that $\E$ is multiplicative on $\mathcal{P}$, if $\E{x_mx_n}=0$ for all $m,n\ge1$.
\par
 Now, we assume $\E$ is multiplicative on $\mathcal{P}$ and let $m,n \in \mathbb{N}$. If $m+n\le1$, then $x_mx_n=\binom{m+n}{m}x_{m+n}$ and $x_n=\frac{1}{n!}x_1^n$ hold trivially. For $m+n\ge2$, it follows inductively by \eqref{eq:mRB} that $x_mx_n = \Int\binom{m+n-1}{m-1}x_{m+n-1} + \Int\binom{m+n-1}{m}x_{m+n-1} = \binom{m+n}{m}x_{m+n}$. In particular, for $n\ge2$, $x_1x_{n-1}=nx_n$ implies $x_n=\frac{1}{n!}x_1^n$ inductively, if $\mathbb{Q} \subseteq \mathcal{R}$.
\end{proof}

Even if $\E$ is not multiplicative on $\mathcal{P}$, we can analyze some properties of the sequence of constants $c_{m,n}:=\E{x_mx_n}$ with $n,m\ge1$.
By $\mathcal{C}$-linearity of $\Int$ and $\E$, it trivially follows that all $c_{m,n}$ are in $Z(\mathcal{C})$.
It can be shown that all $x_n$ commute with each other w.r.t.\ multiplication if and only if the sequence $c_{m,n}$ is symmetric.
If $\mathbb{Q} \subseteq \mathcal{R}$, it can be shown by lengthy computation that the constants $c_{m,n}$ are determined by all $c_{1,n}$ via the recursion
\[
 c_{m,n} = \frac{1}{m}\Bigg(\binom{m{+}n{-}1}{m{-}1}c_{1,m+n-1}+\sum_{j=0}^{m-2}\sum_{k=1}^{n-1}\binom{j{+}k}{j}c_{1,j+k}c_{m-j-1,n-k}\Bigg).
\]
To express $x_n$ in terms of powers of $x_1$ in general, for $\mathbb{Q} \subseteq \mathcal{R}$, we obtain the recursion $x_n=\frac{1}{n!}x_1^n-\sum_{i=2}^n\frac{1}{i!}x_{n-i}\E{x_1^i}$ from Theorem~6.2.2 in \cite{Przeworska-Rolewicz1988}.

\subsection{Generalized shuffle relations}

In this section, we let $(\mathcal{R},\Der,\Int)$ be a commutative integro-differential ring.
Iterating the standard Rota-Baxter identity \eqref{eq:RB} leads to shuffle relations for nested integrals expressing a product of two nested integrals of depth $m$ and $n$ as a sum of nested integrals of depth exactly $m+n$.
Also by recursively applying the Rota-Baxter identity with evaluation \eqref{eq:mRB}, products of nested integrals can be rewritten in $\mathcal{R}$ as sums of nested integrals where also terms of lower depth may occur.
For convenient notation of the formulae involved, it is standard to work in tensor products of $\mathcal{R}$ and to use the shuffle product, which we recall in the following (see \cite{Guo2012} for example).

We consider the $\mathcal{C}$-module $\mathcal{C}\langle\mathcal{R}\rangle:=\bigoplus_{n=0}^\infty\mathcal{R}^{\otimes n}$, where the tensor product is taken over $\mathcal{C}$ and the empty tensor is denoted by $\varepsilon$ so that we have $\mathcal{R}^{\otimes0}=\mathcal{C}\varepsilon$.
Let pure tensors $a_1{\otimes}\dots{\otimes}a_n \in \mathcal{C}\langle\mathcal{R}\rangle$ represent nested integrals $\Int{a_1}\Int{a_2}\dots\Int{a_n} \in \mathcal{R}$.
More formally, by $\mathcal{C}$-linearity of $\Int$, we consider the unique $\mathcal{C}$-module homomorphism $\varphi: \mathcal{C}\langle\mathcal{R}\rangle\to\mathcal{R}$ such that
\[
 \varphi(a_1{\otimes}\dots{\otimes}a_n) = \Int{a_1}\Int{a_2}\dots\Int{a_n} \in \mathcal{R}
\]
and $\varphi(\varepsilon)=1 \in \mathcal{R}$. For pure tensors $a \in \mathcal{C}\langle\mathcal{R}\rangle$, we denote shortened versions of them by $a_i^j:=a_i{\otimes}a_{i+1}{\otimes}\dots{\otimes}a_j$, where $a_i^j:=\varepsilon \in \mathcal{R}^{\otimes0}$ if $i=j+1$.
The shuffle product on $\mathcal{C}\langle\mathcal{R}\rangle$ can be recursively defined as follows. For pure tensors $a,b \in \mathcal{C}\langle\mathcal{R}\rangle$ of length $m$ and $n$, respectively, we set
\[
 a\shuffle{b} :=
 \begin{cases}
  a \otimes b & \text{if }m=0 \vee n=0\\
  a_1\otimes(a_2^m\shuffle{b})+b_1\otimes(a\shuffle{b_2^n}) & \text{otherwise}
 \end{cases}
\]
in $\mathcal{R}^{\otimes(m+n)}$.
Extending this definition to $\mathcal{C}\langle\mathcal{R}\rangle$ by $\mathcal{C}$-linearity, the shuffle product turns $\mathcal{C}\langle\mathcal{R}\rangle$ into a commutative $\mathcal{C}$-algebra.

Using the shuffle product for pure tensors, the product of nested integrals in $\mathcal{R}$ can now be represented as sum of nested integrals as follows.
The constant coefficients of nested integrals of lower depth are evaluations of products of integrals.
Consequently, if $\E$ is multiplicative, then we recover the standard shuffle relations \cite{Ree1958} with all these constant coefficients equal zero.
\begin{theorem}\label{thm:GeneralizedShuffle}
 Let $(\mathcal{R},\Der,\Int)$ be a commutative integro-differential ring with constants $\mathcal{C}$.
 Let $f,g \in \mathcal{C}\langle\mathcal{R}\rangle$ be pure tensors of length $m$ and $n$, respectively.
 Then, the product of the nested integrals $\varphi(f)=\Int{f_1}\Int{f_2}\dots\Int{f_m}$ and $\varphi(g)=\Int{g_1}\Int{g_2}\dots\Int{g_n}$ is given by
 \begin{equation}\label{eq:GeneralizedShuffle}
 \varphi(f)\varphi(g) = \varphi(f\shuffle{g})+\sum_{i=0}^{m-1}\sum_{j=0}^{n-1}e(f_{i+1}^m,g_{j+1}^n)\varphi(f_1^i\shuffle{g_1^j}) \in \mathcal{R}
\end{equation}
with constants
$
 e(f_{i+1}^m,g_{j+1}^n):=\E\varphi(f_{i+1}^m)\varphi(g_{j+1}^n) \in \mathcal{C}.
$
\end{theorem}
\begin{proof}
 Without loss of generality, assume $m \le n$.
 We proceed by induction on $m$.
 If $m=0$, then $f=c\varepsilon$ for some $c \in \mathcal{C}$ and the equation \eqref{eq:GeneralizedShuffle} reads $c\varphi(g) = \varphi(c\varepsilon \otimes g)$, which is trivially true since $c\varepsilon \otimes g=cg$ and $\varphi$ is $\mathcal{C}$-linear.
 For $m\ge1$, we proceed by induction on $n$.
 By virtue of \eqref{eq:mRB}, for $n \ge m$, we have $\varphi(f)\varphi(g) = \Int{f_1}\varphi(f_2^m)\varphi(g) + \Int\varphi(f)g_1\varphi(g_2^n) + e(f,g)$.
 The product $\varphi(f_2^m)\varphi(g)$ is covered by the induction hypothesis on $m$ so that we obtain
 \[
  \Int{f_1}\varphi(f_2^m)\varphi(g) = \Int{f_1}\varphi(f_2^m\shuffle{g})+\sum_{i=1}^{m-1}\sum_{j=0}^{n-1}e(f_{i+1}^m,g_{j+1}^n)\Int{f_1}\varphi(f_2^i\shuffle{g_1^j})
 \]
 by \eqref{eq:GeneralizedShuffle}. The product $\varphi(f)\varphi(g_2^n)$ is covered by the induction hypothesis on $n$ (or on $m$, if $n=m$) so that \eqref{eq:GeneralizedShuffle} yields
 \[
  \Int\varphi(f)g_1\varphi(g_2^n) = \Int{g_1}\varphi(f\shuffle{g_2^n})+\sum_{i=0}^{m-1}\sum_{j=1}^{n-1}e(f_{i+1}^m,g_{j+1}^n)\Int{g_1}\varphi(f_1^i\shuffle{g_2^j}).
 \]
 By definition of $\varphi$ and $\shuffle$, we have $\Int{f_1}\varphi(f_2^m\shuffle{g})+\Int{g_1}\varphi(f\shuffle{g_2^n})=\varphi(f\shuffle{g})$ and similarly $\Int{f_1}\varphi(f_2^i\shuffle{g_1^j})+\Int{g_1}\varphi(f_1^i\shuffle{g_2^j})=\varphi(f_1^i\shuffle{g_1^j})$ for $i,j\ge1$ as well as $\Int{f_1}\varphi(f_2^i\shuffle{g_1^0})=\varphi(f_1^m\shuffle{g_1^0})$ and $\Int{g_1}\varphi(f_1^0\shuffle{g_2^j})=\varphi(f_1^0\shuffle{g_1^j})$.
 Altogether, this yields
 \begin{multline*}
  \varphi(f)\varphi(g) = \varphi(f\shuffle{g})+\sum_{i=1}^{m-1}\sum_{j=1}^{n-1}e(f_{i+1}^m,g_{j+1}^n)\varphi(f_1^i\shuffle{g_1^j})\\
   + \sum_{i=1}^{m-1}e(f_{i+1}^m,g_1^n)\varphi(f_1^m\shuffle{g_1^0}) + \sum_{j=1}^{n-1}e(f_{i+1}^m,g_{j+1}^n)\varphi(f_1^0\shuffle{g_1^j}) + e(f,g),
 \end{multline*}
 which proves \eqref{eq:GeneralizedShuffle}.
\end{proof}

\section{Integro-differential operators}
\label{sec:IDO}

In the following, starting from a given integro-differential ring $(\mathcal{R}, \Der, \Int)$, we define the corresponding ring of operators by generators $\Der, \Int, \E$ and relations.
As additive maps on $\mathcal{R}$, any $f \in \mathcal{R}$ acts as multiplication operator $g \mapsto fg$ and satisfies certain identities together with the maps $\Der, \Int, \E$. Those identities of additive maps that correspond to the defining properties of the operations on $\mathcal{R}$ will be used as defining relations for the abstract ring of operators below.

In particular, the Leibniz rule $\Der fg = f\Der g + (\Der f)g$ of the derivation $\Der$ on $\mathcal{R}$ implies the identity $\Der\circ f = f\circ\Der+\Der{f}$ of additive maps for every multiplication operator $f \in \mathcal{R}$. This motivates the identity \eqref{eq:opLeibniz} in the definition below. Similarly, the identities \eqref{eq:SectionAxiom} and \eqref{eq:EvaluationDef} in $\mathcal{R}$ give rise to the identities \eqref{eq:opFTC1} and \eqref{eq:opFTC2}.
Moreover, from $\mathcal{C}$-linearity of the operations $\Der, \Int, \E$ we obtain $\Int cg = c\Int g$ for all $c \in \mathcal{C}$ and $g \in \mathcal{R}$, for example. In addition, we also obtain $\Int f\E g = (\Int f)\E g$ for all $f,g \in \mathcal{R}$, since $\E g \in \mathcal{C}$. Hence, we also impose commutativity of $\Der,\Int,\E$ with elements of $\mathcal{C}$ and the identities \eqref{eq:opLinDer}--\eqref{eq:opLinEval} in the following definition.

\begin{definition}
\label{def:OperatorRing}
 Let $(\mathcal{R}, \Der, \Int)$ be an integro-differential ring and let $\mathcal{C}$ be its ring of constants. We let
 \[
  \mathcal{R}\langle{\Der,\Int,\E}\rangle
 \]
 be the (noncommutative unital) ring extension of $\mathcal{R}$ generated by indeterminates $\Der,\Int,\E$, where $\Der,\Int,\E$ commute with all elements of $\mathcal{C}$ and the following identities hold for all $f \in \mathcal{R}$.
 \begin{align}
  \Der\cdot{f} &= f\cdot\Der+\Der{f} \label{eq:opLeibniz}\\
  \Der\cdot\Int &= 1 \label{eq:opFTC1}\\
  \Int\cdot\Der &= 1-\E \label{eq:opFTC2}\\
  \Der\cdot{f}\cdot\E &= \Der{f}\cdot\E \label{eq:opLinDer}\\
  \Int\cdot{f}\cdot\E &= \Int{f}\cdot\E \label{eq:opLinInt}\\
  \E\cdot{f}\cdot\E &= \E{f}\cdot\E \label{eq:opLinEval}
 \end{align}
 We call $\mathcal{R}\langle{\Der,\Int,\E}\rangle$ the \emph{ring of (generalized) integro-differential operators (IDO)}.
\end{definition}

Note that, for multiplication in $\mathcal{R}\langle{\Der,\Int,\E}\rangle$, we always explicitly write $\cdot$ when one of $\Der,\Int,\E$ is involved. This is necessary in order to distinguish the product $\Der\cdot{f}$ of operators from the multiplication operator $\Der{f}$, for example.
By construction, the ring $\mathcal{R}\langle{\Der,\Int,\E}\rangle$ has a natural action on $\mathcal{R}$, where the elements of $\mathcal{R}$ act as multiplication operators and $\Der, \Int, \E$ act as the corresponding operations. With this action, $\mathcal{R}$ becomes a left $\mathcal{R}\langle{\Der,\Int,\E}\rangle$-module and multiplication in $\mathcal{R}\langle{\Der,\Int,\E}\rangle$ corresponds to composition of additive maps on $\mathcal{R}$.

Since elements of $\mathcal{C}$ always commute with $\Der,\Int,\E$ in $\mathcal{R}\langle{\Der,\Int,\E}\rangle$, it follows that $\mathcal{R}\langle{\Der,\Int,\E}\rangle$ is a $\mathcal{C}$-algebra whenever $\mathcal{R}$ is commutative.
For computing in the ring $\mathcal{R}\langle{\Der,\Int,\E}\rangle$ of IDO, we use the identities \eqref{eq:opLeibniz}--\eqref{eq:opLinEval} as rewrite rules in the following way.
If the left hand side of one of these identities appears in an expression of an operator, we replace it by the right hand side to obtain a new expression for the same operator.

If rewrite rules can be applied to a given expression in different ways, then it may happen that useful consequences of the defining relations \eqref{eq:opLeibniz}--\eqref{eq:opLinEval} are discovered.
A simple instance starts with the expression $\Int\cdot\Der\cdot\Int$, to which we can apply either \eqref{eq:opFTC1} or \eqref{eq:opFTC2} to obtain the expressions $\Int$ and $\Int-\E\cdot\Int$ for the same operator. Hence, by taking their difference, we see that the identity
\begin{equation}\label{eq:opEvalInt}
 \E\cdot\Int=0
\end{equation}
holds in $\mathcal{R}\langle{\Der,\Int,\E}\rangle$.
Moreover, for every $f \in \mathcal{R}$, the expression $\Int\cdot\Der\cdot{f}$ can be rewritten by \eqref{eq:opLeibniz} and by \eqref{eq:opFTC2}. Thereby we obtain the expressions $\Int\cdot{f}\cdot\Der+\Int\cdot\Der{f}$ and $f-\E\cdot{f}$ for the same operator, which implies the identity
\begin{equation}\label{eq:opIBP}
 \Int\cdot{f}\cdot\Der = f-\E\cdot{f}-\Int\cdot\Der{f}.
\end{equation}
By letting both sides of this identity in $\mathcal{R}\langle{\Der,\Int,\E}\rangle$ act on any $g \in \mathcal{R}$, we show that \emph{integration by parts}
\[
 \Int{f}\Der{g} = fg-\E{fg}-\Int(\Der{f})g
\]
holds in $\mathcal{R}$.
Furthermore, by considering also the newly obtained identity \eqref{eq:opIBP} as rewrite rule, for every $f \in \mathcal{R}$, we can rewrite the expression $\Int\cdot{f}\cdot\Der\cdot\Int$ by \eqref{eq:opFTC1} and by \eqref{eq:opIBP}. Substituting $f$ by $\Int f$ in the difference $f\cdot\Int-\E\cdot{f}\cdot\Int-\Int\cdot\Der{f}\cdot\Int-\Int\cdot{f}$ of the results, we obtain the identity
\begin{equation}\label{eq:opRB}
 \Int\cdot{f}\cdot\Int = \Int{f}\cdot\Int-\Int\cdot\Int{f}-\E\cdot\Int{f}\cdot\Int.
\end{equation}
By acting with both sides of this identity on any $g \in \mathcal{R}$, we obtain an alternative proof for the Rota-Baxter identity with evaluation \eqref{eq:mRB}.
Note that, for $f=1$, we also obtain the following identities from \eqref{eq:opLinDer}--\eqref{eq:opLinEval} and \eqref{eq:opRB}.
\begin{gather}
 \Der\cdot\E=0, \quad\quad \Int\cdot\E=\Int1\cdot\E, \quad\quad \E\cdot\E=\E,\label{eq:opSpecial}\\
 \Int\cdot\Int=\Int1\cdot\Int-\Int\cdot\Int1-\E\cdot\Int1\cdot\Int\label{eq:opSpecialRB}
\end{gather}

In Table~\ref{tab:op}, we collect the identities \eqref{eq:opLeibniz}--\eqref{eq:opSpecialRB} as a rewrite system for expressions of operators in the ring of IDO. In fact, we drop \eqref{eq:opLinDer} since it is redundant in the presence of \eqref{eq:opLeibniz} and \eqref{eq:opSpecial}.

\begin{table}[hb]
\begin{center}
\begin{tabular}{|r@{\ }c@{\ }l|r@{\ }c@{\ }l|}
\hline
 $\Der\cdot{f}$&$=$&$f\cdot\Der+\Der{f}$ & $\Int\cdot{f}\cdot\Der$&$=$&$f-\E\cdot{f}-\Int\cdot\Der{f}$\\
 $\Der\cdot\E$&$=$&$0$ & $\Int\cdot{f}\cdot\E$&$=$&$\Int{f}\cdot\E$\\
 $\Der\cdot\Int$&$=$&$1$ & $\Int\cdot{f}\cdot\Int$&$=$&$\Int{f}\cdot\Int-\Int\cdot\Int{f}-\E\cdot\Int{f}\cdot\Int$\\
 $\E\cdot{f}\cdot\E$&$=$&$\E{f}\cdot\E$ & $\Int\cdot\Der$&$=$&$1-\E$\\
 $\E\cdot\E$&$=$&$\E$ & $\Int\cdot\E$&$=$&$\Int1\cdot\E$\\
 $\E\cdot\Int$&$=$&$0$ & $\Int\cdot\Int$&$=$&$\Int1\cdot\Int-\Int\cdot\Int1-\E\cdot\Int1\cdot\Int$\\
\hline
\end{tabular}
\caption{Rewrite rules for operator expressions}
\label{tab:op}
\end{center}
\vspace*{-\bigskipamount}
\end{table}

\begin{theorem}\label{thm:IDO}
 Let $(\mathcal{R},\Der,\Int)$ be an integro-differential ring. Then, by repeatedly applying the rewrite rules of Table~\ref{tab:op} in any order, every element of the ring $\mathcal{R}\langle{\Der,\Int,\E}\rangle$ can be written as a sum of expressions of the form
 \[
 f \cdot \Der^j, \quad f \cdot \Int \cdot g, \quad f \cdot \E \cdot g \cdot \Der^j, \quad \text{or} \quad f \cdot \E \cdot h \cdot \Int \cdot g
 \]
 where $j\in\mathbb{N}_{0}$, $f,g\in\mathcal{R}$, and $h\in\Int\mathcal{R}$.
\end{theorem}

Note that the expressions specified in the above theorem are irreducible in the sense that they cannot be rewritten any further by any rewrite rules from Table~\ref{tab:op}.
The above derivation of identities \eqref{eq:opEvalInt}--\eqref{eq:opSpecialRB} is similar to Knuth-Bendix completion \cite{KnuthBendix1970} and Buchberger's algorithm for computing Gr\"obner bases \cite{Buchberger1965,Mora1994}.
So, one can show that Table~\ref{tab:op} represents all consequences of Definition~\ref{def:OperatorRing} in the sense that every identity in $\mathcal{R}\langle{\Der,\Int,\E}\rangle$ can be proven by applying the rewrite rules in the table and by exploiting identities in $\mathcal{R}$.
Moreover, the irreducible forms of operators specified in the above theorem are unique up to multiadditivity and commutativity.

In the appendix, we will give a precise statement (Theorem~\ref{thm:IDOtensor}) of this by giving an explicit construction of the ring $\mathcal{R}\langle{\Der,\Int,\E}\rangle$ as a quotient of an appropriate tensor ring.
Then, the translation of Table~\ref{tab:op} into a tensor reduction system facilitates the proof.
In fact, this proof is carried out in Theorem~\ref{thm:IDOPhi} for a more general class of operators including linear functionals, which are useful for dealing with boundary problems, for instance.

\begin{remark}
In the literature, integro-differential operators were considered only with multiplicative evaluation so far.
Integro-differential operators were first introduced in \cite{Rosenkranz2005,RosenkranzRegensburger2008a} over a field of constants using a parametrized Gr\"obner basis in infinitely many variables and a basis of the commutative coefficient algebra.
Integro-differential operators with polynomial coefficients over a field of characteristic zero were also studied using generalized Weyl algebras \cite{Bavula2011}, skew polynomials \cite{RegensburgerRosenkranzMiddeke2009}, and noncommutative Gr\"obner bases \cite{QuadratRegensburger2020}.
A general construction of rings of linear operators over commutative operated algebras is presented in \cite{GaoGuoRosenkranz2018}.
In particular, integro-differential operators are discussed in that setting and also differential Rota-Baxter operators are investigated.
Tensor reduction systems have already been used in \cite{HosseinPoorRaabRegensburger2018,HosseinPoor2018} for the construction of IDO including functionals in the special case that the evaluation and all functionals are multiplicative.
There, also additional operators arising from linear substitutions were included.
In particular, these cover integro-differential-time-delay operators, which were already constructed algebraically in \cite{Quadrat2015}, see also \cite{CluzeauHosseinPoorQuadratRaabRegensburger2018}.
\qed
\end{remark}

To refer to integro-differential operators of special form, we use the following notions.
\begin{definition}
 Let $L \in \mathcal{R}\langle{\Der,\Int,\E}\rangle$, then we call $L$
 \begin{enumerate}
  \item a \emph{differential operator}, if there are $f_0,\dots,f_n \in \mathcal{R}$ such that
   \[L=\sum_{i=0}^nf_i\cdot\Der^i,\]
   where we call $L$ \emph{monic}, if $f_n=1$,
  \item an \emph{integral operator}, if there are $f_1,\dots,f_n,g_1,\dots,g_n \in \mathcal{R}$ such that
   \[L=\sum_{i=1}^nf_i\cdot\Int\cdot{g_i},\]
  \item an \emph{initial operator}, if there are $f_1,\dots,f_n \in \mathcal{R}$ and differential and integral operators $L_1,\dots,L_n$ such that
   \[L=\sum_{i=1}^nf_i\cdot\E\cdot{L_i}.\]
   In particular, we call an initial operator $L$ \emph{monic}, if there is a differential operator $L_1$ and an integral operator $L_2$ such that
   \[L=\E\cdot(L_1+L_2).\]
 \end{enumerate}
\end{definition}
In $\mathcal{R}\langle{\Der,\Int,\E}\rangle$, using Table~\ref{tab:op}, one can check that the differential operators are the elements of the subring generated by $\mathcal{R}$ and $\Der$, the integral operators are the elements of the $\mathcal{R}$-bimodule generated by $\Int$, the initial operators are the elements of the two-sided ideal generated by $\E$, and the monic initial operators are the elements of the right ideal generated by $\E$.
Theorem~\ref{thm:IDO} says that every integro-differential operator can be written as the sum of a differential operator, an integral operator, and an initial operator.
In fact, by the stronger results in the appendix, this decomposition of integro-differential operators even is unique.

\begin{remark}
\label{rem:multiplicative}
 We outline how the construction of integro-differential operators changes when the evaluation is multiplicative, see also \cite{RosenkranzRegensburgerTecBuchberger2012,HosseinPoorRaabRegensburger2018,HosseinPoor2018}.
 For multiplicative evaluation, we have $\E{fg}=(\E{f})\E{g}$ for all $f,g \in \mathcal{R}$.
 So, for the algebraic construction of corresponding operators as in Definition~\ref{def:OperatorRing}, we need to impose in addition that
 \begin{equation}
 \label{eq:opMulEval}
  \E\cdot{f}=\E{f}\cdot\E
 \end{equation}
 for all $f\in\mathcal{R}$.
 This does not give rise to new consequences other than \eqref{eq:opEvalInt}--\eqref{eq:opSpecialRB}, but together with \eqref{eq:opSpecial} it makes \eqref{eq:opLinEval} redundant.
 Hence, in Table~\ref{tab:op}, we can replace \eqref{eq:opLinEval} by \eqref{eq:opMulEval}.
 Moreover, \eqref{eq:opMulEval} allows to reduce the evaluation term in \eqref{eq:opIBP} and, since $\E\Int{f}=0$ for all $f\in\mathcal{R}$, to omit the evaluation terms in \eqref{eq:opRB} and \eqref{eq:opSpecialRB}.
 Consequently, the irreducible forms of operators given in Theorem~\ref{thm:IDO} can be simplified to
 \[
  f \cdot \Der^j, \quad f \cdot \Int \cdot g, \quad f \cdot \E \cdot \Der^j,
 \]
 where $j\in\mathbb{N}_{0}$ and $f,g\in\mathcal{R}$.
 Evidently, this ring of operators is isomorphic to $\mathcal{R}\langle{\Der,\Int,\E}\rangle$ factored by the two-sided ideal $(\E\cdot{f}-\E{f}\cdot\E\ |\ f\in\mathcal{R})$.
 \qed
\end{remark}

\subsection{IDO over integral domains}
\label{sec:IDOaction}

In many concrete situations, when computing with differential operators or differential equations with scalar coefficients, the order of a product of differential operators is the sum of the orders of the factors.
This is equivalent to the ring $\mathcal{R}$ of coefficients being an integral domain, i.e.\ a commutative ring without nontrivial zero divisors.
For the rest of this section, we only consider integro-differential rings that are integral domains.
Under this assumption, we can investigate further properties of computations in the ring of IDO.
For instance, an integro-differential equation can be reduced to a differential equation by differentiation, see Lemma~\ref{lem:OperatorIdealsLeft} below.
Additionally, investigating the action of IDO on integro-differential rings or $\mathcal{R}\langle{\Der,\Int,\E}\rangle$-modules algebraically leads to analyzing two-sided ideals.

As explained earlier, elements of $\mathcal{R}\langle{\Der,\Int,\E}\rangle$ naturally act as additive maps on $\mathcal{R}$.
Likewise, they act naturally on any integro-differential ring extension of $\mathcal{R}$.
In Theorem~\ref{thm:FaithfulAction}, we also provide conditions when this action is faithful, i.e.\ $0 \in \mathcal{R}\langle{\Der,\Int,\E}\rangle$ is the only element that induces the zero map.

Let $(\mathcal{S},\Int,\Der)$ be the integro-differential ring defined in Example~\ref{ex:Laurent} and assume $\mathcal{C}$ is an integral domain.
Then, $\mathcal{S}\langle{\Der,\Int,\E}\rangle$ does not act faithfully on $\mathcal{S}$, since $\E\cdot\ln(x)$ acts like zero.
However, with the integro-differential subring $\mathcal{R}:=\mathcal{C}[x]$ of $\mathcal{S}$, $\mathcal{R}\langle{\Der,\Int,\E}\rangle$ acts faithfully on $\mathcal{S}$ by Theorem~\ref{thm:FaithfulAction} below.
To see this, let $L:=\sum_{i=0}^nf_i\cdot\Der^i \in \mathcal{R}\langle{\Der,\Int,\E}\rangle$ and let $k\in\mathbb{N}$ such that $\coeff(f_n,x^k)\neq0$.
Applying $\E\cdot{L}$ to $x^{n-k}\ln(x)^n \in \mathcal{S}$, we obtain $(\E\cdot{L})x^{n-k}\ln(x)^n = n!\coeff(f_n,x^k) \neq 0$.

As a preparatory step, we need some basic statements involving multiplication with constants that are valid for integro-differential rings that are integral domains.
Linear independence over constants is tied to the Wronskian.
Following the analytic definition, the Wronskian of elements $f_1,\dots,f_n$ of a commutative differential ring is defined by
\[
 W(f_1,\dots,f_n):=\det\left(\left(\Der^{i-1}f_j\right)_{i,j=1,\dots,n}\right).
\]

\begin{lemma}\label{lem:DifferentialIntegralDomain}
 Let $(\mathcal{R},\Der)$ be a differential ring that is an integral domain with ring of constants $\mathcal{C}$ and let $f_1,\dots,f_n \in \mathcal{R}$, $n\ge1$.
 \begin{enumerate}
  \item $f_1,\dots,f_n$ are linearly independent over $\mathcal{C}$ if and only if their Wronskian $W(f_1,\dots,f_n)$ is zero.
  \item If $f_1,\dots,f_n$ are linearly independent over $\mathcal{C}$, then $g_1,\dots,g_{n-1}$ are linearly independent over $\mathcal{C}$, where $g_i:=W(f_i,f_n)$.
 \end{enumerate}
\end{lemma}
\begin{proof}
 For showing the first statement, we note that neither linear independence nor zeroness of the Wronskian changes, if we replace $\mathcal{R}$ and hence $\mathcal{C}$ by their quotient fields.
 For differential fields, a proof can be found in \cite{Kaplansky1957}, which implies the statement given here.
\par
 For showing the second statement, we assume that $f_1,\dots,f_n$ are linearly independent over $\mathcal{C}$ and we let $c_1,\dots,c_{n-1} \in \mathcal{C}$ such that $\sum_{i=1}^{n-1}c_ig_i=0$.
 By definition of $g_i$ and multilinearity of the Wronskian over $\mathcal{C}$, we conclude $W(\sum_{i=1}^{n-1}c_if_i,f_n)=0$.
 By assumption on $f_1,\dots,f_n$, this implies $\sum_{i=1}^{n-1}c_if_i=0$ and hence $c_1=\ldots=c_{n-1}=0$.
\end{proof}

\begin{lemma}\label{lem:ConstantMultiples}
 Let $(\mathcal{R},\Der,\Int)$ be an integro-differential ring that is an integral domain and let $\mathcal{C}$ be its ring of constants.
 Then, elements of $\mathcal{C}$ commute with all elements of $\mathcal{R}\langle{\Der,\Int,\E}\rangle$.
 Moreover, for any $L \in \mathcal{R}\langle{\Der,\Int,\E}\rangle$ and any nonzero $c\in\mathcal{C}$, we have that $L=0$ if and only if $c\cdot{L}=0$.
\end{lemma}
\begin{proof}
 By construction, constants commute with $\Der,\Int,\E \in \mathcal{R}\langle{\Der,\Int,\E}\rangle$.
 Since $\mathcal{R}$ is commutative, it follows that elements of $\mathcal{C}$ commute with all elements of $\mathcal{R}\langle{\Der,\Int,\E}\rangle$.
 Therefore, $\mathcal{R}\langle{\Der,\Int,\E}\rangle$ is a unital $\mathcal{C}$-algebra and hence $Q(\mathcal{C})\otimes_{\mathcal{C}}\mathcal{R}\langle{\Der,\Int,\E}\rangle$ is a unital $\mathcal{C}$-algebra as well, with multiplication $(c_1\otimes{L_1})\cdot(c_2\otimes{L_2})=(c_1c_2)\otimes(L_1\cdot{L_2})$.
 Now, $1\otimes{L}=c^{-1}\otimes(c\cdot{L})$ implies $L=0$ if $c\cdot{L}=0$.
\end{proof}

Now, we are ready to have a closer look at certain computations with IDO.
The following two lemmas construct left or right multiples of IDO by differential operators such that the product does not involve the integration operator anymore.
The tricky part will be to ensure that the product is a nonzero operator again.

\begin{lemma}\label{lem:OperatorIdealsLeft}
 Let $(\mathcal{R},\Der,\Int)$ be an integro-differential ring that is an integral domain.
 Let $\mathcal{C}$ be the ring of constants of $\mathcal{R}$ and let $L \in \mathcal{R}\langle{\Der,\Int,\E}\rangle$ be not an initial operator.
 Then, there exist nonzero $h_1,\dots,h_n \in \mathcal{R}$ such that the following product is a nonzero differential operator.
 \[
  (h_1\cdot\Der-\Der{h_1})\cdot\ldots\cdot(h_n\cdot\Der-\Der{h_n})\cdot{L}
 \]
\end{lemma}
\begin{proof}
 There are $n \in \mathbb{N}$, a nonzero $c \in \mathcal{C}$, differential operators $L_0,\dots,L_n \in \mathcal{R}\langle{\Der,\Int,\E}\rangle$, and $f_1,\dots,f_n,g_1,\dots,g_n \in \mathcal{R}$ such that $f_1,\dots,f_n$ are linearly independent over $\mathcal{C}$ and $c\cdot L = L_0+\sum_{i=1}^nf_i\cdot\left(\Int\cdot{g_i}+\E\cdot{L_i}\right)$.
 Since $L$ is not an initial operator, $L_0$ and $g_1,\dots,g_n$ are not all zero by Lemma~\ref{lem:ConstantMultiples}.
 If $L_0=0$, then we assume without loss of generality that $g_1\neq0$.
\par
 To remove the sum $\sum_{i=1}^nf_i\cdot\left(\Int\cdot{g_i}+\E\cdot{L_i}\right)$, we will multiply $c\cdot L$ iteratively by $n$ first-order differential operators from the left as follows.
 If $n>0$, then by \eqref{eq:opLeibniz} and \eqref{eq:opFTC1} we have $(f_n\cdot\Der-\Der{f_n})\cdot{f_i}\cdot\Int\cdot{g_i} = f_nf_ig_i+(f_n\Der{f_i}-(\Der{f_n})f_i)\cdot\Int\cdot{g_i}$ for all $i$.
 Therefore, using \eqref{eq:opLinDer}, we obtain
 \begin{multline*}
  (f_n\cdot\Der-\Der{f_n})\cdot{c}\cdot{L} = (f_n\cdot\Der-\Der{f_n})\cdot{L_0} + \sum_{i=1}^nf_nf_ig_i\\
  + \sum_{i=1}^{n-1}(f_n\Der{f_i}-(\Der{f_n})f_i)\cdot\left(\Int\cdot{g_i} + \E\cdot{L_i}\right),
 \end{multline*}
 which has the form $\tilde{L}_0+\sum_{i=1}^{n-1}\tilde{f}_i\cdot\left(\Int\cdot{g_i} + \E\cdot{L_i}\right)$ similar to $c\cdot L$.
 Note that also the following two properties of the above representation of $c\cdot L$ are preserved.
 First, if $L_0$ is nonzero, the differential operator $\tilde{L}_0:=(f_n\cdot\Der-\Der{f_n})\cdot{L_0} + \sum_{i=1}^nf_nf_ig_i$ is nonzero too, since $f_n\neq0$.
 Second, $\tilde{f}_i:=f_n\Der{f_i}-(\Der{f_n})f_i$, for $i \in \{1,\dots,n-1\}$, are again linearly independent over $\mathcal{C}$ by Lemma~\ref{lem:DifferentialIntegralDomain}.
 Now, we let $h_n:=f_nc$ such that $(h_n\cdot\Der-\Der{h_n})\cdot{L} = (f_n\cdot\Der-\Der{f_n})\cdot{c}\cdot{L}$.
\par
 If $n=1$, we only need to show that $(h_n\cdot\Der-\Der{h_n})\cdot{L} = \tilde{L}_0$ is nonzero.
 As remarked above, if $L_0\neq0$ then $\tilde{L}_0\neq0$.
 On the other hand, if $L_0=0$, then $\tilde{L}_0=f_1^2g_1$, which is nonzero as well due to the assumptions.
\par
 If $n>1$, we continue by repeating what we did with $c\cdot{L}$ above.
 That is, noting $\tilde{f}_1,\dots,\tilde{f}_{n-1}$ are linearly independent over $\mathcal{C}$, we multiply $(h_n\cdot\Der-\Der{h_n})\cdot{L} = \tilde{L}_0+\sum_{i=1}^{n-1}\tilde{f}_i\cdot\left(\Int\cdot{g_i} + \E\cdot{L_i}\right)$ from the left by $h_{n-1}\cdot\Der-\Der{h_{n-1}}$, where $h_{n-1}:=\tilde{f}_{n-1}$.
 We iterate this a total of $n-1$ times, the result is a differential operator.
 To see that it is also nonzero, we focus on the last step, i.e.\ we refer to the case $n=1$ above.
\end{proof}

An immediate consequence of the previous lemma is that a left ideal in $\mathcal{R}\langle{\Der,\Int,\E}\rangle$ is nontrivial if and only if it contains a nonzero differential operator or a nonzero initial operator.
Left ideals in $\mathcal{R}\langle{\Der,\Int,\E}\rangle$ arise, for instance, as sets of operators annihilating a given subset of a left $\mathcal{R}\langle{\Der,\Int,\E}\rangle$-module.
As another consequence, for two-sided ideals, we obtain Lemma~\ref{lem:OperatorIdealsTwoSided} below.

\begin{lemma}\label{lem:OperatorIdealsRight}
 Let $(\mathcal{R},\Der,\Int)$ be an integro-differential ring that is an integral domain.
 Let $\mathcal{C}$ be the ring of constants of $\mathcal{R}$ and let $L \in \mathcal{R}\langle{\Der,\Int,\E}\rangle$ be a nonzero monic initial operator.
 Then, there exist nonzero $h_1,\dots,h_n \in \mathcal{R}$ and a nonzero differential operator $\tilde{L} \in \mathcal{R}\langle{\Der,\Int,\E}\rangle$ such that
 \[
  L\cdot(h_n\cdot\Der+2\Der{h_n})\cdot\ldots\cdot(h_1\cdot\Der+2\Der{h_1}) = \E\cdot\tilde{L}.
 \]
\end{lemma}
\begin{proof}
 There are $n \in \mathbb{N}$, a nonzero $c \in \mathcal{C}$, a differential operator $L_0 \in \mathcal{R}\langle{\Der,\Int,\E}\rangle$, nonzero $f_1,\dots,f_n \in \Int\mathcal{R}$, and nonzero $g_1,\dots,g_n \in \mathcal{R}$ such that $g_1,\dots,g_n$ are linearly independent over $\mathcal{C}$ and $L\cdot{c}=\E\cdot{L_0}+\sum_{i=1}^n\E\cdot{f_i}\cdot\Int\cdot{g_i}$.
 By Lemma~\ref{lem:ConstantMultiples}, $L\cdot{c}=c\cdot{L}$ is nonzero.
\par
 To remove the sum $\sum_{i=1}^n\E\cdot{f_i}\cdot\Int\cdot{g_i}$, we will multiply $L\cdot{c}$ iteratively by $n$ first-order differential operators from the right as follows.
 If $n>0$, then we have
 \[
  \E\cdot{f_i}\cdot\Int\cdot{g_ig_n}\cdot\Der = \E\cdot{f_ig_ig_n} - \E{f_i}\cdot\E\cdot{g_ig_n} - \E\cdot{f_i}\cdot\Int\cdot((\Der{g_i})g_n+g_i\Der{g_n})
 \]
 for all $i$ by \eqref{eq:opIBP} and by \eqref{eq:opLinEval}. Therefore, by $\E{f_i}=0$, we obtain
 \begin{align*}
  L\cdot{c}\cdot(g_n\cdot\Der+2\Der{g_n}) &= \E\cdot{L_0}\cdot(g_n\cdot\Der+2\Der{g_n}) + \sum_{i=1}^n\E\cdot{f_i}\cdot\Int\cdot(g_ig_n\cdot\Der+2g_i\Der{g_n})\\
  &= \E\cdot{L_1} + \sum_{i=1}^{n-1}\E\cdot{f_i}\cdot\Int\cdot(g_i\Der{g_n}-(\Der{g_i})g_n),
 \end{align*}
 with $L_1:=L_0\cdot(g_n\cdot\Der+2\Der{g_n}) + \sum_{i=1}^nf_ig_ig_n$.
 Note that $g_i\Der{g_n}-(\Der{g_i})g_n$, for $i \in \{1,\dots,n-1\}$, are again linearly independent over $\mathcal{C}$ by Lemma~\ref{lem:DifferentialIntegralDomain}.
 Consequently, with $h_n:=cg_n$ being nonzero, the operator $L\cdot(h_n\cdot\Der+2\Der{h_n})$ is expressed like $L\cdot{c}$ above, just with $L_0$ replaced by $L_1$, with different $g_i$, and without the last summand.
\par
 By iterating this process of multiplying by a differential operator that is chosen as above, we obtain a right multiple of $L$ that is of the form $\E\cdot{L_n}$ with a differential operator $L_n$.
 It remains to show that $L_n$ is nonzero.
 After the last step, i.e.\ passing from $\E\cdot{L_{n-1}}+\E\cdot{f_1}\cdot\Int\cdot\tilde{g}_1$, with differential operator $L_{n-1}$ and nonzero $\tilde{g}_1 \in \mathcal{R}$, to $\E\cdot{L_n}$, we have $L_n=L_{n-1}\cdot(\tilde{g}_1\cdot\Der+2\Der\tilde{g}_1)+f_1\tilde{g}_1^2$.
 If the differential operator $L_{n-1}$ is nonzero, then also $L_n$ is nonzero, since $\tilde{g}_1\neq0$.
 Otherwise, we have $L_n=f_1\tilde{g}_1^2$, which is nonzero as well due to the assumptions.
\end{proof}

\begin{lemma}\label{lem:OperatorIdealsTwoSided}
 Let $(\mathcal{R},\Der,\Int)$ be an integro-differential ring that is an integral domain and let $I \subseteq \mathcal{R}\langle{\Der,\Int,\E}\rangle$ be a two-sided ideal. If $I$ contains an operator that is not an initial operator, then $I$ contains an operator of the form $f\cdot\E$ with nonzero $f \in \mathcal{R}$.
\end{lemma}
\begin{proof}
 By Lemma~\ref{lem:OperatorIdealsLeft}, $I$ contains a nonzero differential operator $L=\sum_{i=0}^nf_i\cdot\Der^i$. Let $k \in \{0,\dots,n\}$ be minimal such that $f_k\neq0$. Then, $f_k\cdot\E \in I$ since
 \[
  L\cdot\Int^k\cdot\E = \sum_{i=k}^nf_i\cdot\Der^i\cdot\Int^k\cdot\E = \sum_{i=k}^nf_i\cdot\Der^{i-k}\cdot\E = f_k\cdot\E. \qedhere
 \]
\end{proof}

Based on the above lemmas, we can find conditions when the action of $\mathcal{R}\langle{\Der,\Int,\E}\rangle$ on an $\mathcal{R}\langle{\Der,\Int,\E}\rangle$-module $M$ is faithful.
In this case, the annihilator
\[
 \ann_{\mathcal{R}\langle{\Der,\Int,\E}\rangle}(M)=\{L \in \mathcal{R}\langle{\Der,\Int,\E}\rangle\ |\ \forall{f\in{M}}:Lf=0\}
\]
is a two-sided ideal in $\mathcal{R}\langle{\Der,\Int,\E}\rangle$.
We are particularly interested in the action of IDO on integro-differential ring extensions of $(\mathcal{R},\Der,\Int)$.
In the theorem and its corollaries below, we impose conditions on the module $M$ by requiring some regularity on the action of $\mathcal{R}$ on $M$ and relating the action of $\E$ on $M$ to the constants $\mathcal{C}$.
By choosing $e=1$, modules $M=\mathcal{S}$ that are integro-differential ring extensions $(\mathcal{S},\Der,\Int)$ of $(\mathcal{R},\Der,\Int)$ and have the same constants $\mathcal{C}$ trivially satisfy these conditions.
For modules $M$ with these properties, we show that the annihilator only contains initial operators and, if $\mathcal{C}$ is a field, can be generated by monic initial operators.

\begin{theorem}\label{thm:FaithfulAction}
 Let $\mathcal{R}$ be an integro-differential ring that is an integral domain and let $\mathcal{C}$ be its ring of constants.
 Let $M$ be an $\mathcal{R}\langle{\Der,\Int,\E}\rangle$-module such that there exists a nonzero element $e \in \E M$ such that $\E M\subseteq\mathcal{C}e$ and $fe\neq0$ for all nonzero $f\in\mathcal{R}$.
 Then, $\mathcal{R}\langle{\Der,\Int,\E}\rangle$ acts faithfully on $M$ if and only if there is no nonzero differential operator $L \in \mathcal{R}\langle{\Der,\Int,\E}\rangle$ such that $\E\cdot{L}$ vanishes on all of $M$.
 In any case, the two-sided ideal $I = \ann_{\mathcal{R}\langle{\Der,\Int,\E}\rangle}(M)$ only contains initial operators.
 Moreover, if $\mathcal{C}$ is a field, $I$ has a generating set consisting only of monic initial operators.
\end{theorem}
\begin{proof}
 Let $m \in M$ such that $e=\E{m}$.
 No operator of the form $f\cdot\E$, with nonzero $f \in \mathcal{R}$, vanishes on $m\neq0$, since $(f\cdot\E)m=fe$ is nonzero by assumption.
 Hence, the ideal $I$ only contains initial operators by Lemma~\ref{lem:OperatorIdealsTwoSided}.
 If there is a nonzero differential operator $L \in \mathcal{R}\langle{\Der,\Int,\E}\rangle$ such that $\E\cdot{L}$ vanishes on all of $M$, then $\mathcal{R}\langle{\Der,\Int,\E}\rangle$ obviously does not act faithfully on $M$, since $\E\cdot{L}$ is nonzero as well.
\par
 Now, we assume that there is no nonzero differential operator $L$ such that $\E\cdot{L} \in I$.
 Let $L \in I$, then there are nonzero $c \in \mathcal{C}$, $f_1,\dots,f_n \in \mathcal{R}$, and nonzero monic initial operators $L_1,\dots,L_n \in \mathcal{R}\langle{\Der,\Int,\E}\rangle$ such that $f_1,\dots,f_n$ are $\mathcal{C}$-linearly independent and $c\cdot{L}=\sum_{i=1}^nf_i\cdot{L_i}$.
 For every $g \in M$, by the assumptions, there are $c_1,\dots,c_n \in \mathcal{C}$ such that $L_ig=c_ie$ for all $i$.
 For $f:=\sum_{i=1}^nf_ic_i \in \mathcal{R}$, since $Lg=0$, we have $fe = (\sum_{i=1}^nf_i\cdot{L_i})g = (c\cdot{L})g = 0$ and hence $f=0$ by assumption on $e$.
 Then, the assumption on $f_1,\dots,f_n$ implies $c_1=\dots=c_n=0$ and so $L_1g=\dots=L_ng=0$.
 Since $g$ was arbitrary, it follows that $L_1,\dots,L_n \in I$.
 If we would have $n\neq0$, then, by Lemma~\ref{lem:OperatorIdealsRight}, there would be a nonzero differential operator $\tilde{L} \in \mathcal{R}\langle{\Der,\Int,\E}\rangle$ such that $\E\cdot\tilde{L} \in I$ in contradiction to the assumption.
 Hence, $n=0$, which implies $c\cdot{L}=0$.
 By Lemma~\ref{lem:ConstantMultiples}, it follows that $L=0$, which shows that $\mathcal{R}\langle{\Der,\Int,\E}\rangle$ acts faithfully on $M$.
\par
 Finally, if $\mathcal{C}$ is a field, we show that $I$ has a generating set consisting only of monic initial operators.
 Let $L \in I$, then $L$ is an initial operator.
 Now, there are nonzero $f_1,\dots,f_n \in \mathcal{R}$, and nonzero monic initial operators $L_1,\dots,L_n \in \mathcal{R}\langle{\Der,\Int,\E}\rangle$ such that $f_1,\dots,f_n$ are $\mathcal{C}$-linearly independent and $L=\sum_{i=1}^nf_i\cdot{L_i}$.
 For all $g \in M$, there are $c_1,\dots,c_n \in \mathcal{C}$ with $L_ig=c_ie$ for all $i$ and, as above, we conclude $L_1g=\dots=L_ng=0$.
 Therefore, $L_1,\dots,L_n \in I$, which shows that $I$ is generated by nonzero monic initial operators.
\end{proof}

Assuming additional properties of $\mathcal{R}$ resp.\ of the action of $\mathcal{R}\langle{\Der,\Int,\E}\rangle$, generating sets of annihilators can be narrowed down even more.
In particular, the following two corollaries consider the situation when $\mathcal{R}$ is a field or $\E$ acts multiplicatively.

\begin{corollary}\label{cor:AnnihilatorGeneratorsField}
 Let $\mathcal{R}$ be an integro-differential field and let $\mathcal{C}$ be its ring of constants.
 Let $M$ be an $\mathcal{R}\langle{\Der,\Int,\E}\rangle$-module such that there exists a nonzero element $e \in \E M$ such that $\E M\subseteq\mathcal{C}e$ and $fe\neq0$ for all nonzero $f\in\mathcal{R}$.
 Then, the two-sided ideal $\ann_{\mathcal{R}\langle{\Der,\Int,\E}\rangle}(M)$ has a generating set consisting only of operators of the form $\E\cdot\sum_{i=0}^nf_i\cdot\Der^i$, with $f_0,\dots,f_n \in \mathcal{R}$ where $f_0 \in \Int\mathcal{R}$ is nonzero and $f_n=1$.
\end{corollary}
\begin{proof}
 By Theorem~\ref{thm:FaithfulAction}, $I:=\ann_{\mathcal{R}\langle{\Der,\Int,\E}\rangle}(M)$ has a generating set consisting only of monic initial operators.
 Now, let $L \in I$ be a monic initial operator and let $A$ be the set of all elements in $I$ that have the form specified in the statement of the corollary.
 By Lemma~\ref{lem:OperatorIdealsRight}, we obtain nonzero $h_1,\dots,h_m \in \mathcal{R}$ such that $L\cdot R = \E\cdot\tilde{L}$ with $R:=(h_m\cdot\Der+2\Der{h_m})\cdot\ldots\cdot(h_1\cdot\Der+2\Der{h_1}) \in \mathcal{R}\langle{\Der,\Int,\E}\rangle$ and some nonzero differential operator $\tilde{L} = \sum_{i=0}^n\tilde{f}_i\cdot\Der^i \in \mathcal{R}\langle{\Der,\Int,\E}\rangle$ with $\tilde{f}_n\neq0$.
 Then, by repeated use of \eqref{eq:opLeibniz}, there are $f_0,\dots,f_n\in\mathcal{R}$ with $f_n=1$ such that $\tilde{L}\cdot\tilde{f}_n^{-1} = \sum_{i=0}^nf_i\cdot\Der^i$.
 Let $k$ be minimal such that $f_k\neq0$, then $\tilde{L}\cdot\tilde{f}_n^{-1}\cdot\Int^k = \sum_{i=0}^{n-k}f_{i+k}\cdot\Der^i$.
 With $g \in M$ such that $e=\E{g}$, we compute $(\E\cdot\tilde{L}\cdot\tilde{f}_n^{-1}\cdot\Int^k)e = (\E\cdot\sum_{i=0}^{n-k}f_{i+k}\cdot\Der^i\cdot\E)g = (\E{f_k}\cdot\E)g = (\E{f_k})e$ by \eqref{eq:opLinEval} and \eqref{eq:opSpecial}.
 Since $\E\cdot\tilde{L}\cdot\tilde{f}_n^{-1}\cdot\Int^k \in I$ vanishes on $e \in M$, we conclude $f_k \in \Int\mathcal{R}$ by assumption on $e$ and hence $\E\cdot\sum_{i=0}^{n-k}f_{i+k}\cdot\Der^i \in A$.
 Since $\mathcal{R}$ is a field, we can verify by \eqref{eq:opLeibniz} and \eqref{eq:opFTC1} that $h_i^{-2}\cdot\Int\cdot{h_i} \in \mathcal{R}\langle{\Der,\Int,\E}\rangle$ is a right inverse of $h_i\cdot\Der+2\Der{h_i}$ for every $i$.
 Hence, there exists $\tilde{R} \in \mathcal{R}\langle{\Der,\Int,\E}\rangle$ such that $R\cdot\tilde{R} = 1$.
 Then, $\E\cdot\left(\sum_{i=0}^{n-k}f_{i+k}\cdot\Der^i\right)\cdot\Der^k\cdot\tilde{f}_n\cdot\tilde{R} = \E\cdot\tilde{L}\cdot\tilde{R} = L\cdot{R}\cdot\tilde{R} = L$, i.e.\ $L$ is in the ideal generated by $A$.
 Altogether, this shows that $I$ is generated by $A$.
\end{proof}

Whenever $\E$ is multiplicative on $\mathcal{R}$, the action of $\mathcal{R}\langle{\Der,\Int,\E}\rangle$ cannot be faithful on $M=\mathcal{R}$ since the operator $\E\cdot{f}$ acts as zero map for any $f\in\mathcal{R}$ with $\E{f}=0$.
By Lemma~\ref{lem:decompositionofring}, $\E{f}=0$ is equivalent to $f\in\Int\mathcal{R}$ and one can always choose $f=\Int1$, for example.
In Remark~\ref{rem:multiplicative}, we already discussed factoring the ring of operators by the relations \eqref{eq:opMulEval} immediately arising from multiplicativity of $\E$.
The following corollary shows that the resulting quotient always acts faithfully.

\begin{corollary}\label{cor:AnnihilatorGeneratorsMultiplicative}
 Let $\mathcal{R}$ be an integro-differential ring that is an integral domain and let $\mathcal{C}$ be its ring of constants.
 Let $M$ be an $\mathcal{R}\langle{\Der,\Int,\E}\rangle$-module such that there exists a nonzero element $e \in \E M$ such that $\E M\subseteq\mathcal{C}e$ and $fe\neq0$ for all nonzero $f\in\mathcal{R}$.
 If $(\E\cdot{f})g=(\E{f}\cdot\E)g$ for all $f\in\mathcal{R}$ and $g\in{M}$, then the two-sided ideal $\ann_{\mathcal{R}\langle{\Der,\Int,\E}\rangle}(M)$ is generated by the set $\{\E\cdot{f}-\E{f}\cdot\E\ |\ f\in\mathcal{R}\}$.
\end{corollary}
\begin{proof}
 Let $I:=\ann_{\mathcal{R}\langle{\Der,\Int,\E}\rangle}(M)$ and let $J\subseteq\mathcal{R}\langle{\Der,\Int,\E}\rangle$ be the two-sided ideal generated by the set $\{\E\cdot{f}-\E{f}\cdot\E\ |\ f\in\mathcal{R}\}$.
 By assumption on $\E$, we have that $J \subseteq I$, so it only remains to show $I \subseteq J$.
\par
 By Theorem~\ref{thm:FaithfulAction}, every $L \in I$ is an initial operator.
 Following the form of initial operators in Remark~\ref{rem:multiplicative}, there are $f_0,\dots,f_n\in\mathcal{R}$ such that $L=\sum_{i=0}^nf_i\cdot\E\cdot\Der^i$ modulo $J$.
 By $J \subseteq I$, we have that $\sum_{i=0}^nf_i\cdot\E\cdot\Der^i \in I$.
 Therefore, $\sum_{i=0}^nf_i\cdot\E\cdot\Der^i$ vanishes on $\Int^ke \in M$ for all $k\in\{0,\dots,n\}$.
 This implies inductively that $f_0,\dots,f_n$ are all zero, since, by \eqref{eq:opFTC1} and \eqref{eq:opSpecial}, we can compute $(\sum_{i=k}^nf_i\cdot\E\cdot\Der^i)\Int^ke = (\sum_{i=0}^{n-k}f_{i+k}\cdot\E\cdot\Der^{i+k}\cdot\Int^k\cdot\E)m = (f_k\cdot\E)m = f_ke$ with $m \in M$ such that $e=\E{m}$.
 Consequently, we have $L \in J$.
\end{proof}

\section{Equational prover in calculus}
\label{sec:equationalprover}

In this section, we illustrate how results from analysis can be proven via computations in the ring of integro-differential operators.
The general approach consists in formulating an analytic statement as an identity of integro-differential operators and then prove this identity algebraically in $\mathcal{R}\langle{\Der,\Int,\E}\rangle$, with minimal assumptions on the integro-differential ring $(\mathcal{R},\Der,\int)$ used in the coefficients.
Note that we prove results directly by a computation with integro-differential operators, instead of doing the whole computation with elements of $\mathcal{R}$ or any other left $\mathcal{R}\langle{\Der,\Int,\E}\rangle$-module.

An identity in $\mathcal{R}\langle{\Der,\Int,\E}\rangle$ can be proven by comparing irreducible forms of the left hand side and right hand side. Recall that such irreducible forms are given by Theorem~\ref{thm:IDO} and can be computed systematically by the rewrite rules in Table~\ref{tab:op}.
In this sense, the rewrite rules provide an equational prover for integro-differential operators provided one can decide equality of irreducible forms in $\mathcal{R}\langle{\Der,\Int,\E}\rangle$, which includes deciding identities in $\mathcal{R}$.
In practice, this is often possible for concrete irreducible forms.

Once an identity $L_1=L_2$ is proven in $\mathcal{R}\langle{\Der,\Int,\E}\rangle$, we immediately infer the corresponding identity in $\mathcal{R}$, i.e.\ $L_1f=L_2f$ for all $f \in \mathcal{R}$, by the canonical action of $\mathcal{R}\langle{\Der,\Int,\E}\rangle$ on $\mathcal{R}$. Moreover, by acting on any other left $\mathcal{R}\langle{\Der,\Int,\E}\rangle$-module, we obtain the analogous identity also in those modules.
Furthermore, any concrete computation with operators acting on functions uses only finitely many derivatives.
Consequently, by inspecting every step of the computation, an identity proven for infinitely differentiable functions can also be proven for functions that are only sufficiently often differentiable.

\subsection{Variation of constants for scalar equations}

In this section, we deal with the method of variation of constants for computing solutions of inhomogeneous ODEs in terms of integro-differential operators. First, we recall the analytic statement, see e.g.\ Theorem~6.4 in Chapter~3 of \cite{CoddingtonLevinson1955}.
Consider the inhomogeneous linear ODE
\[
 y^{(n)}(x)+a_{n-1}(x)y^{(n-1)}(x)+\dots+a_0(x)y(x)=f(x).
\]
Assume that $z_1(x),\dots,z_n(x)$ is a fundamental system of the homogeneous equation, i.e.\ $z_i^{(n)}(x)+a_{n-1}(x)z_i^{(n-1)}(x)+\dots+a_0(x)z_i(x)=0$ and the Wronskian $w(x):=W(z_1(x),\dots,z_n(x))$ is nonzero. Then,
\begin{equation}\label{eq:varconstsol}
 z^*(x) := \sum_{i=1}^n(-1)^{n-i}z_i(x)\int_{x_0}^x\frac{W(z_1(t),\dots,z_{i-1}(t),z_{i+1}(t),\dots,z_n(t))}{w(t)}f(t)\,dt
\end{equation}
is a particular solution of the inhomogeneous equation above.

Algebraically, we model scalar functions by fixing a commutative integro-differential ring $(\mathcal{R},\Der,\Int)$ and, as indicated above, we model computations with scalar equations by computations in the corresponding ring of integro-differential operators $\mathcal{R}\langle{\Der,\Int,\E}\rangle$.
From the operator viewpoint, mapping the inhomogeneous part $f$ to a solution of the equation $L y = f$ amounts to constructing a right inverse of the differential operator $L$. In general, for an operator $L$ and a right inverse $H$, a particular solution of the inhomogeneous equation is given by $z^*= H f$, since
\[
 Lz^* = L(Hf) = (L \cdot H)f = 1f = f.
\]

Recall that $\Int$ is a right inverse of $\Der$ by definition. Now, we consider an arbitrary monic first-order differential operator
\[
 L=\Der+a,
\]
with $a\in\mathcal{R}$, and assume that $z \in \mathcal{R}$ is a solution of $Ly=0$, i.e.\ $\Der z + a z=0$.
Then, using \eqref{eq:opLeibniz}, we compute
\[
 (\Der + a) \cdot z = \Der \cdot z + a \cdot z = z \cdot \Der + \Der z + a z = z \cdot \Der.
\]
If, moreover, we assume that $z$ has a multiplicative inverse $z^{-1} \in \mathcal{R}$, we obtain
\[
 (\Der + a) \cdot (z \cdot \Int \cdot z^{-1}) =  z \cdot \Der \cdot \Int \cdot z^{-1} = z \cdot z^{-1} = 1.
\]
Hence 
\begin{equation*}
 H = z \cdot \Int \cdot z^{-1}
\end{equation*}
is a right inverse of $L$ in $\mathcal{R}\langle{\Der,\Int,\E}\rangle$.

We also outline the computation for a second order differential operator
\[
 L=\Der^2+a_1\cdot\Der+a_0,
\]
with $a_1,a_0\in\mathcal{R}$.  If $z \in \mathcal{R}$ is a solution of $Ly=0$,  then
\[
 L \cdot z = z \cdot \Der^2 + (2\Der{z} + a_1z) \cdot \Der.
\]
Assume that there exist two solutions $z_1,z_2 \in \mathcal{R}$ of $Ly=0$ such that their Wronskian
\[
 w = z_1 \Der z_2 - z_2 \Der z_1
\]
has a multiplicative inverse $\tfrac{1}{w} \in \mathcal{R}$. Let
\[
 H = -z_1 \cdot \Int \cdot \tfrac{z_2}{w} + z_2 \cdot \Int \cdot \tfrac{z_1}{w}.
\]
In the ring of operators, we can compute
\[
 \Der \cdot \tfrac{z_i}{w} = \tfrac{z_i}{w} \cdot \Der + \tfrac{\Der z_i}{w} - \tfrac{z_i\Der w}{w^2}.
\]
Then, using the normal forms of $L \cdot z_i$ and $\Der \cdot \tfrac{z_i}{w}$, we obtain
\begin{multline*}
 L \cdot H = -z_1 \cdot \Der \cdot \tfrac{z_2}{w} - 2(\Der z_1)\tfrac{z_2}{w} - a_1z_1\tfrac{z_2}{w} + z_2 \cdot \Der \cdot \tfrac{z_1}{w} + 2(\Der z_2)\tfrac{z_1}{w} + a_1z_2\tfrac{z_1}{w}\\
 = -z_1 \cdot \Der \cdot \tfrac{z_2}{w} + z_2 \cdot \Der \cdot \tfrac{z_1}{w} + 2 \tfrac{w}{w}
 = -z_1 \Der \tfrac{z_2}{w} + z_2 \Der \tfrac{z_1}{w} + 2 = 1.
\end{multline*}
Hence $H$ is a right inverse of $L$ in $\mathcal{R}\langle{\Der,\Int,\E}\rangle$.

More generally, we have the following formulation of the method of variation of constants for integro-differential operators.
Recall from Section~\ref{sec:IDOaction}, that the Wronskian $W(f_1,\dots,f_n)$ of elements $f_1,\dots,f_n \in \mathcal{R}$ is defined completely analogous to the analytic situation.

\begin{theorem}\label{thm:varconstCommutative}
 Let $(\mathcal{R},\Der,\Int)$ be a commutative integro-differential ring and let $L=\Der^n+\sum_{i=0}^{n-1}a_i\cdot\Der^i \in \mathcal{R}\langle{\Der,\Int,\E}\rangle$, $n\ge1$, with $a_0,\dots,a_{n-1} \in \mathcal{R}$. Assume that $z_1,\dots,z_n \in \mathcal{R}$ are such that $Lz_i=0$ and $w:=W(z_1,\dots,z_n)$ has a multiplicative inverse $\tfrac{1}{w} \in \mathcal{R}$. Then, with
 \[
  H := \sum_{i=1}^n(-1)^{n-i}z_i\cdot\Int\cdot\frac{W(z_1,\dots,z_{i-1},z_{i+1},\dots,z_n)}{w}
 \]
 we have that $L \cdot H = 1$ in $\mathcal{R}\langle{\Der,\Int,\E}\rangle$.
\end{theorem}
\begin{proof}
 The cases $n=1$ and $n=2$ have been shown above. In principle, an analogous computation could be done for any concrete $n\ge3$ as well. To obtain a finite proof for all $n\ge3$ at once, we will utilize a more general framework in Section~\ref{sec:varconst}.
\end{proof}

\subsection{Linear systems and operators with matrix coefficients}
\label{sec:linearsystems}

Algebraically, computing with linear systems of differential equations can be modelled by integro-differential operators over some noncommutative integro-differential ring, whose elements correspond to matrices.
Such noncommutative integro-differential rings can be obtained from any scalar integro-differential ring, since one can equip the ring of $n\times n$ matrices with a derivation and an integration defined entrywise, see Lemma~\ref{lem:matrixring}.

Recall that Definition~\ref{def:OperatorRing} and Theorem~\ref{thm:IDO} hold also for operators with coefficients in noncommutative integro-differential rings, in particular for operators with matrix coefficients. Consequently, any computation using a general noncommutative integro-differential ring is automatically valid for concrete matrices of any size, without the need for entrywise computations. This allows for compact proofs for matrices of general size.
For switching to entrywise computations, one can view operators with matrix coefficients equivalently also as matrices of operators with scalar coefficients by identifying operators $\Der,\Int,\E$ with corresponding diagonal matrices of operators.
The following lemma shows that also computations are equivalent in both viewpoints.
Here, we use the notation $E_{i,j}(L):=(\delta_{i,k}\delta_{j,l}L)_{k,l=1,\dots,n}$ for matrices with only one nonzero entry $L \in \mathcal{R}\langle{\Der,\Int,\E}\rangle$.

\begin{lemma}\label{lem:matrixisomorphism}
 Let $(\mathcal{R},\Der,\Int)$ an integro-differential ring and let $n\ge1$. Then, there is exactly one ring homomorphism
 \[
  \varphi:\mathcal{R}^{n\times n}\langle{\Der,\Int,\E}\rangle \to \mathcal{R}\langle{\Der,\Int,\E}\rangle^{n\times n}
 \]
 with $\varphi(A)=A$ for $A \in \mathcal{R}^{n\times n}$ and $\varphi(L)=\diag(L,\dots,L)$ for $L \in \{\Der,\Int,\E\}$.
 This $\varphi$ is an isomorphism and its inverse homomorphism
 \[
  \psi:\mathcal{R}\langle{\Der,\Int,\E}\rangle^{n\times n} \to \mathcal{R}^{n\times n}\langle{\Der,\Int,\E}\rangle
 \]
 can be given by $\psi(E_{i,j}(f))=E_{i,j}(f)$ and $\psi(E_{i,j}(L))=E_{i,j}(1)\cdot{L}$ for all $f \in \mathcal{R}$ resp.\ $L \in \{\Der,\Int,\E\}$.
\end{lemma}
\begin{proof}
 First, we check that the definition of $\varphi$ indeed provides a unique and well-defined homomorphism of rings. Uniqueness follows from the fact that $\varphi$ is defined on a generating set of $\mathcal{R}^{n\times n}\langle{\Der,\Int,\E}\rangle$. For proving well-definedness, we need to verify that the definition of $\varphi$ respects all identities of generators given in Definition~\ref{def:OperatorRing} as follows.
 Evidently, we have $\varphi(1)=I_n$.
 It is immediate to see that $\varphi(\Der)\cdot\varphi(\Int)=\varphi(1)$ and $\varphi(\Int)\cdot\varphi(\Der)=\varphi(1)-\varphi(\E)$ hold.
 For $L \in \{\Der,\Int,\E\}$, we verify in $\mathcal{R}\langle{\Der,\Int,\E}\rangle^{n\times n}$ that $\varphi(L)$ commutes with all elements of $\mathcal{C}^{n\times n}$ and that, for all $A \in \mathcal{R}^{n\times n}$, we have $\varphi(L)\cdot\varphi(A)\cdot\varphi(\E)=\varphi(LA)\cdot\varphi(\E)$.
 More explicitly, using \eqref{eq:opLinDer}--\eqref{eq:opLinEval} in $\mathcal{R}\langle{\Der,\Int,\E}\rangle$, we compute
 \begin{multline*}
  \varphi(L)\cdot\varphi(A)\cdot\varphi(\E) = \diag(L,\dots,L)\cdot{A}\cdot\diag(\E,\dots,\E)\\
  = \left(L\cdot{a_{i,j}}\cdot\E\right)_{i,j=1,\dots,n} = \left(L{a_{i,j}}\cdot\E\right)_{i,j=1,\dots,n} = \varphi(LA)\cdot\varphi(\E).
 \end{multline*}
 Analogously, using \eqref{eq:opLeibniz} in $\mathcal{R}\langle{\Der,\Int,\E}\rangle$, we can also verify that $\varphi(\Der)\cdot\varphi(A)=\varphi(A)\cdot\varphi(\Der)+\varphi(\Der{A})$ for all $A \in \mathcal{R}^{n\times n}$.\par
 Similarly, we need to verify that the definition of $\psi$ on a generating set of $\mathcal{R}\langle{\Der,\Int,\E}\rangle^{n\times n}$ gives rise to a well-defined homomorphism. For example, using \eqref{eq:opLinDer}--\eqref{eq:opLinEval} in $\mathcal{R}^{n\times n}\langle{\Der,\Int,\E}\rangle$, we compute
 \begin{multline*}
  \psi(E_{i,j}(L))\cdot\psi(E_{p,q}(f))\cdot\psi(E_{k,l}(\E)) = E_{i,j}(1)\cdot{L}\cdot{E_{p,q}(f)}\cdot{E_{k,l}(1)}\cdot\E\\
  = E_{i,j}(1)\cdot{L}\cdot\delta_{q,k}E_{p,l}(f)\cdot\E = E_{i,j}(1)\cdot\delta_{q,k}LE_{p,l}(f)\cdot\E = \delta_{j,p}\delta_{q,k}E_{i,l}(Lf)\cdot\E\\
  = E_{i,q}(\delta_{j,p}Lf)\cdot{E_{k,l}(1)}\cdot\E = \psi(E_{i,q}(\delta_{j,p}Lf))\cdot\psi(E_{k,l}(\E))
 \end{multline*}
 for all $f \in \mathcal{R}$, $L \in \{\Der,\Int,\E\}$, and all $i,j,k,l,p,q \in \{1,\dots,n\}$.\par
 Finally, $\psi$ is the inverse of $\varphi$, since we have $\psi(\varphi(A))=A$, $\varphi(\psi(E_{i,j}(f)))=E_{i,j}(f)$, $\psi(\varphi(L))=L$, and $\varphi(\psi(E_{i,j}(L)))=E_{i,j}(L)$ for the generators.
\end{proof}

Integro-differential operators with coefficients from $\mathcal{R}^{n\times n}$ constructed this way have natural actions on $\mathcal{R}^{n\times n}$ and on $\mathcal{R}^n$.
By the above lemma, for any left $\mathcal{R}\langle{\Der,\Int,\E}\rangle$-module $M$, there is a natural way of viewing $M^n$ as a left $\mathcal{R}^{n\times n}\langle{\Der,\Int,\E}\rangle$-module.
It is straightforward to show that the action of $\mathcal{R}^{n\times n}\langle{\Der,\Int,\E}\rangle$ on $M^n$ is faithful if and only if $\mathcal{R}\langle{\Der,\Int,\E}\rangle$ acts faithfully on $M$.

\subsection{Variation of constants for first-order systems}
\label{sec:varconst}

Instead of higher order scalar ODEs, we now consider variation of constants for first-order systems, see Theorem~3.1 in Chapter~3 of \cite{CoddingtonLevinson1955} for example. Analytically, it can be stated as follows.
For an $n\times n$ matrix $A(x)$ and a vector $f(x)$ of size $n$, we consider the first-order system given by
\[
 y^\prime(x)+A(x)y(x) = f(x).
\]
If $\Phi(x)$ is a fundamental matrix of the homogeneous system, i.e.\ it satisfies $\Phi^\prime(x)+A(x)\Phi(x) = 0$ and $\det(\Phi(x))\neq0$, then a particular solution of the inhomogeneous system is given by
\[
 z^*(x) = \Phi(x) \int_{x_0}^x\Phi(t)^{-1}f(t)\,dt.
\]

\begin{theorem}\label{thm:varconst}
 Let $(\mathcal{R},\Der,\Int)$ be an integro-differential ring. Assume that $a\in\mathcal{R}$ is such that there exists $z \in \mathcal{R}$ that satisfies $\Der z + a z=0$ and has a multiplicative right inverse $z^{-1} \in \mathcal{R}$.
 Then, in $\mathcal{R}\langle{\Der,\Int,\E}\rangle$, the operators $L := \Der + a$ and $H := z \cdot \Int \cdot z^{-1}$ satisfy
 \[
  L \cdot H = 1.
 \]
\end{theorem}
\begin{proof}
 The same computation as in the commutative case above can be done also for noncommutative $\mathcal{R}$ without any changes.
 Note that $zz^{-1}=1$ was the only property of $z^{-1}$ used there, so it suffices that $z^{-1}$ is a right inverse of $z$.
\end{proof}

Observe that from this statement over an abstract integro-differential ring $(\mathcal{R},\Der,\Int)$, which was proven without referring to matrices at all, the analytic statement follows for arbitrary size $n$ of the matrix $A(x)$, provided we assume sufficient regularity of the functions involved.

Based on Theorem~\ref{thm:varconst}, we now can complete the proof of Theorem~\ref{thm:varconstCommutative} for arbitrary $n\ge3$.
Before doing so, we first detail the required translation from a first-order system to the scalar equation entirely at the operator level. For shorter notation, in the following, we again use the symbol $\Der$ also for the matrix $\diag(\Der,\dots,\Der) \in \mathcal{R}\langle{\Der,\Int,\E}\rangle^{n\times n}$.

\begin{lemma}\label{lem:companionsolution}
 Let $(\mathcal{R},\Der,\Int)$ be an integro-differential ring and let
 \[
  L=\Der^n+\sum_{i=0}^{n-1}a_i\cdot\Der^i \in \mathcal{R}\langle{\Der,\Int,\E}\rangle
 \]
 with $a_0,\dots,a_{n-1} \in \mathcal{R}$.
 Moreover, let $H\in\mathcal{R}\langle{\Der,\Int,\E}\rangle^{n\times n}$ satisfy $(\Der+A) \cdot H = I_n$ in $\mathcal{R}\langle{\Der,\Int,\E}\rangle^{n\times n}$, where
 \[
  A :=
  \begin{pmatrix}
   0&-1&0&\cdots&0\\
   \vdots&\ddots&\ddots&\ddots&\vdots\\
   \vdots&&\ddots&\ddots&0\\
   0&\cdots&\cdots&0&-1\\
   a_0&a_1&\cdots&\cdots&a_{n-1}
  \end{pmatrix}.
 \]
 Then, the upper right entry $H_{1,n}$ of $H$ satisfies $L \cdot H_{1,n} = 1$ in $\mathcal{R}\langle{\Der,\Int,\E}\rangle$ and $H_{i,n}=\Der^{i-1}\cdot{H}_{1,n}$ for $i\in\{1,\dots,n\}$.
\end{lemma}
\begin{proof}
 For $n=1$, the statement is trivial. So, we let $n\ge2$ in the following.
 In short, starting from the identity $(\Der+A) \cdot H = I_n$ of $n\times n$ matrices of operators, we multiply both sides from the left with a suitable $S \in \mathcal{R}\langle{\Der,\Int,\E}\rangle^{n\times n}$ and inspect the last column. In particular, for elimination in the last $n-1$ columns of $\Der+A$, we use
 \[
  S := \begin{pmatrix}1&0&\cdots&\cdots&\cdots&0\\
  \Der&\ddots&\ddots&&&\vdots\\
  \Der^2&\ddots&\ddots&\ddots&&\vdots\\
  \vdots&\ddots&\ddots&\ddots&\ddots&\vdots\\
  \Der^{n-2}&\cdots&\Der^2&\Der&1&0\\
  s_1&\cdots&\cdots&\cdots&s_{n-1}&1\end{pmatrix}
 \]
 with $s_j := \Der^{n-j}+\sum_{i=0}^{n-j-1}a_{i+j}\cdot\Der^i \in \mathcal{R}\langle{\Der,\Int,\E}\rangle$ for $j\in\{1,\dots,n-1\}$. Observing that, in $\mathcal{R}\langle{\Der,\Int,\E}\rangle$, we have $-s_j+s_{j+1}\cdot\Der=-a_j$ for $j\in\{1,\dots,n-2\}$, we compute $S\cdot(\Der+A)$ explicitly:
 \begin{multline*}
  S \cdot
  \begin{pmatrix}
   \Der&-1&0&\cdots&0\\
   0&\ddots&\ddots&\ddots&\vdots\\
   \vdots&\ddots&\ddots&\ddots&0\\
   0&\cdots&0&\Der&-1\\
   a_0&a_1&\cdots&a_{n-2}&\Der+a_{n-1}
  \end{pmatrix}\\
  =
  \begin{pmatrix}
   \Der&-1&0&\cdots&0\\
   \Der^2&0&-1&\ddots&\vdots\\
   \vdots&\vdots&\ddots&\ddots&0\\
   \Der^{n-1}&0&\cdots&0&-1\\
   s_1\cdot\Der+a_0 & 0 & \cdots & 0 & -s_{n-1}+\Der+a_{n-1}
  \end{pmatrix}\\
  =
  \begin{pmatrix}
   \Der&-1&0&\cdots&0\\
   \Der^2&0&-1&\ddots&\vdots\\
   \vdots&\vdots&\ddots&\ddots&0\\
   \Der^{n-1}&\vdots&&\ddots&-1\\
   L & 0 & \cdots & \cdots & 0
  \end{pmatrix}.
 \end{multline*}
 Altogether, we obtain that
 \[
  \begin{pmatrix}
   \Der&-1&0&\cdots&0\\
   \Der^2&0&-1&\ddots&\vdots\\
   \vdots&\vdots&\ddots&\ddots&0\\
   \Der^{n-1}&\vdots&&\ddots&-1\\
   L & 0 & \cdots & \cdots & 0
  \end{pmatrix}
  \cdot H = S \cdot (\Der+A) \cdot H = S,
 \]
 where comparison of the entries in the last column of the left hand side and right hand side yields $\Der \cdot H_{1,n} - H_{2,n} = 0,\dots,\Der^{n-1} \cdot H_{1,n} - H_{n,n} = 0$ and $L \cdot H_{1,n} = 1$.
\end{proof}

Finally, we proceed with the proof of Theorem~\ref{thm:varconstCommutative}.
Recalling the assumptions, we fix a commutative integro-differential ring $(\mathcal{R},\Der,\Int)$ and we let $L=\Der^n+\sum_{i=0}^{n-1}a_i\cdot\Der^i \in \mathcal{R}\langle{\Der,\Int,\E}\rangle$ with $a_0,\dots,a_{n-1} \in \mathcal{R}$.
We also let $z_1,\dots,z_n \in \mathcal{R}$ such that $Lz_i=0$ and such that $w:=W(z_1,\dots,z_n)$ has a multiplicative inverse $\tfrac{1}{w} \in \mathcal{R}$.

\begin{proof}[Proof of Theorem~\ref{thm:varconstCommutative}]
 Let $n\ge3$.
 In the (noncommutative) integro-differential ring $(\mathcal{R}^{n\times n},\Der,\Int)$, we consider the matrices
 \[
  A:=
  \begin{pmatrix}
   0&-1&0&\cdots&0\\
   \vdots&\ddots&\ddots&\ddots&\vdots\\
   \vdots&&\ddots&\ddots&0\\
   0&\cdots&\cdots&0&-1\\
   a_0&a_1&\cdots&\cdots&a_{n-1}
  \end{pmatrix}
  \quad\text{and}\quad
  Z:=
  \begin{pmatrix}
   z_1&\cdots&z_n\\
   \Der{z_1}&\cdots&\Der{z_n}\\
   \vdots&&\vdots\\
   \Der^{n-1}z_1&\cdots&\Der^{n-1}z_n
  \end{pmatrix}
 \]
 and note that $(\Der+A)Z=0$ and that $Z^{-1} \in \mathcal{R}^{n\times n}$ exists, since $\det(Z)=w$ was assumed to be invertible in $\mathcal{R}$.
 Then, we apply Theorem~\ref{thm:varconst} to obtain $(\Der+A)\cdot(Z\cdot\Int\cdot{Z}^{-1})=1$ in $\mathcal{R}^{n\times n}\langle{\Der,\Int,\E}\rangle$.
 Instead of integro-differential operators with matrix coefficients in $\mathcal{R}^{n\times n}$, we now consider the objects as $n\times n$ matrices whose entries are integro-differential operators with coefficients in $\mathcal{R}$, cf.\ Lemma~\ref{lem:matrixisomorphism}.
 Then, by Lemma~\ref{lem:companionsolution}, we obtain that $L\cdot(Z\cdot\Int\cdot{Z}^{-1})_{1,n}=1$ holds in $\mathcal{R}\langle{\Der,\Int,\E}\rangle$.
 In order to compute the entry $(Z\cdot\Int\cdot{Z}^{-1})_{1,n}$, we can easily determine a general form of the entries of the matrix product $Z\cdot\Int\cdot{Z}^{-1} \in \mathcal{R}\langle{\Der,\Int,\E}\rangle^{n\times n}$:
 \begin{multline*}
  Z \cdot \begin{pmatrix}\Int&0&\cdots&0\\0&\ddots&\ddots&\vdots\\\vdots&\ddots&\ddots&0\\0&\cdots&0&\Int\end{pmatrix} \cdot{Z}^{-1} = Z\cdot\begin{pmatrix}\Int\cdot(Z^{-1})_{1,1}&\cdots&\Int\cdot(Z^{-1})_{1,n}\\\vdots&&\vdots\\\Int\cdot(Z^{-1})_{n,1}&\cdots&\Int\cdot(Z^{-1})_{n,n}\end{pmatrix}\\
  = \left(\sum_{k=1}^n(\Der^{i-1}z_k)\cdot\Int\cdot(Z^{-1})_{k,j}\right)_{i,j=1,\dots,n}.
 \end{multline*}
 Then, we compute all entries
 \[
  (Z^{-1})_{i,n}=(-1)^{n+i}\frac{W(z_1,\dots,z_{i-1},z_{i+1},\dots,z_n)}{w}
 \]
 in the last column of $Z^{-1} \in \mathcal{R}^{n\times n}$ via Cramer's rule (or via the cofactor matrix) so that we recognize $H \in \mathcal{R}\langle{\Der,\Int,\E}\rangle$ defined in the statement of Theorem~\ref{thm:varconstCommutative} as the top right entry $(Z\cdot\Int\cdot{Z}^{-1})_{1,n}$ of the above matrix. This concludes the proof that $L \cdot H = 1$.
\end{proof}

\section{Generalizing identities from calculus}
\label{sec:generalizations}

The generalizations of well-known identities presented in this section introduce additional terms involving the induced evaluation, which vanish if the evaluation is multiplicative.
As explained in the previous section, the proofs mostly rely on computing irreducible forms for IDO.

\subsection{Initial value problems}
\label{sec:IVP}

Recall that the formula given in \eqref{eq:varconstsol} provides a particular solution of the inhomogeneous linear ODE
\[
 y^{(n)}(x)+a_{n-1}(x)y^{(n-1)}(x)+\dots+a_0(x)y(x)=f(x).
\]
Moreover, this particular solution also satisfies the homogeneous initial conditions $y(x_0)=0,y^\prime(x_0)=0,\dots,y^{(n-1)}(x_0)=0$, see e.g.\ Theorem~6.4 in Chapter~3 of \cite{CoddingtonLevinson1955}.
The induced evaluation of an integro-differential ring allows us to model such properties algebraically. Similarly to the general proof of Theorem~\ref{thm:varconstCommutative}, we first investigate the general solution formula for homogeneous initial value problems of first-order systems. To this end, we fix a (not necessarily commutative) integro-differential ring $\mathcal{R}$ and we work in the corresponding ring of IDO $\mathcal{R}\langle{\Der,\Int,\E}\rangle$.

As a general principle, not only in the ring $\mathcal{R}\langle{\Der,\Int,\E}\rangle$ but in any ring with unit element, if some $H$ is a right inverse of some $L$ and $B,P$ are such that $L \cdot P = 0$ and $B \cdot P = B$, then $G := (1-P) \cdot H$ satisfies
\[
 L \cdot G = 1 \quad\text{and}\quad B \cdot G = 0.
\]
So, by choosing an appropriate operator $P \in \mathcal{R}\langle{\Der,\Int,\E}\rangle$ that satisfies $(\Der+a) \cdot P = 0$ and $\E \cdot P = \E$, we can obtain the following version of Theorem~\ref{thm:varconst} that also solves the homogeneous initial condition.

\begin{theorem}\label{thm:IVP}
 Let $(\mathcal{R},\Der,\Int)$ be an integro-differential ring and let $L = \Der+a \in \mathcal{R}\langle{\Der,\Int,\E}\rangle$ with $a\in\mathcal{R}$.
 Assume that $z\in\mathcal{R}$ is such that $\Der z + a z=0$ and in addition to a multiplicative right inverse $z^{-1}\in\mathcal{R}$ also a right inverse $(\E{z})^{-1}\in\mathcal{C}$ exists.
 Then, in $\mathcal{R}\langle{\Der,\Int,\E}\rangle$, the operator
 \[
  G := (1-z(\E{z})^{-1} \cdot \E) \cdot z \cdot \Int \cdot z^{-1}
 \]
 satisfies $L \cdot G = 1$ and $\E \cdot G = 0$.
\end{theorem}
\begin{proof}
 By Theorem~\ref{thm:varconst}, we have that $H := z \cdot \Int \cdot z^{-1}$ satisfies $L \cdot H = 1$. Now, letting $P := z(\E{z})^{-1} \cdot \E$, we compute the normal forms of $L \cdot P$ and $\E \cdot P$:
 \begin{align*}
  (\Der+a) \cdot z(\E{z})^{-1} \cdot \E &= \left(z(\E{z})^{-1} \cdot \Der + \Der{z}(\E{z})^{-1}\right) \cdot \E + az(\E{z})^{-1} \cdot \E\\
  &= z(\E{z})^{-1} \cdot \Der \cdot \E = 0,\\
  \E \cdot z(\E{z})^{-1} \cdot \E &= \E{z}(\E{z})^{-1} \cdot \E = \E.
 \end{align*}
 Then, from $L \cdot P = 0$ and $\E \cdot P = \E$, it follows straightforwardly that $G = (1-P) \cdot H$ has the properties claimed.
\end{proof}

\begin{remark}
 If the evaluation $\E$ is multiplicative, then the existence of $(\E{z})^{-1}\in\mathcal{C}$ follows form the existence of $z^{-1}\in\mathcal{R}$ and, with \eqref{eq:opMulEval}, we also obtain $\E \cdot z \cdot \Int = (\E{z})\E \cdot \Int = 0$, which implies $P \cdot H = 0$ and hence $G=H$.
 Therefore, the right inverse $H$ obtained in Theorem~\ref{thm:varconst} satisfies the initial value condition $\E \cdot H=0$ already.
 \qed
\end{remark}

For conversion to the scalar case, we need to compute the element in the top right corner of $G = (1-Z(\E{Z})^{-1}\cdot\E)\cdot{Z}\cdot\Int\cdot{Z}^{-1}$ and show that it has the desired properties.

\begin{theorem}
 Let $(\mathcal{R},\Der,\Int)$ be a commutative integro-differential ring and let $L=\Der^n+\sum_{i=0}^{n-1}a_i\cdot\Der^i \in \mathcal{R}\langle{\Der,\Int,\E}\rangle$, $n\ge1$, with $a_0,\dots,a_{n-1} \in \mathcal{R}$. Assume that $z_1,\dots,z_n \in \mathcal{R}$ are such that $Lz_i=0$ and $w:=W(z_1,\dots,z_n)$ has a multiplicative inverse $\tfrac{1}{w} \in \mathcal{R}$. Moreover, assume that there are $c_{i,j}\in\mathcal{C}$ such that $\E\Der^k\sum_{i=1}^nc_{i,j}z_i=\delta_{j,k}$ for all $j,k\in\{0,\dots,n-1\}$. Then, with these $c_{i,j}$ and
 \[
  G := \sum_{k=1}^n(-1)^{n-k}\bigg(z_k-\sum_{i,j=1}^nz_ic_{i,j-1}\cdot\E\cdot(\Der^{j-1}z_k)\bigg)\cdot\Int\cdot\tfrac{W(z_1,\dots,z_{k-1},z_{k+1},\dots,z_n)}{w}
 \]
 we have that $L \cdot G = 1$ and $\E \cdot G = 0$ in $\mathcal{R}\langle{\Der,\Int,\E}\rangle$.
\end{theorem}
\begin{proof}
 As in the proof of Theorem~\ref{thm:varconstCommutative}, we consider the matrices
 \[
  A:=
  \begin{pmatrix}
   0&-1&0&\cdots&0\\
   \vdots&\ddots&\ddots&\ddots&\vdots\\
   \vdots&&\ddots&\ddots&0\\
   0&\cdots&\cdots&0&-1\\
   a_0&a_1&\cdots&\cdots&a_{n-1}
  \end{pmatrix}
  \quad\text{and}\quad
  Z:=
  \begin{pmatrix}
   z_1&\cdots&z_n\\
   \Der{z_1}&\cdots&\Der{z_n}\\
   \vdots&&\vdots\\
   \Der^{n-1}z_1&\cdots&\Der^{n-1}z_n
  \end{pmatrix}
 \]
 and note that $(\Der+A)Z=0$ and that $Z^{-1} \in \mathcal{R}^{n\times n}$ exists, since $\det(Z)=w$ was assumed to be invertible in $\mathcal{R}$.
 Furthermore, we note that
 \[
  \begin{pmatrix}
   c_{1,0} & \cdots & c_{1,n-1}\\
   \vdots && \vdots\\
   c_{n,0} & \cdots & c_{n,n-1}
  \end{pmatrix}
 \]
 is the multiplicative (right) inverse of $\E{Z}$ in $\mathcal{C}^{n\times n}$.
 By Theorem~\ref{thm:IVP}, we conclude that $\tilde{G} := (1-Z(\E{Z})^{-1} \cdot \E) \cdot Z \cdot \Int \cdot Z^{-1} \in \mathcal{R}^{n\times n}\langle{\Der,\Int,\E}\rangle$ satisfies $(\Der+A)\cdot\tilde{G}=1$ and $\E\cdot\tilde{G}=0$. Passing to $\mathcal{R}\langle{\Der,\Int,\E}\rangle^{n\times n}$ via Lemma~\ref{lem:matrixisomorphism}, the latter identity implies $\E\cdot\tilde{G}_{i,j}=0$ and the former implies $L\cdot\tilde{G}_{1,n}=1$ and $\tilde{G}_{i,n} = \Der^{i-1}\cdot\tilde{G}_{1,n}$ by Lemma~\ref{lem:companionsolution}.
 Hence, we also have $\E\cdot\Der^{i-1}\cdot\tilde{G}_{1,n}=0$ for $i \in \{1,\dots,n\}$.\par
 Finally, we verify that $\tilde{G}_{1,n}=G$.
 With $H:=Z \cdot \Int \cdot Z^{-1} \in \mathcal{R}\langle{\Der,\Int,\E}\rangle^{n\times n}$, we have $\tilde{G}_{1,n} = H_{1,n}-\sum_{i,j=1}^nZ_{1,i}((\E{Z})^{-1})_{i,j}\cdot\E\cdot{H}_{j,n}$. From the proof of Theorem~\ref{thm:varconstCommutative}, we obtain that $H_{j,n} = (-1)^{n+k}(\Der^{j-1}z_k)\cdot\Int\cdot\frac{W(z_1,\dots,z_{i-1},z_{i+1},\dots,z_n)}{w}$. Altogether, this yields $\tilde{G}_{1,n}=G$.
\end{proof}

\subsection{Taylor formula}

Usually, Taylor's theorem is only considered for sufficiently smooth functions.
In integro-differential rings, an analog of the Taylor formula with integral remainder term
\[
 f(x) = \sum_{k=0}^n\frac{f^{(k)}(x_0)}{k!}x^k + \int_{x_0}^x\frac{(x-t)^n}{n!}f^{(n+1)}(t)\,dt
\]
can be formulated.
While the formula arising from the identity of operators in Theorem~\ref{thm:Taylor} is more complicated, it is also valid if singularities are present.

We start by giving a first version of the Taylor formula where the remainder term is given as repeated integral.
It simply follows by iterating \eqref{eq:EvaluationDef} resp.~\eqref{eq:opFTC2}.
See also Corollary 2.2.1 in \cite{Przeworska-Rolewicz1988}.

\begin{lemma}\label{lem:TaylorFirst}
 Let $(\mathcal{R},\Der,\Int)$ be an integro-differential ring.
 In $\mathcal{R}\langle{\Der,\Int,\E}\rangle$ we have for any $n \in \mathbb{N}$ that
 \[
  1 = \sum_{i=0}^n\Int^i\cdot\E\cdot\Der^i + \Int^{n+1}\cdot\Der^{n+1}.
 \]
\end{lemma}
\begin{proof}
 For any $n \in \mathbb{N}$, we can use \eqref{eq:opFTC2} to rewrite the right hand side:
 \begin{multline*}
  \sum_{i=0}^n\Int^i\cdot\E\cdot\Der^i + \Int^{n+1}\cdot\Der^{n+1}
  = \sum_{i=0}^n\Int^i\cdot\E\cdot\Der^i + \Int^n\cdot(1-\E)\cdot\Der^n\\
  = \sum_{i=0}^{n-1}\Int^i\cdot\E\cdot\Der^i + \Int^n\cdot\Der^n.
 \end{multline*}
 Setting $n=0$ here gives $\E + \Int\cdot\Der = 1$.
 Hence, the claim follows by induction.
\end{proof}

Using \eqref{eq:opLinInt} and the related identity in \eqref{eq:opSpecial}, we can always write operators $\Int^i\cdot\E$ as $x_i\cdot\E$, where $x_i:=\Int^i1$ as in Theorem~\ref{thm:polynomials}.
By \eqref{eq:opRB} and \eqref{eq:opSpecialRB}, we can always write $\Int^{n+1}$ without higher powers of $\Int$.
To make the resulting expressions simpler, we restrict to the case that $\E$ is multiplicative on polynomials, i.e.\ $\E x_mx_n=0$ for all $m,n\ge1$.

\begin{lemma}\label{lem:RepeatedIntegralOperator}
 Let $(\mathcal{R},\Der,\Int)$ be an integro-differential ring such that the induced evaluation $\E$ satisfies $\E x_mx_n=0$ for all $m,n\ge1$.
 Then, in $\mathcal{R}\langle{\Der,\Int,\E}\rangle$, we have for any $n \in \mathbb{N}$ that
 \[
  \Int^{n+1} = \sum_{k=0}^n(-1)^{n-k}{x_k}\cdot\Int\cdot{x_{n-k}}
  -\sum_{k=0}^{n-1}\sum_{j=1}^{n-k}(-1)^{n-k-j}{x_k}\cdot\E\cdot{x_j}\cdot\Int\cdot{x_{n-k-j}}.
 \]
\end{lemma}
\begin{proof}
 We prove this identity by induction on $n \in \mathbb{N}$. For $n=0$, the right hand side directly yields $\Int$ in agreement with the left hand side. Assuming the identity holds for some $n \in \mathbb{N}$, we multiply both sides by $\Int$ from the left and we rewrite the right hand side using \eqref{eq:opRB} and \eqref{eq:opLinInt} to obtain
 \begin{multline*}
  \Int^{n+2} = \sum_{k=0}^n(-1)^{n-k}\Big({\Int x_k}\cdot\Int-\Int\cdot{\Int x_k}-\E\cdot{\Int x_k}\cdot\Int\Big)\cdot{x_{n-k}}\\
  -\sum_{k=0}^{n-1}\sum_{j=1}^{n-k}(-1)^{n-k-j}{\Int x_k}\cdot\E\cdot{x_j}\cdot\Int\cdot{x_{n-k-j}}.
 \end{multline*}
  We have $x_{k+1}x_{n-k} = \binom{n+1}{k+1}x_{n+1}$ by Theorem~\ref{thm:polynomials}, so we can expand the right hand side into
 \begin{multline*}
  \sum_{k=0}^n(-1)^{n-k}{x_{k+1}}\cdot\Int\cdot{x_{n-k}}
  +\sum_{k=0}^n(-1)^{n-k+1}\binom{n+1}{k+1}\Int\cdot{x_{n+1}}\\
  -\sum_{k=0}^n(-1)^{n-k}\E\cdot{x_{k+1}}\cdot\Int\cdot{x_{n-k}}
  -\sum_{k=0}^{n-1}\sum_{j=1}^{n-k}(-1)^{n-k-j}{x_{k+1}}\cdot\E\cdot{x_j}\cdot\Int\cdot{x_{n-k-j}}.
 \end{multline*}
 Exploiting $\sum_{k=0}^n(-1)^{n-k+1}\binom{n+1}{k+1}=(-1)^{n+1}$ in the second sum, we can regroup terms to obtain the right hand side of the claimed identity for $n+1$.
 This completes the induction.
\end{proof}

Altogether, we obtain the following identity in $\mathcal{R}\langle{\Der,\Int,\E}\rangle$ generalizing the usual Taylor formula with integral remainder term.
This identity holds over any integro-differential ring in which the induced evaluation is multiplicative on the integro-differential subring generated by $1$.

\begin{theorem}[Taylor formula]
\label{thm:Taylor}
 Let $(\mathcal{R},\Der,\Int)$ be an integro-differential ring such that $\E$ is multiplicative on the integro-differential subring generated by $1$, i.e.\ $\E x_mx_n=0$ for all $m,n\ge1$.
 Then, for all $f \in \mathcal{R}$ and all $n \in \mathbb{N}$ we have
 \begin{multline*}
  1 = \sum_{k=0}^nx_k\cdot\E\cdot\Der^k + \sum_{k=0}^n(-1)^{n-k}x_k\cdot\Int\cdot{x_{n-k}}\cdot\Der^{n+1}\\
  -\sum_{k=0}^{n-1}\sum_{j=1}^{n-k}(-1)^{n-k-j}x_k\cdot\E\cdot{x_j}\cdot\Int\cdot{x_{n-k-j}}\cdot\Der^{n+1}
 \end{multline*}
\end{theorem}
\begin{proof}
 Follows from Lemma~\ref{lem:TaylorFirst} using $\Int^i\cdot\E=x_i\cdot\E$ and Lemma~\ref{lem:RepeatedIntegralOperator}, as explained above.
\end{proof}

In particular, if in addition to the assumptions of the theorem we have $\mathbb{Q}\subseteq\mathcal{R}$, then, with operators acting on some $f$, we have that
\begin{multline*}
 f = \sum_{k=0}^n\frac{x_1^k}{k!}\E\Der^kf + \sum_{k=0}^n\frac{(-1)^{n-k}}{k!(n-k)!}x_1^k\Int{x_1^{n-k}}\Der^{n+1}f\\
 -\sum_{k=0}^{n-1}\sum_{j=1}^{n-k}\frac{(-1)^{n-k-j}}{k!j!(n-k-j)!}x_1^k\E{x_1^j}\Int{x_1^{n-k-j}}\Der^{n+1}f
\end{multline*}
While the first and the second sum correspond to the Taylor polynomial and the integral remainder term, the third sum corresponds to an additional polynomial that arises from our general setting allowing non-multiplicative evaluations.
It can be viewed as an integro-differential algebraic version of the analytic formula
\[
 -\sum_{k=0}^{n-1}\frac{(x-x_0)^k}{k!}\left[\int_{x_0}^x\frac{(x-t)^{n-k}-(x_0-t)^{n-k}}{(n-k)!}f^{(n+1)}(t)dt\right]_{x=x_0},
\]
which vanishes for smooth functions $f(x)$.
Although in practice these additional terms yield zero even for many elements $f$ that model singular functions in concrete integro-differential rings, they cannot be dropped in general, as illustrated in $\mathcal{C}[x,x^{-1},\ln(x)]$ by considering $n=1$ and $f=\ln(x)$, for example.
With integration defined as in Example~\ref{ex:Laurent}, we have $x_1=x$ and with $f=\ln(x)$ the Taylor polynomial $\E{f}+x\E\Der{f}$ vanishes.
By $\Der^2f=-\frac{1}{x^2}$, the second sum yields $-\Int{x\Der^2f}+x\Int\Der^2f=\ln(x)+1$ and the third sum $-\E{x\Int\Der^2f}=-1$ compensates the constant term.
More generally, we can characterize the integro-differential rings where the additional polynomial does not play a role in the Taylor formula.
\begin{corollary}
 Let $(\mathcal{R},\Der,\Int)$ be an integro-differential ring such that $\E$ is multiplicative on the integro-differential subring generated by $1$, i.e.\ $\E x_mx_n=0$ for all $m,n\ge1$.
 Then, we have $\E x_ng=0$ for all $n\ge1$ and $g\in\mathcal{R}$ if and only if we have
 \[
  f = \sum_{k=0}^nx_k\E\Der^kf + \sum_{k=0}^n(-1)^{n-k}x_k\Int{x_{n-k}}\Der^{n+1}f
 \]
 for all $n\in\mathbb{N}$ and $f\in\mathcal{R}$.
\end{corollary}
\begin{proof}
 If $\E x_ng=0$ for all $n\ge1$ and $g\in\mathcal{R}$, then we have in particular $\E{x_j}\Int{x_{n-k-j}}\Der^{n+1}f=0$ for all $j,k,n\in\mathbb{N}$ and $f\in\mathcal{R}$ with $1 \le j \le n-k$.
 Hence, Theorem~\ref{thm:Taylor} implies the claimed identity for all $n\in\mathbb{N}$ and $f\in\mathcal{R}$.
\par
 For the converse, let $n\ge1$ be minimal such that $\E x_ng\neq0$ for some $g\in\mathcal{R}$.
 With such $g$, we let $f:=\Int^ng$ and, by minimality of $n$, we obtain that
 \[
  \sum_{k=0}^{n-1}\sum_{j=1}^{n-k}(-1)^{n-k-j}x_k\E{x_j}\Int{x_{n-k-j}}\Der^{n+1}f = \E{x_n}\Int\Der{g} = \E{x_n}g - \E{x_n}\E{g} = \E{x_n}g
 \]
 is nonzero.
 Consequently, Theorem~\ref{thm:Taylor} implies that the claimed identity does not hold for this $n$ and $f$.
\end{proof}

\subsubsection*{Acknowledgements}
This work was supported by the Austrian Science Fund (FWF): P~27229, P~31952, and P~32301.
Part of this work was done while both authors were at the Radon Institute of Computational and Applied Mathematics (RICAM) of the Austrian Academy of Sciences.
The authors would like to thank Alban Quadrat for bringing the book \cite{Przeworska-Rolewicz1988} to the attention of the second author during his stay at INRIA Saclay.

\appendix

\section{Normal forms for IDO in tensor rings}

The goal of this appendix is to state and prove a refinement of Theorem~\ref{thm:IDO} providing uniqueness of normal forms.
Uniqueness is achieved by representing operators by elements of a tensor ring, which is formed on a module of basic operators generated by $\mathcal{R}$, $\Der$, $\Int$, $\E$.
The ring of operators can be constructed as quotient of the tensor ring, where relations of basic operators are encoded by tensor reduction rules.
The main technical tool for proving uniqueness of normal forms is a generalization of Bergman's Diamond Lemma in tensor rings \cite{Bergman1978}.
For the convenience of the reader, we give a formal and largely self-contained summary of tensor reduction systems and we explain the translation of identities of operators into this framework.
In this appendix, $\mathcal{K}$ denotes a ring (not necessarily commutative) with unit element.

In Section~\ref{sec:tensorrings}, we start by recalling basic properties of bimodules and tensor rings on them. For further details on tensor rings and proofs see, for example, \cite{Rowen1991,Cohn2003b}.
Then, we recall decompositions with specialization, which are used for defining reduction rules.
In Section~\ref{sec:TenReS}, we state the Diamond Lemma for tensor reduction systems with specialization from \cite{HosseinPoorRaabRegensburger2018} and provide a summary of the relevant notions.
In Section~\ref{sec:TensorIDO}, we construct an appropriate tensor ring along with a tensor reduction system for dealing with IDO.
We use these to state Theorem~\ref{thm:IDOtensor}, which provides a precise formulation of uniqueness of normal forms of IDO in terms of tensors.
In Section~\ref{sec:TensorIDOFunctionals}, we give a complete proof of the theorem, where the necessary computations for verifying uniqueness of normal forms in the tensor ring are done in an automated way by our package \texttt{TenReS} in the computer algebra system Mathematica.
These computations are contained in the Mathematica file accompanying this paper, which includes a log of the reduction steps and is available at \url{http://gregensburger.com/softw/tenres/}.
The theorem proved in that section even covers more general rings of operators, which allow to deal with additional functionals besides the induced evaluation.

\subsection{Tensor rings on bimodules and decompositions}
\label{sec:tensorrings}

A $\mathcal{K}$-bimodule is a left $\mathcal{K}$-module $M$ which is also a right $\mathcal{K}$-module satisfying the associativity condition $(km)l = k(ml)$ for all $m \in M$ and $k,l \in \mathcal{K}$.
A ring $\mathcal{R}$ that is a $\mathcal{K}$-bimodule such that $(xy)z = x(yz)$ for any $x,y,z$ in $\mathcal{R}$ or $\mathcal{K}$ is called a $\mathcal{K}$-ring.
In particular, if $\mathcal{K}$ is a subring of some ring $\mathcal{R}$, then $\mathcal{R}$ is a $\mathcal{K}$-ring.

We first recall basic properties of the tensor product on $\mathcal{K}$-bimodules.
For $\mathcal{K}$-bimodules $M,N$, their $\mathcal{K}$-tensor product $M\otimes N$ is a $\mathcal{K}$-bimodule generated by the pure tensors $\{m\otimes n \mid m \in M, n \in N\}$ with relations
\begin{gather*}
 (m+\tilde{m})\otimes n = m\otimes n+ \tilde{m}\otimes n,\quad
 m\otimes(n+\tilde{n}) = m\otimes n+ m\otimes \tilde{n},\quad\text{and}\\
 mk\otimes n = m\otimes kn
\end{gather*}
having scalar multiplications
\[
 k(m\otimes n) = (km)\otimes n \quad\text{and}\quad (m\otimes n)k = m\otimes (nk)
\]
for all $m,\tilde{m} \in M$, $n,\tilde{n} \in N$, and $k \in \mathcal{K}$.

We denote the tensor product of $M$ with itself over $\mathcal{K}$ by $M^{\otimes{n}}=M\otimes\dots\otimes M$ ($n$ factors).
In particular, $M^{\otimes1}=M$ and we interpret $M^{\otimes0}$ as the $\mathcal{K}$-module $\mathcal{K}\varepsilon$, where $\varepsilon$ denotes the empty tensor and right scalar multiplication satisfies $(k_1\varepsilon)k_2=(k_1k_2)\varepsilon$ for $k_1,k_2\in\mathcal{K}$.
As a $\mathcal{K}$-bimodule, the tensor ring $\mathcal{K}\langle{M}\rangle$ is defined as the direct sum
\[
 \mathcal{K}\langle{M}\rangle = \bigoplus_{n=0}^\infty M^{\otimes{n}}.
\]
It can be turned into a $\mathcal{K}$-ring with unit element $\varepsilon$ where multiplication $M^{\otimes{r}} \times M^{\otimes{s}} \to M^{\otimes(r+s)}$ is defined via
\[
 (m_1\otimes\dots\otimes m_r,\tilde{m}_1\otimes\dots\otimes\tilde{m}_s)\mapsto m_1\otimes\dots\otimes m_r \otimes\tilde{m}_1\otimes\dots\otimes\tilde{m}_s.
\]

Via the tensor product, any decomposition of the module $M$ carries over to a decomposition of the tensor ring $\mathcal{K}\langle{M}\rangle$.
In particular, we use \emph{decompositions with specialization}, which were introduced in \cite{HosseinPoorRaabRegensburger2018}.
These are given by a family $(M_z)_{z \in Z}$ of $\mathcal{K}$-subbimodules of $M$ and a subset $X \subseteq Z$ with $M = \sum_{z \in Z}M_z = \bigoplus_{x \in X}M_x$ such that every module $M_z$, $z \in Z$, satisfies
\[
 M_z = \bigoplus_{x \in S(z)}M_x
\]
where $S(z) := \{x \in X \mid M_x \subseteq M_z\}$ is the set of \emph{specializations} of $z$.
Note that this definition implies $S(x)=\{x\}$ for $x \in X$.
For words $W=w_1\dots{w_n}$ in the word monoid $\langle{Z}\rangle$, we define the corresponding $\mathcal{K}$-subbimodule of $\mathcal{K}\langle{M}\rangle$ by
\[
 M_W := M_{w_1} \otimes\dots\otimes M_{w_n}.
\]
The notion of specialization extends from the alphabet $Z$ to the whole word monoid $\langle{Z}\rangle$ by
\[
 S(W) := \{v_1\dots{v_n}\in\langle{X}\rangle \mid \forall{i}:v_i\in S(w_i)\}
\]
such that $S(W) = \{V\in\langle{X}\rangle \mid M_V\subseteq M_W\}$.
We have the following generalization
\[
M_W = \bigoplus_{V \in S(W)} M_V
\]
of the direct sum above and the decomposition
\[
 \mathcal{K}\langle{M}\rangle = \sum_{W\in\langle{Z}\rangle}M_W = \bigoplus_{W\in\langle{X}\rangle}M_W
\]
of the tensor ring.

\subsection{Tensor reduction systems with specialization}
\label{sec:TenReS}

Fixing a decomposition with specialization of $M$, a \emph{reduction rule} for $\mathcal{K}\langle{M}\rangle$ is given by a pair $r=(W,h)$ of a word $W\in\langle{Z}\rangle$ and a $\mathcal{K}$-bimodule homomorphism $h\colon M_W \to \mathcal{K}\langle{M}\rangle$.
It acts on tensors of the form $a\otimes w\otimes b$ with $a\in M_A$, $w\in M_W$, and $b\in M_B$ for some $A,B\in\langle{Z}\rangle$ by $a\otimes w\otimes b \rightarrow_r a\otimes h(w)\otimes b$.
Later, we will specify homomorphisms $h$ in concrete reduction rules $(W,h)$ via their values on a generating set of $M_W$.
Formally, well-definedness of such homomorphisms can be ensured by the universal property of the tensor product.
A set $\Sigma$ of such reduction rules is called a reduction system over $Z$ on $\mathcal{K}\langle{M}\rangle$ and induces the two-sided \emph{reduction ideal}
\[
 I_\Sigma := (t-h(t)\mid(W,h)\in\Sigma\text{ and }t\in M_{W}) \subseteq \mathcal{K}\langle{M}\rangle.
\]
For computing in the factor ring $\mathcal{K}\langle{M}\rangle/I_{\Sigma}$, we apply the reduction relation $\rightarrow_\Sigma$ induced by $\Sigma$ on $\mathcal{K}\langle{M}\rangle$.
It reduces a tensor $t\in\mathcal{K}\langle{M}\rangle$ to a tensor $s\in\mathcal{K}\langle{M}\rangle$ if there is an $r\in\Sigma$ such that $t\rightarrow_rs$.
We say that $t$ \emph{can be reduced} to $s$ by $\Sigma$ if $t=s$ or there exists a finite sequence of reduction rules $r_1,\ldots,r_n$ in $\Sigma$ such that
\[
 t \rightarrow_{r_1} t_1\rightarrow_{r_2} \cdots \rightarrow_{r_{n-1}} t_{n-1} \rightarrow_{r_n} s.
\]
If one tensor can be reduced to another, then their difference is contained in $I_\Sigma$ and they represent the same element of $\mathcal{K}\langle{M}\rangle/I_{\Sigma}$.
The $\mathcal{K}$-subbimodule of \emph{irreducible tensors}
\[
 \mathcal{K}\langle{M}\rangle_\mathrm{irr} = \bigoplus\limits_{W\in\langle{X}\rangle_\mathrm{irr}} M_W
\]
can be characterized by the set of \emph{irreducible words} $\langle{X}\rangle_\mathrm{irr}\subseteq\langle{X}\rangle$, which consists of those words that avoid subwords arising as specializations $S(W)$ of words occurring in reduction rules $(W,h) \in \Sigma$.
The irreducible tensors to which a given tensor $t$ can be reduced, are called its normal forms. If $t$ has a unique normal form, it is denoted by $t\!\downarrow_\Sigma$.

An \emph{ambiguity} is a minimal situation where two (not necessarily distinct) reduction rules can be applied to tensors in different ways.
For each pure tensor of this kind, the corresponding \emph{S-polynomial} is the difference of the results of the two reduction steps.
For example, an overlap ambiguity arises from two reduction rules $(AB_1,h),(B_2C,\tilde{h})\in\Sigma$, where $A,B_1,B_2,C\in\langle{Z}\rangle$ are nonempty such that $B_1,B_2$ are equal or have a common specialization, and corresponding S-polynomials are referred to by $\Spol(A\underline{B_1},\underline{B_2}C)$.
An ambiguity is called \emph{resolvable}, if all its S-polynomials can be reduced to zero by $\Sigma$.
If all ambiguities of $\Sigma$ are resolvable, then the reduction relation induced by $\Sigma$ on $\mathcal{K}\langle{M}\rangle$ is confluent and, by abuse of language, we also call $\Sigma$ \emph{confluent}.
This means that there are no hidden consequences implied in $\mathcal{K}\langle{M}\rangle/I_{\Sigma}$ by the identities explicitly specified by $\Sigma$.

The following theorem relies on the existence of a partial order of words in $\langle{Z}\rangle$ that has certain properties, which are briefly explained now.
A partial order $\le$ on $\langle{Z}\rangle$ is called a \emph{semigroup partial order} if it is compatible with concatenation of words.
If in addition the empty word $\epsilon$ is the least element of $\langle{Z}\rangle$, then $\le$ is called a \emph{monoid partial order}.
It is called \emph{Noetherian} if there are no infinite descending chains.
We call a partial order $\le$ on $\langle{Z}\rangle$ \emph{consistent with specialization} if every strict inequality $V < W$ implies $\tilde{V} < \tilde{W}$ for all specializations $\tilde{V} \in S(V)$ and $\tilde{W} \in S(W)$.
A partial order $\le$ on $\langle{Z}\rangle$ is \emph{compatible} with a reduction system $\Sigma$ over $Z$ on $\mathcal{K}\langle{M}\rangle$ if for all $(W,h) \in \Sigma$ the image of $h$ is contained in the sum of modules $M_V$ where $V \in \langle{Z}\rangle$ satisfies $V<W$.

\begin{theorem}\label{thm:DiamondLemma}
 \cite[Thm.~20]{HosseinPoorRaabRegensburger2018}
 Let $M$ be a $\mathcal{K}$-bimodule and let $(M_z)_{z\in Z}$ be a decomposition with specialization.
 Let $\Sigma$ be a reduction system over $Z$ on $\mathcal{K}\langle{M}\rangle$ and let $\le$ be a Noetherian semigroup partial order on $\langle{Z}\rangle$ consistent with specialization and compatible with $\Sigma$.
 Then, the following are equivalent:
 \begin{enumerate}
  \item All ambiguities of $\Sigma$ are resolvable.
  \item Every $t\in\mathcal{K}\langle{M}\rangle$ has a unique normal form $t\!\downarrow_\Sigma$.
  \item $\mathcal{K}\langle{M}\rangle/I_\Sigma$ and $\mathcal{K}\langle{M}\rangle_\mathrm{irr}$ are isomorphic as $\mathcal{K}$-bimodules.
 \end{enumerate}
 Moreover, if these conditions are satisfied, then we can define a multiplication on $\mathcal{K}\langle{M}\rangle_\mathrm{irr}$ by $s \cdot t := (s \otimes t) \! \downarrow_\Sigma$ so that $\mathcal{K}\langle{M}\rangle/I_\Sigma$ and $\mathcal{K}\langle{M}\rangle_\mathrm{irr}$ are isomorphic as $\mathcal{K}$-rings.
\end{theorem}

\subsection{Tensor reduction systems for IDO}
\label{sec:TensorIDO}

Before using the tensor setting to construct the ring of integro-differential operators $\mathcal{R}\langle{\Der,\Int,\E}\rangle$, we illustrate this construction on the well-known ring of differential operators $\mathcal{R}\langle{\Der}\rangle$ to highlight some of the special properties of the construction.
Usually, the ring of differential operators with coefficients from $\mathcal{R}$ is constructed via skew polynomials $\sum_i f_i\Der^i$ over $\mathcal{R}$ in one indeterminate $\Der$, with commutation rule $\Der \cdot f = f \cdot \Der+\Der{f}$.
Let $(\mathcal{R},\Der,\Int)$ be an integro-differential ring and let $\mathcal{C}$ denote its ring of constants.
In the following, we consider $\mathcal{K}$-tensor rings with $\mathcal{K}:=\mathcal{C}$.

\begin{example}
The module of basic operators that generates all differential operators is given by
\[
 M := M_\mathsf{R} \oplus M_\mathsf{D},
\]
with $\mathcal{K}$-bimodules
\begin{equation}
\label{eq:module_R_full}
 M_\mathsf{R} := \mathcal{R}
\end{equation}
and $M_\mathsf{D}$ defined as the free left $\mathcal{K}$-module
\begin{equation}
\label{eq:module_D}
 M_\mathsf{D} := \mathcal{K}\Der
\end{equation}
generated by the symbol $\Der$, which we view as a $\mathcal{K}$-bimodule with the right multiplication $c\Der\cdot d = cd\Der$ for all $c,d \in \mathcal{K}$. This definition is based on left $\mathcal{K}$-linearity of the derivation $\Der$ on $\mathcal{R}$.
The commutation rule $\Der \cdot f = f \cdot \Der+\Der{f}$ coming from the Leibniz rule in $\mathcal{R}$ translates to the tensor reduction rule
\[
 (\mathsf{DR},\Der{\otimes}{f} \mapsto f{\otimes}\Der+\Der{f}).
\]
This rule is formalized by the $\mathcal{K}$-bimodule homomorphism $M_\mathsf{DR}=M_\mathsf{D}\otimes{M_\mathsf{R}} \to \mathcal{K}\langle{M}\rangle$ defined by $\Der{\otimes}{f} \mapsto f{\otimes}\Der+\Der{f}$ on tensors $\Der{\otimes}{f}$ generating $M_\mathsf{DR}$.
Note that this homomorphism represents a parameterized family of identities.
Reduction of tensors by the rule above allows to syntactically move all multiplication operators to the left of any differentiation operator.

Note that, in order to correctly model differential operators as equivalence classes of tensors in $\mathcal{K}\langle{M}\rangle$, other relations among operators need to be phrased as tensor reduction rules as well. This is because the tensor ring $\mathcal{K}\langle{M}\rangle$ itself is constructed without respecting relations coming from multiplication in $\mathcal{R}$.
For instance, the composition of two multiplication operators $f$ and $g$ is a multiplication operator again, which leads to the reduction rule $(\mathsf{RR},f{\otimes}g \mapsto fg)$ defined on the module $M_\mathsf{RR}=M_\mathsf{R}\otimes{M_\mathsf{R}}$.
Moreover, the multiplication operator that multiplies by $1$ acts like the identity operator, which is represented by the empty tensor $\varepsilon$. To define reduction rules that act only on $\mathcal{C}$ instead of all of $\mathcal{R}$, we need a direct decomposition of the $\mathcal{K}$-bimodule
\begin{equation}
\label{eq:module_R}
 M_\mathsf{R} = M_\mathsf{K} \oplus M_{\tilde{\mathsf{R}}},
\end{equation}
which in our case can be given by
\begin{equation}
\label{eq:special_module_R}
 M_\mathsf{K} := \mathcal{K} \quad\text{and}\quad M_{\tilde{\mathsf{R}}} := \Int\mathcal{R}
\end{equation}
based on Lemma~\ref{lem:decompositionofring}. Then, we can define a reduction rule on $M_\mathsf{K}=\mathcal{C}$ by $1\mapsto\varepsilon$. Altogether, the tensor reduction system $\Sigma_\mathrm{Diff}=\{r_\mathsf{K},r_\mathsf{RR},r_\mathsf{DR}\}$ for differential operators is given by the three reduction rules
\begin{equation*}\label{eq:SigmaDiff}
 \{
(\mathsf{K},{1}\mapsto \varepsilon), \quad (\mathsf{RR}, {f}{\otimes}{g}\mapsto {fg}), \quad
   (\mathsf{DR},\Der{\otimes}{f} \mapsto f{\otimes}\Der+\Der{f})
 \}
\end{equation*}
defined above. It induces the two-sided ideal
\begin{align*}
 I_\mathrm{Diff} :=&\ (t-h(t)\mid (W,h)\in \Sigma_\mathrm{Diff} \text{ and }t\in M_W)\\
 =&\ (1-\varepsilon, f{\otimes}g-fg, \Der{\otimes}{f}-f{\otimes}\Der-\Der{f} \mid f,g \in \mathcal{R})
\end{align*}
in the ring $\mathcal{K}\langle{M}\rangle$ and computations with differential operators are modelled in the quotient ring
\[
 \mathcal{R}\langle{\Der}\rangle:=\mathcal{K}\langle{M}\rangle/I_\mathrm{Diff}.
\]
Tensors that are not reducible w.r.t.\ $\Sigma_\mathrm{Diff}$ are precisely $\mathcal{K}$-linear combinations of pure tensors of the form $\Der^{\otimes{i}}$ and $f\otimes\Der^{\otimes{i}}$, where $i\in\mathbb{N}_0$ and $f\in\Int\mathcal{R}$.
With alphabets $X=\{\mathsf{K},\tilde{\mathsf{R}},\mathsf{D}\}$ and $Z=X\cup\{\mathsf{R}\}$, one can check that all ambiguities of $\Sigma_\mathrm{Diff}$ are resolvable.
Using an appropriate ordering of words, one can show that the conditions of Theorem~\ref{thm:DiamondLemma} are satisfied by this construction of $\mathcal{R}\langle{\Der}\rangle$.
\qed
\end{example}

Proceeding to an analogous construction of the ring $\mathcal{R}\langle{\Der,\Int,\E}\rangle$ from Definition~\ref{def:OperatorRing}, we also require tensor reduction rules corresponding to the identities \eqref{eq:opFTC1}--\eqref{eq:opLinEval}, in addition to the three rules from the example above.
To this end, analogous to $M_\mathsf{D}$ above, we introduce the $\mathcal{K}$-bimodules
\begin{equation}
\label{eq:module_IandE}
 M_\mathsf{I}:=\mathcal{K}\Int \quad\text{and}\quad M_\mathsf{E}:=\mathcal{K}\E,
\end{equation}
which are freely generated as left $\mathcal{K}$-modules by the symbols $\Int$ and $\E$, respectively.
Then, we consider the $\mathcal{K}$-tensor ring on the $\mathcal{K}$-bimodule
\[
 M:=M_\mathsf{R} \oplus M_\mathsf{D} \oplus M_\mathsf{I} \oplus M_\mathsf{E}.
\]
Altogether, we obtain the tensor reduction system given in Table~\ref{tab:SigmaIDOdef}. The two-sided ideal $I_\mathrm{IDO}$ induced by it allows to construct the ring $\mathcal{R}\langle{\Der,\Int,\E}\rangle$ as the quotient $\mathcal{K}\langle{M}\rangle/I_\mathrm{IDO}$.
\begin{table}[h]
\begin{center}
\begin{tabular}{|cr@{\ $\mapsto$\ }l|cr@{\ $\mapsto$\ }l|}
\hline
   $\mathsf{K}$ & $1$&$\varepsilon$ & $\mathsf{ID}$ & $\Int{\otimes}\Der$&$\varepsilon-\E$\\
   $\mathsf{RR}$ & $f{\otimes}g$&$fg$ & $\mathsf{DRE}$ & $\Der{\otimes}f{\otimes}\E$&$\Der{f}{\otimes}\E$\\
   $\mathsf{DR}$ & $\Der{\otimes}f$&$f{\otimes}\Der+\Der{f}$ & $\mathsf{IRE}$ & $\Int{\otimes}f{\otimes}\E$&$\Int{f}{\otimes}\E$\\
   $\mathsf{DI}$ & $\Der{\otimes}\Int$&$\varepsilon$ & $\mathsf{ERE}$ & $\E{\otimes}f{\otimes}\E$&$(\E{f})\E$
\\\hline
\end{tabular}
\caption{Defining reduction system for integro-differential operators}
\label{tab:SigmaIDOdef}
\end{center}
\vspace*{-\bigskipamount}
\end{table}
However, this reduction system is not confluent. In order to obtain normal forms that are unique as tensors in $\mathcal{K}\langle{M}\rangle$, we need a confluent tensor reduction system on $\mathcal{K}\langle{M}\rangle$ that induces the same ideal $I_\mathrm{IDO}$.
The confluent tensor reduction system given in Table~\ref{tab:SigmaIDOcomplete} can be obtained by turning the identities of Table~\ref{tab:op} into tensor reduction rules and including the rules $r_\mathsf{K}$ and $r_\mathsf{RR}$ from above.
\begin{table}[h]
\begin{center}
\begin{tabular}{|cr@{\ $\mapsto$\ }l|cr@{\ $\mapsto$\ }l|}
\hline
   $\mathsf{K}$ & $1$&$\varepsilon$ & $\mathsf{EI}$ & $\E{\otimes}\Int$&$0$\\
   $\mathsf{RR}$ & $f{\otimes}g$&$fg$ & $\mathsf{IRD}$ & $\Int{\otimes}f{\otimes}\Der$&$f-\E{\otimes}f-\Int{\otimes}\Der{f}$\\
   $\mathsf{DR}$ & $\Der{\otimes}f$&$f{\otimes}\Der+\Der{f}$ & $\mathsf{IRE}$ & $\Int{\otimes}f{\otimes}\E$&$\Int{f}{\otimes}\E$\\
   $\mathsf{DE}$ & $\Der{\otimes}\E$&$0$ & $\mathsf{IRI}$ & $\Int{\otimes}f{\otimes}\Int$&$\Int{f}{\otimes}\Int-\E{\otimes}\Int{f}{\otimes}\Int-\Int{\otimes}\Int{f}$\\
   $\mathsf{DI}$ & $\Der{\otimes}\Int$&$\varepsilon$ & $\mathsf{ID}$ & $\Int{\otimes}\Der$&$\varepsilon-\E$\\
   $\mathsf{ERE}$ & $\E{\otimes}f{\otimes}\E$&$(\E{f})\E$ & $\mathsf{IE}$ & $\Int{\otimes}\E$&$\Int1{\otimes}\E$\\
   $\mathsf{EE}$ & $\E{\otimes}\E$&$\E$ & $\mathsf{II}$ & $\Int{\otimes}\Int$&$\Int1{\otimes}\Int-\E{\otimes}\Int1{\otimes}\Int-\Int{\otimes}\Int1$
\\\hline
\end{tabular}
\caption{Confluent reduction system $\Sigma_\mathrm{IDO}$ for integro-differential operators}
\label{tab:SigmaIDOcomplete}
\end{center}
\vspace*{-\bigskipamount}
\end{table}
This allows us to state the following more precise version of Theorem~\ref{thm:IDO}.
Note that multiples like $f\cdot\Der$ are treated differently now, due to the fact that we are working in the tensor ring $\mathcal{K}\langle{M}\rangle$ and we have the reduction rule $(\mathsf{K},1\mapsto\varepsilon)$, which splits $f\in\mathcal{R}$ according to Lemma~\ref{lem:decompositionofring}.

\begin{theorem}
\label{thm:IDOtensor}
 Let $(\mathcal{R},\Der,\Int)$ be an integro-differential ring with constants $\mathcal{K}=\mathcal{C}$. Let $M$ be given as above in terms of the modules defined in Eqs.\ \eqref{eq:module_R}, \eqref{eq:special_module_R}, \eqref{eq:module_D}, and \eqref{eq:module_IandE} and let the tensor reduction system $\Sigma_\mathrm{IDO}$ be defined by Table~\ref{tab:SigmaIDOcomplete}.\par
 Then every $t\in \mathcal{K}\langle{M}\rangle$ has a unique normal form $t\! \downarrow_{\Sigma_\mathrm{IDO}} \in \mathcal{K}\langle{M}\rangle$, which can be written as a $\mathcal{K}$-linear combination of pure tensors of the form
 \[
  f \otimes \Der^{\otimes j}, \quad f \otimes \Int \otimes g, \quad   f \otimes \E \otimes g \otimes \Der^{\otimes j}, \quad \text{or} \quad  f \otimes \E \otimes h \otimes \Int \otimes g
 \]
 where $j\in\mathbb{N}_{0}$, $f,g,h\in\Int\mathcal{R}$, and each $f$ and $g$ may be absent. Moreover,
 \[
  \mathcal{R}\langle{\Der,\Int,\E}\rangle \cong \mathcal{K}\langle{M}\rangle_\mathrm{irr}
 \]
 as $\mathcal{K}$-rings, where multiplication on $\mathcal{K} \langle M \rangle_\mathrm{irr}$ is defined by $s \cdot t := (s \otimes t) \! \downarrow_{\Sigma_\mathrm{IDO}}$.
\end{theorem}

Instead of proving this theorem, we will prove a more general one below, which allows to include additional functionals into the construction of the ring of operators.
Theorem~\ref{thm:IDOtensor} follows from Theorem~\ref{thm:IDOPhi} by specializing $\Phi=\{\E\}$.

\subsection{Proof of normal forms for IDO with functionals}
\label{sec:TensorIDOFunctionals}

To treat more general problems than the initial value problems in Section~\ref{sec:IVP}, it is useful to include additional functionals into the ring of operators. For instance, dealing with boundary problems requires evaluations at more than one point.
In general, we consider a set $\Phi$ of $\mathcal{K}$-linear functionals $\mathcal{R}\to\mathcal{K}$ including $\E$. We consider the $\mathcal{K}$-bimodule $M_{\mathsf{\Phi}}$ defined as free left $\mathcal{K}$-module
\begin{equation}\label{eq:module_Phi}
 M_{\mathsf{\Phi}} := \mathcal{K}\Phi
\end{equation}
generated by the elements of $\Phi$, where we define right multiplication in terms of $\Phi$ by $\left(\sum_{\varphi\in\Phi}c_\varphi\varphi\right)\cdot d = \sum_{\varphi\in\Phi}c_\varphi d\varphi$ for $c_\varphi,d\in\mathcal{K}$.
Note that $\mathcal{K}$-linear combinations of $\mathcal{K}$-linear maps are not necessarily $\mathcal{K}$-linear again, if $\mathcal{K}$ is not commutative.
Since $\Der$, $\Int$, and all elements of $\Phi$ are $\mathcal{K}$-linear, we have, for instance, that $\varphi f\psi g = (\varphi f)\psi g$ for all $f,g\in\mathcal{R}$ and $\varphi,\psi\in\Phi$.
This allows to extend the last three reduction rules of Table~\ref{tab:SigmaIDOdef} to cover all elements of $\Phi$ instead of the evaluation $\E$ only.
Altogether, we consider the $\mathcal{K}$-tensor ring on the $\mathcal{K}$-bimodule
\begin{equation}\label{eq:M_IDOPhi}
 M := M_\mathsf{R} \oplus M_\mathsf{D} \oplus M_\mathsf{I} \oplus M_\mathsf{\Phi}
\end{equation}
and we have the defining reduction system given in Table~\ref{tab:SigmaIDOPhiDef}.
\begin{table}[h]
\begin{center}
\begin{tabular}{|cr@{\ }c@{\ }l|cr@{\ }c@{\ }l|}
\hline
   $\mathsf{K}$ & $1$&$\mapsto$&$\varepsilon$ & $\mathsf{ID}$ & $\Int{\otimes}\Der$&$\mapsto$&$\varepsilon-\E$\\
   $\mathsf{RR}$ & $f{\otimes}g$&$\mapsto$&$fg$ & $\mathsf{DR\Phi}$ & $\Der{\otimes}f{\otimes}\varphi$&$\mapsto$&$\Der{f}{\otimes}\varphi$\\
   $\mathsf{DR}$ & $\Der{\otimes}f$&$\mapsto$&$f{\otimes}\Der+\Der{f}$ & $\mathsf{IR\Phi}$ & $\Int{\otimes}f{\otimes}\varphi$&$\mapsto$&$\Int{f}{\otimes}\varphi$\\
   $\mathsf{DI}$ & $\Der{\otimes}\Int$&$\mapsto$&$\varepsilon$ & $\mathsf{\Phi R\Phi}$ & $\varphi{\otimes}f{\otimes}\psi$&$\mapsto$&$(\varphi{f})\psi$
\\\hline
\end{tabular}
\caption{Defining reduction system for integro-differential operators with functionals}
\label{tab:SigmaIDOPhiDef}
\end{center}
\vspace*{-\bigskipamount}
\end{table}
\begin{definition}\label{def:IDOPhi_tensor}
 Let $(\mathcal{R}, \Der, \Int)$ be an integro-differential ring with constants $\mathcal{K}=\mathcal{C}$ and let $\Phi$ be a set of $\mathcal{K}$-linear functionals $\mathcal{R}\to\mathcal{K}$ including $\E$. Let the $\mathcal{K}$-bimodule $M$ be defined as above in Eqs.\ \eqref{eq:module_R_full}, \eqref{eq:module_D}, \eqref{eq:module_IandE}, \eqref{eq:module_Phi}, and \eqref{eq:M_IDOPhi}. We call
 \[
  \mathcal{R}\langle{\Der,\Int,\Phi}\rangle:=\mathcal{K}\langle{M}\rangle /I_\mathrm{IDO\Phi}
 \]
 the \emph{ring of integro-differential operators with functionals $\Phi$}, where $I_\mathrm{IDO\Phi}$ is the two-sided ideal induced by the reduction system obtained from Table~\ref{tab:SigmaIDOPhiDef} using also submodules of $M$ defined in Eq.~\eqref{eq:special_module_R}.
\end{definition}
We use Theorem~\ref{thm:DiamondLemma} above to determine unique normal forms of tensors.
Again, the reduction system given by Table~\ref{tab:SigmaIDOPhiDef} is not confluent and we need to construct a confluent reduction system, like $\Sigma_\mathrm{IDO\Phi}$ given in Table~\ref{tab:SigmaIDOPhi}, by a completion process similar to how Table~\ref{tab:SigmaIDOcomplete} was obtained.
Observe that, whenever $\Phi=\{\E\}$, the ring $\mathcal{R}\langle{\Der,\Int,\Phi}\rangle$ and the relation $\rightarrow_{\Sigma_\mathrm{IDO\Phi}}$ specialize to $\mathcal{R}\langle{\Der,\Int,\E}\rangle$ and $\rightarrow_{\Sigma_\mathrm{IDO}}$, respectively.
\begin{table}[h]
\begin{center}
\begin{tabular}{|cr@{\ }c@{\ }l|cr@{\ }c@{\ }l|}
\hline
   $\mathsf{K}$ & $1$&$\mapsto$&$\varepsilon$ & $\mathsf{EI}$ & $\E{\otimes}\Int$&$\mapsto$&$0$\\
   $\mathsf{RR}$ & $f{\otimes}g$&$\mapsto$&$fg$ & $\mathsf{IRD}$ & $\Int{\otimes}f{\otimes}\Der$&$\mapsto$&$f-\E{\otimes}f-\Int{\otimes}\Der{f}$\\
   $\mathsf{DR}$ & $\Der{\otimes}f$&$\mapsto$&$f{\otimes}\Der+\Der{f}$ & $\mathsf{IR\Phi}$ & $\Int{\otimes}f{\otimes}\varphi$&$\mapsto$&$\Int{f}{\otimes}\varphi$\\
   $\mathsf{D\Phi}$ & $\Der{\otimes}\varphi$&$\mapsto$&$0$ & $\mathsf{IRI}$ & $\Int{\otimes}f{\otimes}\Int$&$\mapsto$&$\Int{f}{\otimes}\Int-\E{\otimes}\Int{f}{\otimes}\Int-\Int{\otimes}\Int{f}$\\
   $\mathsf{DI}$ & $\Der{\otimes}\Int$&$\mapsto$&$\varepsilon$ & $\mathsf{ID}$ & $\Int{\otimes}\Der$&$\mapsto$&$\varepsilon-\E$\\
   $\mathsf{\Phi R\Phi}$ & $\varphi{\otimes}f{\otimes}\psi$&$\mapsto$&$(\varphi{f})\psi$ & $\mathsf{I\Phi}$ & $\Int{\otimes}\varphi$&$\mapsto$&$\Int1{\otimes}\varphi$\\
   $\mathsf{\Phi\Phi}$ & $\varphi{\otimes}\psi$&$\mapsto$&$(\varphi1)\psi$ & $\mathsf{II}$ & $\Int{\otimes}\Int$&$\mapsto$&$\Int1{\otimes}\Int-\E{\otimes}\Int1{\otimes}\Int-\Int{\otimes}\Int1$
\\\hline
\end{tabular}
\caption{Confluent reduction system $\Sigma_\mathrm{IDO\Phi}$ for integro-differential operators with functionals}
\label{tab:SigmaIDOPhi}
\end{center}
\vspace*{-\bigskipamount}
\end{table}

\begin{theorem}
\label{thm:IDOPhi}
 Let $(\mathcal{R}, \Der, \Int)$ be an integro-differential ring with constants $\mathcal{K}=\mathcal{C}$ and let $\Phi$ be a set of $\mathcal{K}$-linear functionals $\mathcal{R}\to\mathcal{K}$ including $\E$. Let $M$ and $\mathcal{R}\langle{\Der,\Int,\Phi}\rangle$ be defined as in Definition~\ref{def:IDOPhi_tensor} above and let the tensor reduction system $\Sigma_\mathrm{IDO\Phi}$ be defined by Table~\ref{tab:SigmaIDOPhi}.\par
 Then every $t\in \mathcal{K}\langle{M}\rangle$ has a unique normal form $t\! \downarrow_{\Sigma_\mathrm{IDO\Phi}} \in \mathcal{K}\langle{M}\rangle$, which can be written as a $\mathcal{K}$-linear combination of pure tensors of the form
 \[
  f \otimes \Der^{\otimes j}, \quad f \otimes \Int \otimes g, \quad   f \otimes \varphi \otimes h \otimes \Der^{\otimes j}, \quad \text{or} \quad  f \otimes \varphi \otimes h \otimes \Int \otimes g
 \]
 where $j\in\mathbb{N}_{0}$, $f,g,h\in\Int\mathcal{R}$, $\varphi\in\Phi$, and each $f,g,h$ may be absent such that $\varphi \otimes h \otimes \Int$ does not specialize to $\E \otimes \Int$. Moreover,
 \[
  \mathcal{R}\langle{\Der,\Int,\Phi}\rangle \cong \mathcal{K}\langle{M}\rangle_\mathrm{irr}
 \]
 as $\mathcal{K}$-rings, where multiplication on $\mathcal{K} \langle M \rangle_\mathrm{irr}$ is defined by $s \cdot t := (s \otimes t) \! \downarrow_{\Sigma_\mathrm{IDO\Phi}}$.
\end{theorem}
\begin{proof}
 We use the alphabets $X:=\{\mathsf{K},\mathsf{\tilde{R}},\mathsf{D},\mathsf{I},\mathsf{E},\mathsf{\tilde{\Phi}}\}$ and $Z:=X\cup\{\mathsf{R},\mathsf{\Phi}\}$, which turns $(M_z)_{z\in Z}$ into a decomposition with specialization for the module $M$, where $S(\mathsf{R})=\{\mathsf{K},\tilde{\mathsf{R}}\}$ and $S(\mathsf{\Phi})=\{\mathsf{E},\tilde{\mathsf{\Phi}}\}$.
 For defining a Noetherian monoid partial order $\le$ on $\langle{Z}\rangle$ consistent with specialization that is compatible with $\Sigma_\mathrm{IDO\Phi}$, it is sufficient to require the order to satisfy
 \[
  \mathsf{DR}>\mathsf{RD} \quad\text{and}\quad \mathsf{I}>\mathsf{E\tilde{R}}.
 \]
 For instance, we could first define a monoid total order on $\langle{\{\mathsf{R},\mathsf{D},\mathsf{I},\mathsf{\Phi}\}}\rangle\subseteq\langle{Z}\rangle$ by counting occurrences of the letter $\mathsf{I}$ and breaking ties with any degree-lexicographic order satisfying $\mathsf{D}>\mathsf{R}$ and then generate from it a partial order on $\langle{Z}\rangle$ that is consistent with specialization.
 \par
 By our Mathematica package \texttt{TenReS}, we generate all ambiguities of $\Sigma_\mathrm{IDO\Phi}$ and verify that they are resolvable.
 There are 54 ambiguities and indeed all S-polynomials reduce to zero, see the accompanying Mathematica file at \url{http://gregensburger.com/softw/tenres/}.
 Here, we just give two short examples.
 The first one illustrates in its last step of computation that also identities in $M_\mathsf{R}$, like the Leibniz rule, need to be used.
 \begin{align*}
  \Spol(\mathsf{IR\underline{D}},\mathsf{\underline{D}R}) &= (f-\E\otimes{f}-\Int\otimes\Der{f})\otimes{g} - \Int\otimes{f}\otimes(g\otimes\Der+\Der{g})\\
  &\rightarrow_{r_\mathsf{RR}} fg-\E\otimes{fg}-\Int\otimes(\Der{f})g - \Int\otimes{fg}\otimes\Der - \Int\otimes{f\Der{g}}\\
  &\rightarrow_{r_\mathsf{IRD}} -\Int\otimes(\Der{f})g + \Int\otimes\Der{fg} - \Int\otimes{f\Der{g}}\\
  &= \Int\otimes(-(\Der{f})g+\Der{fg}-f\Der{g}) = 0
 \end{align*}
 The second one illustrates that ambiguities involving specialization also need to be considered.
 \[
  \Spol(\mathsf{\Phi\underline{\Phi}},\mathsf{\underline{E}I}) = (\varphi1)\E\otimes\Int - \varphi\otimes0 \rightarrow_{r_\mathsf{EI}} 0
 \]
 Since all ambiguities are resolvable, by Theorem~\ref{thm:DiamondLemma} every element in $\mathcal{K}\langle{M}\rangle$ has a unique normal form and $\mathcal{R}\langle{\Der,\Int,\Phi}\rangle \cong \mathcal{K}\langle{M}\rangle_\mathrm{irr}$ as $\mathcal{K}$-rings.
 \par
 It remains to determine the explicit form of elements in $\mathcal{K}\langle{M}\rangle_\mathrm{irr}$. In order to do so, we determine the set of irreducible words $\langle{X}\rangle_\mathrm{irr}$ in $\langle{X}\rangle$. Irreducible words containing only the letters $\mathsf{K},\tilde{\mathsf{R}},\mathsf{E},\tilde{\mathsf{\Phi}}$ have to avoid the subwords arising from the reduction rules $\mathsf{K}$, $S(\mathsf{RR})=\{\mathsf{KK},\mathsf{K\tilde{R}},\mathsf{\tilde{R}K},\mathsf{\tilde{R}\tilde{R}}\}$, $S(\mathsf{\Phi\Phi})=\{\mathsf{EE},\mathsf{E\tilde{\Phi}},\mathsf{\tilde{\Phi}E},\mathsf{\tilde{\Phi}\tilde{\Phi}}\}$, and $S(\mathsf{\Phi R\Phi})$.
 Hence they are given by
 \[
  \epsilon, \mathsf{\tilde{R}}, \mathsf{E}, \mathsf{\tilde{\Phi}}, \mathsf{\tilde{R}E}, \mathsf{\tilde{R}\tilde{\Phi}}, \mathsf{E\tilde{R}}, \mathsf{\tilde{\Phi}\tilde{R}}, \mathsf{\tilde{R}E\tilde{R}}, \mathsf{\tilde{R}\tilde{\Phi}\tilde{R}}
 \]
 Allowing also the letter $\mathsf{D}$, we have to avoid the subwords coming from $S(\mathsf{DR})=\{\mathsf{DK},\mathsf{D\tilde{R}}\}$ and $S(\mathsf{D\Phi})=\{\mathsf{DE},\mathsf{D\tilde{\Phi}}\}$. Therefore, we can only append words $\mathsf{D}^j$ with $j\in\mathbb{N}_0$ to the irreducible words determined so far, in order to obtain all elements of $\langle{X}\rangle_\mathrm{irr}$ not containing the letter $\mathsf{I}$.
 Finally, we also consider the letter $\mathsf{I}$. Since we have to avoid the subwords $S(\mathsf{I\Phi})=\{\mathsf{IE},\mathsf{I\tilde{\Phi}}\}$, $\mathsf{ID}$, and $\mathsf{II}$, any letter immediately following $\mathsf{I}$ has to be $\mathsf{\tilde{R}}$. In addition, we have to avoid the subwords $S(\mathsf{IR\Phi}) = \{\mathsf{IKE},\mathsf{IK\tilde{\Phi}},\mathsf{I\tilde{R}E},\mathsf{I\tilde{R}\tilde{\Phi}}\}$, $S(\mathsf{IRD})=\{\mathsf{IKD},\mathsf{I\tilde{R}D}\}$, and $S(\mathsf{IRI})=\{\mathsf{IKI},\mathsf{I\tilde{R}I}\}$, so the letter $\mathsf{I}$ cannot be followed by a subword of length greater than one. Therefore, the letter $\mathsf{I}$ can appear at most once in an element of $\langle{X}\rangle_\mathrm{irr}$ and, since subwords $\mathsf{EI}$ and $\mathsf{DI}$ have to be avoided, it can only be immediately preceded by the letters $\mathsf{\tilde{R}}$ or $\mathsf{\tilde{\Phi}}$. Altogether, the elements of $\langle{X}\rangle_\mathrm{irr}$ are precisely of the form
 \[
  \mathsf{\tilde{R}}U\mathsf{D}^j \quad\text{or}\quad \mathsf{\tilde{R}}V\mathsf{I\tilde{R}},
 \]
 where $j\in\mathbb{N}_0$ and each of $\mathsf{\tilde{R}}$ and $U\in\{\mathsf{E},\mathsf{\tilde{\Phi}},\mathsf{E\tilde{R}},\mathsf{\tilde{\Phi}\tilde{R}}\}$ and $V\in\{\mathsf{\tilde{\Phi}},\mathsf{E\tilde{R}},\mathsf{\tilde{\Phi}\tilde{R}}\}$ may be absent.
\end{proof}

\begin{remark}
 As discussed in Remark~\ref{rem:IDalgebras}, the induced evaluation of an integro-differential ring often is multiplicative in concrete examples. Also when considering several point evaluations of regular functions, the resulting functionals are multiplicative, i.e.\ $\varphi{fg}=(\varphi{f})\varphi{g}$ for all $f,g \in \mathcal{R}$. We discuss how reflecting this additional property of some functionals in the ring of operators influences the normal forms of operators.
 To this end, we consider a subset
 \[
  \Phi_m\subseteq\{\varphi \in \Phi\ |\ \varphi\text{ is multiplicative and }\varphi1=1\}
 \]
 of $\Phi$, which may or may not include the evaluation $\E$. For the corresponding elements $\varphi\in\Phi_m$ in $\mathcal{R}\langle{\Der,\Int,\Phi}\rangle$, we impose
 \begin{equation*}
  \varphi \cdot f = (\varphi{f})\varphi
 \end{equation*}
 for all $f \in \mathcal{R}$. To model this identity by reduction rules, we consider the submodule
 \begin{equation*}
  M_{\mathsf{\Phi_m}} := \mathcal{K}\Phi_m
 \end{equation*}
 of $M_\mathsf{\Phi}$.
 Evidently, we have the decomposition
 \[
  M_\mathsf{\Phi} = M_\mathsf{E} \oplus M_\mathsf{\tilde{\Phi}_m} \oplus M_{\tilde{\mathsf{\Phi}}}
 \]
 where the submodules $M_\mathsf{\tilde{\Phi}_m}$ and $M_{\tilde{\mathsf{\Phi}}}$ are generated by $\Phi_m\setminus\{\E\}$ and $\Phi\setminus(\{\E\}\cup\Phi_m)$, respectively.
We include the reduction rule
 \[
  (\mathsf{\Phi_mR},\varphi{\otimes}f\mapsto(\varphi f)\varphi)
 \]
 into Tables \ref{tab:SigmaIDOPhiDef} and \ref{tab:SigmaIDOPhi}, where $\varphi$ in the formula defining this $\mathcal{K}$-bimodule homomorphism on $M_\mathsf{\Phi_mR}$ is not a general element of $M_\mathsf{\Phi_m}$ but of $\Phi_m$ and the definition needs to be extended by left $\mathcal{K}$-linearity to all of $M_\mathsf{\Phi_mR}$.
 \par
 Theorem~\ref{thm:IDOPhi} generalizes to this situation with the additional restriction on normal forms in $\mathcal{K}\langle{M}\rangle$ that $h$ needs to be absent whenever $\varphi\in\Phi_m$. The proof is analogous by adapting the alphabets and the order accordingly. The additional reduction rule gives rise to additional ambiguities, whose resolvability is also checked in the Mathematica file. Determination of irreducible words also needs to be adapted accordingly.
 The resulting theorem includes Theorem~\ref{thm:IDOPhi} for $\Phi_m=\emptyset$ and it includes Theorem~27 from \cite{HosseinPoorRaabRegensburger2018} for $\Phi_m=\Phi$, i.e.\ when all functionals are multiplicative.
\qed
\end{remark}

\bibliographystyle{amsplain}
\bibliography{ref}

\providecommand{\bysame}{\leavevmode\hbox to3em{\hrulefill}\thinspace}
\providecommand{\MR}{\relax\ifhmode\unskip\space\fi MR }
\providecommand{\MRhref}[2]{%
  \href{http://www.ams.org/mathscinet-getitem?mr=#1}{#2}
}
\providecommand{\href}[2]{#2}
\begin{thebibliography}{10}

\bibitem{AschenbrennerDriesHoeven2017}
Matthias Aschenbrenner, Lou van~den Dries, and Joris van~der Hoeven,
  \emph{Asymptotic differential algebra and model theory of transseries},
  Annals of Mathematics Studies, vol. 195, Princeton University Press,
  Princeton, NJ, 2017.

\bibitem{Bavula2011}
Vladimir~V. Bavula, \emph{The algebra of integro-differential operators on a
  polynomial algebra}, J. Lond. Math. Soc. (2) \textbf{83} (2011), 517--543.

\bibitem{Bergman1978}
George~M. Bergman, \emph{The diamond lemma for ring theory}, Adv. in Math.
  \textbf{29} (1978), 178--218.

\bibitem{Buchberger1965}
Bruno Buchberger, \emph{An algorithm for finding the bases elements of the
  residue class ring modulo a zero dimensional polynomial ideal ({German})},
  Ph.D. thesis, University of Innsbruck, 1965.

\bibitem{CluzeauHosseinPoorQuadratRaabRegensburger2018}
Thomas Cluzeau, Jamal Hossein~Poor, Alban Quadrat, Clemens~G. Raab, and Georg
  Regensburger, \emph{Symbolic computation for integro-differential-time-delay
  operators with matrix coefficients}, IFAC-PapersOnLine \textbf{51} (2018),
  no.~14, 153--158.

\bibitem{CoddingtonLevinson1955}
Earl~A. Coddington and Norman Levinson, \emph{Theory of ordinary differential
  equations}, McGraw-Hill Book Co., Inc., New York-Toronto-London, 1955.

\bibitem{Cohn2003b}
Paul~M. Cohn, \emph{Further algebra and applications}, Springer-Verlag London,
  Ltd., London, 2003.

\bibitem{Edgar2010}
Gerald~A. Edgar, \emph{Transseries for beginners}, Real Anal. Exchange
  \textbf{35} (2010), no.~2, 253--309.

\bibitem{GaoGuoRosenkranz2018}
Xing Gao, Li~Guo, and Markus Rosenkranz, \emph{On rings of differential
  {R}ota-{B}axter operators}, Internat. J. Algebra Comput. \textbf{28} (2018),
  no.~1, 1--36.

\bibitem{Guo2012}
Li~Guo, \emph{An introduction to {R}ota-{B}axter algebra}, Surveys of Modern
  Mathematics, vol.~4, International Press, Somerville, MA; Higher Education
  Press, Beijing, 2012.

\bibitem{GuoKeigher2008}
Li~Guo and William Keigher, \emph{On differential {Rota-Baxter} algebras}, J.
  Pure Appl. Algebra \textbf{212} (2008), 522--540.

\bibitem{GuoRegensburgerRosenkranz2014}
Li~Guo, Georg Regensburger, and Markus Rosenkranz, \emph{On
  integro-differential algebras}, J. Pure Appl. Algebra \textbf{218} (2014),
  456--473.

\bibitem{HosseinPoor2018}
Jamal Hossein~Poor, \emph{Tensor reduction systems for rings of linear
  operators}, Ph.D. thesis, Johannes Kepler University Linz, Austria, 2018.

\bibitem{HosseinPoorRaabRegensburger2018}
Jamal Hossein~Poor, Clemens~G. Raab, and Georg Regensburger, \emph{Algorithmic
  operator algebras via normal forms in tensor rings}, J. Symbolic Comput.
  \textbf{85} (2018), 247--274.

\bibitem{Kaplansky1957}
Irving Kaplansky, \emph{An introduction to differential algebra}, Publ. Inst.
  Math. Univ. Nancago, No. 5, Hermann, Paris, 1957.

\bibitem{Keigher1975}
William~F. Keigher, \emph{Adjunctions and comonads in differential algebra},
  Pacific J. Math. \textbf{59} (1975), no.~1, 99--112.

\bibitem{Keigher1997}
\bysame, \emph{On the ring of {H}urwitz series}, Comm. Algebra \textbf{25}
  (1997), no.~6, 1845--1859.

\bibitem{KeigherPritchard2000}
William~F. Keigher and F.~Leon Pritchard, \emph{Hurwitz series as formal
  functions}, J. Pure Appl. Algebra \textbf{146} (2000), no.~3, 291--304.

\bibitem{KnuthBendix1970}
Donald~E. Knuth and Peter~B. Bendix, \emph{Simple word problems in universal
  algebras}, Computational {P}roblems in {A}bstract {A}lgebra, Pergamon,
  Oxford, 1970, pp.~263--297.

\bibitem{Kummer1840}
Ernst~E. Kummer, \emph{Ueber die {T}ranscendenten, welche aus wiederholten
  {I}ntegrationen rationaler {F}ormeln entstehen}, J. Reine Angew. Math.
  \textbf{21} (1840), 74--90, 193--225, and 328--371.

\bibitem{LappoDanilevski1928}
Ivan~A. Lappo-Danilevski, \emph{R\'{e}solution algorithmique des probl\`{e}mes
  r\'{e}guliers de {P}oincar\'{e} et de riemann ({M}\'{e}moire premier: {L}e
  probl\`{e}me de {P}oincar\'{e}, concernant de la construction d'un groupe de
  monodromie d'un syst\`{e}me donn\'{e} d'\'{e}quations diff\'{e}rentielles
  lin\'{e}aires aux int\'{e}grales r\'{e}guli\`{e}res)}, J. Soc. Phys.-Math.
  L\'{e}ningrade \textbf{2} (1928), no.~1, 94--120.

\bibitem{Mora1994}
Teo Mora, \emph{An introduction to commutative and noncommutative {G}r\"obner
  bases}, Theoret. Comput. Sci. \textbf{134} (1994), no.~1, 131--173.

\bibitem{Przeworska-Rolewicz1988}
Danuta Przeworska-Rolewicz, \emph{Algebraic analysis}, PWN---Polish Scientific
  Publishers, Warsaw; D. Reidel Publishing Co., Dordrecht, 1988.

\bibitem{Przeworska-Rolewicz2000}
\bysame, \emph{Two centuries of the term ``algebraic analysis''}, Algebraic
  analysis and related topics ({W}arsaw, 1999), Banach Center Publ., vol.~53,
  Polish Acad. Sci. Inst. Math., Warsaw, 2000, pp.~47--70.

\bibitem{Quadrat2015}
Alban Quadrat, \emph{A constructive algebraic analysis approach to {A}rtstein's
  reduction of linear time-delay systems}, IFAC-PapersOnLine \textbf{48}
  (2015), no.~12, 209--214.

\bibitem{QuadratRegensburger2020}
Alban Quadrat and Georg Regensburger, \emph{Computing polynomial solutions and
  annihilators of integro-differential operators with polynomial coefficients},
  Algebraic and symbolic computation methods in dynamical systems, Adv. Delays
  Dyn., vol.~9, Springer, Cham, 2020, pp.~87--114.

\bibitem{RaabRegensburger}
Clemens~G. Raab and Georg Regensburger, \emph{The free commutative
  integro-differential ring}, in preparation.

\bibitem{Ree1958}
Rimhak Ree, \emph{Lie elements and an algebra associated with shuffles}, Ann.
  of Math. (2) \textbf{68} (1958), 210--220.

\bibitem{RegensburgerRosenkranzMiddeke2009}
Georg Regensburger, Markus Rosenkranz, and Johannes Middeke, \emph{A skew
  polynomial approach to integro-differential operators}, Proceedings of ISSAC
  '09 (New York, NY, USA), ACM, 2009, pp.~287--294.

\bibitem{Rosenkranz2005}
Markus Rosenkranz, \emph{A new symbolic method for solving linear two-point
  boundary value problems on the level of operators}, J. Symbolic Comput.
  \textbf{39} (2005), 171--199.

\bibitem{RosenkranzRegensburger2008a}
Markus Rosenkranz and Georg Regensburger, \emph{Solving and factoring boundary
  problems for linear ordinary differential equations in differential
  algebras}, J. Symbolic Comput. \textbf{43} (2008), 515--544.

\bibitem{RosenkranzRegensburgerTecBuchberger2012}
Markus Rosenkranz, Georg Regensburger, Loredana Tec, and Bruno Buchberger,
  \emph{Symbolic analysis for boundary problems: {F}rom rewriting to
  parametrized {G}r{\"o}bner bases}, {Numerical and Symbolic Scientific
  Computing: Progress and Prospects} (Ulrich Langer and Peter Paule, eds.),
  Springer Vienna, Vienna, 2012, pp.~273--331.

\bibitem{RosenkranzSerwa2019}
Markus Rosenkranz and Nitin Serwa, \emph{An integro-differential structure for
  {D}irac distributions}, J. Symbolic Comput. \textbf{92} (2019), 156--189.

\bibitem{Rota2001}
Gian-Carlo Rota, \emph{Twelve problems in probability no one likes to bring
  up}, Algebraic combinatorics and computer science, Springer Italia, Milan,
  2001, pp.~57--93.

\bibitem{Rowen1991}
Louis~H. Rowen, \emph{Ring theory}, student ed., Academic Press, Inc., Boston,
  MA, 1991.

\bibitem{Stanley1980}
Richard~P. Stanley, \emph{Differentiably finite power series}, European J.
  Combin. \textbf{1} (1980), no.~2, 175--188.

\bibitem{Hoeven2006}
Joris van~der Hoeven, \emph{Transseries and real differential algebra}, Lecture
  Notes in Mathematics, vol. 1888, Springer-Verlag, Berlin, 2006.

\end{thebibliography}

\end{document}